\documentclass[12pt]{article}
\usepackage[utf8]{inputenc}
\usepackage[T1]{fontenc}
\usepackage[french,english]{babel}
\usepackage{csquotes}
\usepackage{lmodern}
\usepackage[a4paper, top=3cm, bottom=3cm, left=3cm, right=3cm]{geometry}
\usepackage{setspace}
\usepackage{emptypage}
\usepackage{microtype}
\UseMicrotypeSet[protrusion]{basicmath}
\usepackage[defaultlines=3,all]{nowidow}
\usepackage{needspace}
\usepackage{amsmath,amsthm,amssymb,amsfonts}
\usepackage{mathtools}
\usepackage{mathdots}
\usepackage{stmaryrd}
\usepackage{cancel}
\allowdisplaybreaks[2]
\usepackage[all,pdf]{xy}
\usepackage{tikz-cd}
\usepackage{quiver}
\usepackage{xcolor}
\usepackage[backend=biber,style=alphabetic,url=true,nohashothers=false,maxnames=50,safeinputenc,bibencoding=utf8]{biblatex}
\DeclareNameAlias{sortname}{family-given}
\DeclareNameAlias{default}{family-given}
\usepackage{pdfpages}
\usepackage{clipboard}
\usepackage[colorlinks]{hyperref}
\hypersetup{colorlinks=false,pdfborder={0 0 1},urlcolor=blue,linkcolor=blue,citecolor=blue}
\usepackage{bookmark}
\usepackage{cleveref}
\newcommand{\bb}{\mathbb}
\newcommand{\mcal}{\mathcal}

\newcommand{\op}{\mathrm{op}}

\DeclareMathOperator{\Spec}{Spec}
\DeclareMathOperator{\Sym}{Sym}

\DeclareMathOperator{\QCoh}{QCoh}

\DeclareMathOperator{\Perf}{Perf}
\DeclareMathOperator{\dSt}{dSt}

\DeclareMathOperator{\CAlg}{CAlg}
\DeclareMathOperator{\Mod}{Mod}
\DeclareMathOperator{\dAff}{dAff}
\DeclareMathOperator{\dRings}{dRings}
\DeclareMathOperator{\PerfLinSt}{PerfLinSt}
\DeclareMathOperator{\PerfLieSt}{PerfLieSt}

\DeclareMathOperator{\cofib}{cofib}

\newcommand{\DR}{\mathbf{DR}}
\newcommand{\Map}{\mathbf{Map}}
\newcommand{\Arr}{\mathrm{Arr}}

\newcommand{\Aff}{\dAff}

\newcommand{\printbib}{\begin{raggedright}\printbibliography[heading=bibintoc,title={Bibliography}]\end{raggedright}}
\newcommand{\restatable}[3]{\begin{#1}[\Cref{#2-restate}]\label[#1]{#2}\Copy{#2}{#3}\end{#1}}
\newcommand{\restate}[2]{\begin{#1}[\Cref{#2}]\label[#1]{#2-restate} \Paste{#2} \end{#1}}

\newcommand{\thmneedspace}{\Needspace{2\baselineskip}}
\newcommand{\separation}{\par\medskip}
\newcommand{\thmheadingname}{Theorem}
\newcommand{\dfheadingname}{Definition}
\newcommand{\propheadingname}{Proposition}
\newcommand{\corheadingname}{Corollary}
\newcommand{\lemheadingname}{Lemma}
\newcommand{\slemheadingname}{Sub-lemma}
\newcommand{\expectedheadingname}{Expected result}
\newcommand{\exheadingname}{Example}
\newcommand{\rmkheadingname}{Remark}
\newcommand{\constructionheadingname}{Construction}
\newcommand{\notationheadingname}{Notation}
\newcommand{\reminderheadingname}{Reminder}
\newcommand{\conclusionheadingname}{Conclusion}
\newcommand{\setenglishthmnames}{%
	\renewcommand{\thmheadingname}{Theorem}%
	\renewcommand{\dfheadingname}{Definition}%
	\renewcommand{\propheadingname}{Proposition}%
	\renewcommand{\corheadingname}{Corollary}%
	\renewcommand{\lemheadingname}{Lemma}%
	\renewcommand{\slemheadingname}{Sub-lemma}%
	\renewcommand{\expectedheadingname}{Expected result}%
	\renewcommand{\exheadingname}{Example}%
	\renewcommand{\rmkheadingname}{Remark}%
	\renewcommand{\constructionheadingname}{Construction}%
	\renewcommand{\notationheadingname}{Notation}%
	\renewcommand{\reminderheadingname}{Reminder}%
	\renewcommand{\conclusionheadingname}{Conclusion}%
	\crefname{section}{section}{sections}%
	\Crefname{section}{Chapter}{Chapters}%
	\crefname{appendix}{appendix}{appendices}%
	\Crefname{appendix}{Appendix}{Appendices}%
	\crefname{thm}{theorem}{theorems}%
	\Crefname{thm}{Theorem}{Theorems}%
	\crefname{df}{definition}{definitions}%
	\Crefname{df}{Definition}{Definitions}%
	\crefname{prop}{proposition}{propositions}%
	\Crefname{prop}{Proposition}{Propositions}%
	\crefname{cor}{corollary}{corollaries}%
	\Crefname{cor}{Corollary}{Corollaries}%
	\crefname{lem}{lemma}{lemmas}%
	\Crefname{lem}{Lemma}{Lemmas}%
	\crefname{slem}{sub-lemma}{sub-lemmas}%
	\Crefname{slem}{Sub-lemma}{Sub-lemmas}%
	\crefname{expected}{expected result}{expected results}%
	\Crefname{expected}{Expected result}{Expected results}%
	\crefname{ex}{example}{examples}%
	\Crefname{ex}{Example}{Examples}%
	\crefname{rmk}{remark}{remarks}%
	\Crefname{rmk}{Remark}{Remarks}%
	\crefname{construction}{construction}{constructions}%
	\Crefname{construction}{Construction}{Constructions}%
	\crefname{notation}{notation}{notations}%
	\Crefname{notation}{Notation}{Notations}%
	\crefname{reminder}{reminder}{reminders}%
	\Crefname{reminder}{Reminder}{Reminders}%
	\crefname{conclusion}{conclusion}{conclusions}%
	\Crefname{conclusion}{Conclusion}{Conclusions}%
}
\newcommand{\setfrenchthmnames}{%
	\renewcommand{\thmheadingname}{Théorème}%
	\renewcommand{\dfheadingname}{Définition}%
	\renewcommand{\propheadingname}{Proposition}%
	\renewcommand{\corheadingname}{Corollaire}%
	\renewcommand{\lemheadingname}{Lemme}%
	\renewcommand{\slemheadingname}{Sous-lemme}%
	\renewcommand{\expectedheadingname}{Résultat attendu}%
	\renewcommand{\exheadingname}{Exemple}%
	\renewcommand{\rmkheadingname}{Remarque}%
	\renewcommand{\constructionheadingname}{Construction}%
	\renewcommand{\notationheadingname}{Notation}%
	\renewcommand{\reminderheadingname}{Rappel}%
	\renewcommand{\conclusionheadingname}{Conclusion}%
	\crefname{section}{chapitre}{chapitres}%
	\Crefname{section}{Chapitre}{Chapitres}%
	\crefname{appendix}{annexe}{annexes}%
	\Crefname{appendix}{Annexe}{Annexes}%
	\crefname{thm}{théorème}{théorèmes}%
	\Crefname{thm}{Théorème}{Théorèmes}%
	\crefname{df}{définition}{définitions}%
	\Crefname{df}{Définition}{Définitions}%
	\crefname{prop}{proposition}{propositions}%
	\Crefname{prop}{Proposition}{Propositions}%
	\crefname{cor}{corollaire}{corollaires}%
	\Crefname{cor}{Corollaire}{Corollaires}%
	\crefname{lem}{lemme}{lemmes}%
	\Crefname{lem}{Lemme}{Lemmes}%
	\crefname{slem}{sous-lemme}{sous-lemmes}%
	\Crefname{slem}{Sous-lemme}{Sous-lemmes}%
	\crefname{expected}{résultat attendu}{résultats attendus}%
	\Crefname{expected}{Résultat attendu}{Résultats attendus}%
	\crefname{ex}{exemple}{exemples}%
	\Crefname{ex}{Exemple}{Exemples}%
	\crefname{rmk}{remarque}{remarques}%
	\Crefname{rmk}{Remarque}{Remarques}%
	\crefname{construction}{construction}{constructions}%
	\Crefname{construction}{Construction}{Constructions}%
	\crefname{notation}{notation}{notations}%
	\Crefname{notation}{Notation}{Notations}%
	\crefname{reminder}{rappel}{rappels}%
	\Crefname{reminder}{Rappel}{Rappels}%
	\crefname{conclusion}{conclusion}{conclusions}%
	\Crefname{conclusion}{Conclusion}{Conclusions}%
}
\theoremstyle{plain}
\newtheorem{thm}{\thmheadingname}[subsection]
\newtheorem{df}[thm]{\dfheadingname}
\newtheorem{prop}[thm]{\propheadingname}
\newtheorem{cor}[thm]{\corheadingname}
\newtheorem{lem}[thm]{\lemheadingname}

\newtheorem{expected}[thm]{\expectedheadingname}
\theoremstyle{definition}
\newtheorem{ex}[thm]{\exheadingname}
\newtheorem{rmk}[thm]{\rmkheadingname}
\newtheorem{construction}[thm]{\constructionheadingname}
\newtheorem{notation}[thm]{\notationheadingname}
\newtheorem{reminder}[thm]{\reminderheadingname}
\newtheorem{conclusion}[thm]{\conclusionheadingname}
\AtBeginDocument{\setenglishthmnames}
\AddToHook{env/thm/begin}{\thmneedspace}
\AddToHook{env/df/begin}{\thmneedspace}
\AddToHook{env/prop/begin}{\thmneedspace}
\AddToHook{env/cor/begin}{\thmneedspace}
\AddToHook{env/lem/begin}{\thmneedspace}
\AddToHook{env/slem/begin}{\thmneedspace}
\AddToHook{env/expected/begin}{\thmneedspace}
\AddToHook{env/ex/begin}{\thmneedspace}
\AddToHook{env/rmk/begin}{\thmneedspace}
\AddToHook{env/construction/begin}{\thmneedspace}
\AddToHook{env/notation/begin}{\thmneedspace}
\AddToHook{env/reminder/begin}{\thmneedspace}
\AddToHook{env/conclusion/begin}{\thmneedspace}
\title{A Mapping Theorem for Derived Foliations}
\author{Victor Alfieri}
\date{}
\begin{document}
	\maketitle
	
	\begin{abstract}
		In this paper, we construct in characteristic zero a derived foliation on derived mapping stacks $\underline{\Map}_S(X,Y)$ under natural representability assumptions, for $S$ a base derived stack, $X$ a proper schematic, flat, and local complete intersection derived stack over $S$, and $Y$ a relative derived Deligne-Mumford stack locally of finite presentation over $S$, when $Y$ is equipped with a derived foliation relative to $S$.
		
		In the process, given a derived stack $Z$ locally of finite presentation over a derived stack $X$, we will first show that the $\infty$-category of derived foliations over $Z$ relative to $X$ embeds as a full subcategory of derived stacks over $Z$ equipped with extra structure, and describe its essential image explicitly.
		
		We will then show that given a proper schematic, flat, and local complete intersection map of derived stacks $f : X \to Y$, the push-forward or Weil-restriction functor $f_*$ from derived stacks over $X$ to derived stacks over $Y$ preserves the preceding essential images, and thus defines a push-forward, from derived foliations over $Z$ relative to $X$, to derived foliations over $f_* Z$ relative to $Y$, under natural representability assumptions. The aforementioned result on derived mapping stacks is obtained as a special case of this statement.
		
		As example applications, given a smooth projective scheme $X$ equipped with a derived foliation, we obtain derived foliations on the derived moduli stacks $\bb R \overline{\mathbf M}_{g,n}(X)$ and $\bb R \mathbf{Hilb}^{lci}(X)$, which are respectively the derived enhancements of the moduli stack of families of stable curves of genus $g$ with $n$ marked points on $X$, and of the Hilbert scheme of closed subschemes of $X$.
	\end{abstract}
	\tableofcontents
	\section{Introduction} \label[section]{Intro-en}
	In this paper, we will work in derived algebraic geometry in characteristic zero\footnote{We expect results from this paper to extend to positive characteristics. However, some of them will require reformulations, particularly in the distinction between derived foliations and infinitesimal foliations which only coincide in characteristic zero.}, and prove existence theorems for derived foliations. Derived algebraic geometry was introduced by Bertrand Toën and Gabriele Vezzosi in \cite{toenVezzosiHAGI,toenVezzosiHAGII}, and by Jacob Lurie in \cite{lurieDAG}. Foliations in this context were introduced and studied in a series of papers by Bertrand Toën and Gabriele Vezzosi, \cite{AGI} and \cite{AGII}, and later on in their book \cite{toen2025derivedfoliations}.
	\subsection{Classical Foliations}
	In classical differential geometry, a foliation on a differential manifold $M$ is morally the data of a partition of $M$ into sub-manifolds, called the "leaves" of the foliation. Given such a foliation $F$ and a point $x \in M$, one can define the sub-vector space $T_x F \subset T_x M$ of tangent vectors to $x$ in $F_x \subset M$, where $F_x$ is the only leaf of $F$ containing $x$. Collectively they form a sub-bundle $TF \subset TM$ of the tangent bundle of $M$, called the tangent bundle of the foliation, which is moreover stable by the Lie bracket. By dualizing, one obtains instead a quotient of the de Rham algebra $\Omega^\bullet(M) \to \Omega^\bullet(F)$ which inherits an internal differential from the de Rham differential, and $\Omega^1(F)$ is called the cotangent complex of the foliation. Moreover, when properly defined and in the smooth case, all these descriptions of foliations turn out to be equivalent to each other. A reference on the subject covering these different definitions and their equivalence is \cite{classicalfoliations}.
	
	In classical algebraic geometry, the latter point of view has been used to define algebraic foliations, which are defined as differential ideals of the cotangent complex. This allows in particular to consider foliations on singular varieties which was already not possible in the differential case, but they are difficult to construct and cannot often be transported along morphisms. A reference on modern applications of the theory of foliations in algebraic geometry is \cite{AGFoliations}.
	\subsection{Derived Foliations and Main Results}
	In derived geometry, we consider instead a notion of derived foliation which is less rigid and can thus exist in a wider range of contexts, as studied in \cite{toen2025derivedfoliations}.
	Most notably, in this paper we will prove the following main theorem.
	\restatable{thm}{main_theorem}{
		Let $S$ be a base derived stack in characteristic 0. Let $X$ be a proper-schematic, flat, and local complete intersection derived stack over $S$, and $Y$ be a relative derived Deligne-Mumford stack locally of finite presentation over $S$, equipped with a derived foliation $\mcal F$ relative to $S$.
		Then, assuming the derived mapping stack $\underline{\Map}_S(X,Y)$ is a relative derived Artin stack locally of finite presentation over $S$, it inherits a derived foliation $\mcal F_{\Map}$ relative to $S$.
	}
	\begin{rmk}
		Note that in comparison, in the differential case, the statement of \cref{main_theorem} would already require clarifications on its meaning because $\Map(X,Y)$ is an infinite dimensional object. Moreover, even in the case of classical algebraic geometry where the meaning of such a statement is clear, without the tools of derived geometry its conclusion would be false in this generality and require much stronger hypotheses, because this object would be very singular and its de Rham algebra extremely difficult to compute and not bounded.
	\end{rmk}
	
	Furthermore, as in classical foliation theory, derived foliations admit tangent complexes, and the tangent complex of the induced derived foliation on $\underline{\Map}_S(X,Y)$ stated by the above theorem has the following description:
	\restatable{prop}{main_prop}{
		In the context of \cref{main_theorem} (\ref{main_theorem-restate}), let $\bb T_{\mcal F}$ (resp.\ $\bb T_{\mcal F,\Map})$ denote the tangent complex of $\mcal F$ (resp.\ $\mcal F_{\Map}$), and $\bb T_{X/S}$ the tangent complex of $X$ relative to $S$. For any $S$-point of $\underline{\Map}_S(X,Y)$, which corresponds to an $S$-linear map $g : X \to Y$, there is a canonical equivalence:
		\[\bb T_g \mcal F_{\Map} \simeq \Gamma_S(X,g^* \bb T_{\mcal F})\]
		where $\Gamma_S : \QCoh(X) \to \QCoh(S)$ is the global-sections functor relative to the base $S$, i.e.\ the push-forward of quasi-coherent sheaves along the structural map $X \to S$, which recovers the usual global sections when $S = \Spec(k)$.
		More generally, the tangent complex of $\mcal F_{\Map}$ itself is given by:
		\[\bb T_{\mcal F,\Map} \simeq \pi_* ev^* \bb T_{\mcal F}\]
		where $ev : X \times_S \underline{\Map}_S(X,Y) \to Y$ is the evaluation, and $\pi : X \times_S \underline{\Map}_S(X,Y) \to \underline{\Map}_S(X,Y)$ is the projection
	}
	Intuitively, this property means that a leaf $L$ of $\mcal F_{\Map}$ is morally the data of a family of maps $(g_l : X \to Y)_{l \in L}$ such that, for any $x \in X$, the family $(g_l(x))_{l \in L}$ is a leaf of $\mcal F$.
	\separation
	Let us now give some important example applications of this theorem, which are also proved in this paper.
	\begin{ex}[See proof in \cref{foliation_Mgn}] \label[ex]{example_application_1}
		Let $\bb R \mcal C_{g,n}^{pre}$ be the derived enhancement of the universal family of curves of genus $g$ with $n$ marked points, defined over the derived enhancement of the moduli stack of pre-stable curves of genus $g$ with $n$ marked points $\bb R \mathbf M_{g,n}^{pre}$, and $X$ be a smooth projective scheme equipped with a derived foliation.
		In this case, $\bb R \mathbf M_{g,n}^{pre}(X)$, defined as $\underline{\Map}_{\bb R \mathbf M_{g,n}^{pre}}(\bb R \mcal C_{g,n}^{pre},X \times \bb R \mathbf M_{g,n}^{pre})$, inherits a derived foliation, and moreover it contains $\bb R \overline{\mathbf M}_{g,n}(X;\beta)$ as an open substack, which is the derived enhancement of the moduli stack of families of stable curves of genus $g$ with $n$ marked points of class $\beta$ on $X$, which also inherits a derived foliation by restriction. This provides a more systematic construction of the derived foliation introduced in \cite[Section 2.2]{toen2022foliations}, in which only the case of derived foliations arising from vector fields was treated, whereas the results of this paper extend to arbitrary derived foliations even when they do not arise from vector fields or maps.
		In particular, one could consider a new kind of Gromov-Witten invariants, relative to this derived foliation. Note that in the statement of this example, we restricted ourselves to $X$ being a smooth projective variety precisely to match the usual setting of Gromov-Witten theory, however our construction would yield a foliation for more general $X$.
	\end{ex}
	\begin{ex}[See proof in \cref{foliation_Hilb}] \label[ex]{example_application_2}
		Similarly, let $Y$ be a derived Deligne-Mumford stack locally of finite presentation equipped with a derived foliation. We can form the derived stack $\bb R \mathbf{Hilb}^{lci}(Y)$, whose $S$-points for $S$ a derived scheme are closed immersions $i : Z \hookrightarrow Y$ where $Z$ is a derived scheme over $S$ and $Z \to S$ is proper and local complete intersection (lci), thus $\bb R \mathbf{Hilb}^{lci}(Y)$ is canonically embedded in a mapping stack, and inherits a derived foliation from $Y$. Moreover, we can study in particular the integral points of this derived foliation, i.e. the locus where the foliation is $0$. See \cite[Def.~1.9.2(1)]{EDSbook} for the original definition in differential geometry. In our setting, we expect these integral points to correspond to closed derived subschemes $Z \subset Y$ that are stable by the derived foliation on $Y$, i.e.\ closed derived subschemes of algebraic leaves of the derived foliation, but this statement is neither formally stated nor proved in this paper.
	\end{ex}
	For more information on both of these examples, see \cref{applications}. Similar results could be obtained for various derived moduli stacks, provided the universal object is proper-schematic, flat, and local complete intersection over the derived moduli stack in question.
	\separation
	In order to prove \cref{main_theorem}, we will more generally construct a push-forward of derived foliations, as stated next. Note that for a map of derived stacks $Z \to X$, we will denote by $\mcal Fol(Z/X)$ the $\infty$-category of derived foliations defined on $Z$ relative to $X$, and for a map of derived stacks $f : X \to Y$, we will denote by $f_* : \dSt/X \to \dSt/Y$ the push-forward functor, also called direct image or Weil restriction.
	\restatable{thm}{thm_pushforward_fol}{
		Let $Z \to X$ be a map of derived stacks which is a relative derived Deligne-Mumford stack locally of finite presentation, and $f : X \to Y$ a proper-schematic, flat, and local complete intersection map of derived stacks. Assume that $f_* Z \to Y$ is a relative derived Artin stack locally of finite presentation.
		Then, there is a push-forward functor:
		\[f_{*,\mcal Fol} : \mcal Fol(Z/X) \to \mcal Fol(f_* Z/Y)\]
	}
	In order to construct this push-forward, we will show that $\mcal Fol(Z/X)$ is equivalent to an $\infty$-category of derived stacks over $Z$ equipped with extra structure (see \cref{equiv_fol_dst_rel_dm}), and that the canonical push-forward functor of derived stacks $f_* : \dSt/X \to \dSt/Y$ preserves this structure (see \cref{pushforward_preserves_fols}).
	Moreover, the tangent complex of the push-forward of a foliation admits a description very similar to that of \cref{main_prop}, see \cref{tangent_pushforward_fol}.
	Finally, note that the following special case of the above theorem yields the main \cref{main_theorem}.
	\begin{ex} \label[ex]{ex_main_thm}
		Let $S$ be a derived stack, $X$ be a proper-schematic, flat, and local complete intersection derived stack over $S$, and $F$ be a relative derived Deligne-Mumford stack locally of finite presentation over $S$ equipped with a derived foliation relative to $S$.
		Let $Y = S$ and $Z = X \times_S F$, and note that $X \to Y$ is thus proper-schematic, flat, and local complete intersection, and that the projection $Z \to X$ is a relative derived Deligne-Mumford stack locally of finite presentation.
		Moreover, there is a canonical equivalence $f_* Z = \underline{\Map}_S(X,F)$ shown in the proof of \cref{main_theorem-restate}. Assuming $\underline{\Map}_S(X,F)$ is a relative derived Artin stack over $S$, the latter inherits a derived foliation by application of the above theorem.
	\end{ex}
	\subsection{Plan of the Paper}
	To finish this introductory section, we will now outline the plan of this paper.
	
	In \cref{Defs}, we will first recall the algebraic definitions of derived foliations, then rephrase them using geometric objects, by reviewing results from previous papers. The conclusions of this section are summarized in \cref{defs_conclusion}.
	
	In \cref{Foliations_Geometric}, this geometric description will be used to construct a functor from the $\infty$-category of derived foliations to an $\infty$-category of derived stacks equipped with extra structure, as summarized in \cref{main_diagram}. We will then show that this functor is fully faithful, and explicitly describe its essential image, which yields an equivalence of $\infty$-categories between derived foliations and this essential image, as stated in \cref{equiv_fol_dst} for the affine case, and generalized to relative derived Deligne-Mumford stacks in \cref{equiv_fol_dst_rel_dm}. Note that this paper contains an \cref{beck-chevalley-section}, which consists in a collection of results on Beck-Chevalley conditions, precisely because these conditions appear in the proofs of the \cref{equiv-affine-case} of this section.
	
	In \cref{Push-forward}, we will show that the push-forward functor from derived stacks over the base to derived stacks over the target preserves the extra structure of the derived stacks mentioned above, as stated in \cref{pushforward_preserves_fols}. This will yield both the above \cref{thm_pushforward_fol} and finally the main \cref{main_theorem}, as well as the description of the tangent complex of the resulting derived foliations in \cref{main_prop}.
	
	In \cref{applications}, we will prove two applications of the results of this paper, mentioned in \cref{example_application_1} and \cref{example_application_2}, respectively for Gromov-Witten theory and Hilbert schemes.
	
	Finally, in \cref{future_work}, we will talk about ongoing and future work following this paper.
	\separation
	\textbf{Acknowledgements:} I would like to give many thanks to Bertrand Toën and Marco Robalo for their precious help when writing this paper. I also want to thank Damien Calaque and Pavel Safronov for having subsequently reviewed its content and provided many suggestions for improvement as examiners of my thesis. In my thesis, more specifically \cite[Chapters~4-6, Appendix~A]{Alfieri2026ShiftedSymplecticFoliations}, I included and revised content from a previous version of this paper. This version is a further revision of this content in my thesis.
	\separation
	\textbf{Use of generative AI:} no generative AI tool was used to produce content in this paper.
	\section{Derived Foliations} \label[section]{Defs}
	\subsection{Introduction}
	Throughout this paper, we will work in characteristic zero, i.e. over a fixed $\bb Q$-algebra $k$. Since $k$ is fixed, it will be omitted from all notations.
	
	\bigskip
	
	In this section we will first recall the algebraic definitions of graded mixed algebras, some constructions on them, the properties used to define derived foliations, and discuss their stability under the preceding constructions.
	
	Secondly, we will recall results from \cite{Mou} and \cite{HKR} to reformulate the algebraic definitions in a more geometric context, i.e. in terms of quasi-coherent sheaves over some derived stacks.
	\subsection{Algebraic Definitions} \label{alg_defs}
	In this subsection, we will recall the definitions of derived foliations in derived algebraic geometry. In a first part, we will recall the algebraic definitions of graded mixed algebras, some constructions on them, the properties used to define derived foliations, and discuss their stability under the preceding constructions. In a second part, we will recall results from \cite{Mou} and \cite{HKR} to reformulate these algebraic definitions in a more geometric context, i.e.\ in terms of quasi-coherent sheaves over some derived stacks.
	We will first recall the definition of the structure of a graded mixed algebra over the base field $k$, which we will call the absolute definition. Graded mixed algebras are essentially graded algebras equipped with internal differentials shifting the grading, whose main class of examples is de Rham algebras. Secondly, we will recall the definitions of graded and graded mixed objects defined over other $k$-algebras, which we will call the relative definitions.
	
	All the definitions in this subsection can be found in \cite{PTVV} and \cite{toen2025derivedfoliations}.
	\subsubsection{Absolute Definitions}
	\begin{notation}
		Let us list here the notations we use related to graded objects:
		\begin{itemize}
			\item The 1-category of $\bb Z$-graded $k$-dg modules, which are families of $k$-dg-modules $(E(n))_{n \in \bb Z}$, is defined as:
			\[dg^{gr}_{st} := \textup{Fun}(\bb Z^{dis},dg_k)\]
			As in \cite{PTVV}, we use the cohomological conventions, so its objects and all complexes in this paper are to be considered as cochain complexes.
			\item We will refer to its objects as graded modules, and the extra grading will be called the weight grading.
			\item The 1-category $dg^{gr}_{st}$ is canonically equipped with a class of weak equivalences comprised of graded maps that induce quasi-isomorphism of $k$-dg-modules in each weight.
			\item The 1-category $dg^{gr}_{st}$ also has a canonical symmetric monoidal structure, given by a Cauchy product formula, and the 1-category of commutative algebras is denoted by:
			\[cdga^{gr}_{st} := CAlg(dg^{gr}_{st})\]
			\item The 1-category $cdga^{gr}_{st}$ inherits weak equivalences from $dg^{gr}_{st}$.
		\end{itemize}
	\end{notation}
	\begin{df}
		Given a graded module $E$, a mixed structure on $E$ consists in a family of maps of dg-modules $\epsilon_n : E(n) \to E(n+1)[-1]$ (where $[-1]$ denotes shifting the cohomological grading of the complex by -1), in particular satisfying $d \circ \epsilon + \epsilon \circ d = 0$, and moreover satisfying $\epsilon^2 = 0$.
	\end{df}
	This data can be visualized in a diagram of the following shape, in which each column is a given weight, downwards arrows are cohomological differentials within each weight, and diagonal arrows are the mixed differentials:
	
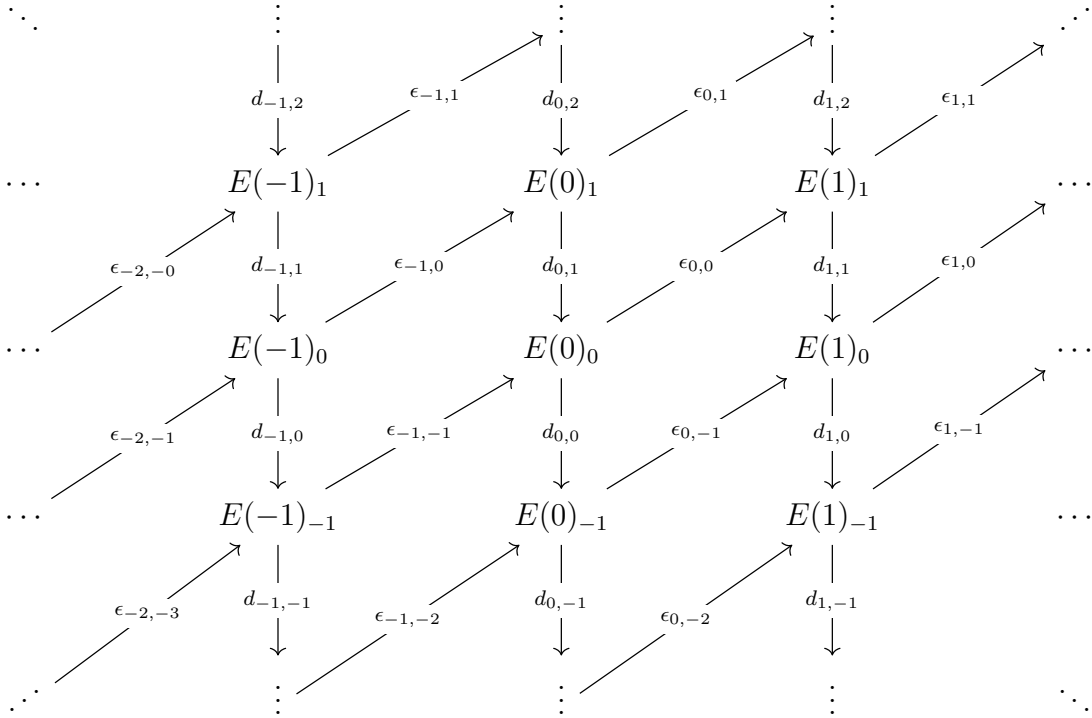
\[\begin{tikzcd}
		\ddots && \vdots && \vdots && \vdots && \iddots \\
		\\
		\cdots && {E(-1)_{-1}} && {E(0)_{-1}} && {E(1)_{-1}} && \cdots \\
		\\
		\cdots && {E(-1)_0} && {E(0)_0} && {E(1)_0} && \cdots \\
		\\
		\cdots && {E(-1)_{1}} && {E(0)_{1}} && {E(1)_{1}} && \cdots \\
		\\
		\iddots && \vdots && \vdots && \vdots && \ddots
		\arrow["{\epsilon_{-2,-2}}"{description}, from=1-1, to=3-3]
		\arrow["{d_{-1,-2}}"{description}, from=1-3, to=3-3]
		\arrow["{\epsilon_{-1,-2}}"{description}, from=1-3, to=3-5]
		\arrow["{d_{0,-2}}"{description}, from=1-5, to=3-5]
		\arrow["{\epsilon_{0,-2}}"{description}, from=1-5, to=3-7]
		\arrow["{d_{1,-2}}"{description}, from=1-7, to=3-7]
		\arrow["{\epsilon_{-2,-1}}"{description}, from=3-1, to=5-3]
		\arrow["{d_{-1,-1}}"{description}, from=3-3, to=5-3]
		\arrow["{\epsilon_{-1,-1}}"{description}, from=3-3, to=5-5]
		\arrow["{d_{0,-1}}"{description}, from=3-5, to=5-5]
		\arrow["{\epsilon_{0,-1}}"{description}, from=3-5, to=5-7]
		\arrow["{d_{1,-1}}"{description}, from=3-7, to=5-7]
		\arrow["{\epsilon_{1,-1}}"{description}, from=3-7, to=5-9]
		\arrow["{\epsilon_{-2,0}}"{description}, from=5-1, to=7-3]
		\arrow["{d_{-1,0}}"{description}, from=5-3, to=7-3]
		\arrow["{\epsilon_{-1,0}}"{description}, from=5-3, to=7-5]
		\arrow["{d_{0,0}}"{description}, from=5-5, to=7-5]
		\arrow["{\epsilon_{0,0}}"{description}, from=5-5, to=7-7]
		\arrow["{d_{1,0}}"{description}, from=5-7, to=7-7]
		\arrow["{\epsilon_{1,0}}"{description}, from=5-7, to=7-9]
		\arrow["{d_{-1,1}}"{description}, from=7-3, to=9-3]
		\arrow["{\epsilon_{-1,-1}}"{description}, from=7-3, to=9-5]
		\arrow["{d_{0,1}}"{description}, from=7-5, to=9-5]
		\arrow["{\epsilon_{0,-1}}"{description}, from=7-5, to=9-7]
		\arrow["{d_{1,1}}"{description}, from=7-7, to=9-7]
		\arrow["{\epsilon_{1,1}}"{description}, from=7-7, to=9-9]
	\end{tikzcd}\]
	\begin{notation}
		Let us list here the notations we use related to graded mixed objects:
		\begin{itemize}
			\item A graded module equipped with a mixed structure is called a graded mixed module.
			\item Morphisms of graded mixed modules are morphisms of graded modules commuting with the differentials.
			\item The 1-category of graded mixed modules and their morphisms is denoted by:
			\[\epsilon-dg^{gr}_{st}\]
			\item The 1-category $\epsilon-dg^{gr}_{st}$ is canonically equipped with a class of weak equivalences comprised of morphisms that induce quasi-isomorphisms in each weight.
			\item This 1-category can also be canonically equipped with a symmetric monoidal structure similar to that of $dg^{gr}_{st}$.
			\item A commutative algebra in this monoidal category is called a graded mixed algebra, and their 1-category will be denoted by:
			\[\epsilon-cdga^{gr}_{st} := CAlg(\epsilon-dg^{gr}_{st})\]
			Note that in such algebras, the mixed differential $\epsilon$ will moreover satisfy Leibniz' rule, as a consequence of the definition of the monoidal structure of $\epsilon-dg^{gr}_{st}$.
			\item The 1-category $\epsilon-cdga^{gr}_{st}$ also comes equipped with a class of weak equivalences inherited from $\epsilon-dg^{gr}_{st}$.
		\end{itemize}
	\end{notation}
	\begin{rmk}
		The 1-categories $dg^{gr}_{st}$, $cdga^{gr}_{st}$, $\epsilon-dg^{gr}_{st}$, and $\epsilon-cdga^{gr}_{st}$ all come equipped with canonical classes of weak equivalences, as described above.
		As in \cite[below Prop.~1.3]{PTVV}, we can localize these categories by their classes of weak equivalences, which yields $\infty$-categories.
	\end{rmk}
	\begin{notation}
		We will call $dg^{gr}$, $cdga^{gr}$, $\epsilon-dg^{gr}$, and $\epsilon-cdga^{gr}$ the $\infty$-categories obtained after localization of the corresponding above 1-categories by their classes of weak equivalences.
	\end{notation}
	\subsubsection{Relative Definitions}
	We will now introduce definitions relative to a different base $k$-algebra than solely $k$. Note that in \cite{PTVV}, only the relative case over $k$ is considered, but this generalization is mild.
	\begin{df} \label[df]{def_dg_gr_B}
		Let $B$ be a connective $k$-cdga, and equip it with the trivial grading, i.e.\ where $B$ is concentrated in weight 0. We can thus canonically consider that $B \in CAlg(dg^{gr}) = cdga^{gr}$.
		We will denote by $dg^{gr}(B)$ (resp.\ $cdga^{gr}(B)$) the $\infty$-category of $B$-modules (resp.\ $B$-algebras) in $dg^{gr}$, i.e.:
		\[dg^{gr}(B) := B-Mod(dg^{gr})\]
		\[cdga^{gr}(B) := B-CAlg(dg^{gr})\]
	\end{df}
	\begin{rmk}
		Note that the $\infty$-category $cdga^{gr}(B)$ is canonically equivalent to the slice $\infty$-category $B/cdga^{gr}$.
	\end{rmk}
	This defines graded objects over $B$. Similarly we can define relative graded mixed objects. For this let us first consider the main example of graded mixed algebras.
	\begin{ex} \label[ex]{deRham_graded_mixed_alg}
		Let $\phi : A \to B$ be any homomorphism of connective $k$-cdga's, making $B$ into an $A$-algebra.
		Recall that $\bb L_{B/A} \in B-Mod$ denotes the relative cotangent module of $\phi$, and consider the relative de Rham algebra:
		\[\DR(B/A) := \Sym_B(\bb L_{B/A}[1])\]
		It is canonically equipped with the structure of a graded mixed algebra.
		\begin{itemize}
			\item The application of $\Sym$\footnote{In this context, we are working in characteristic zero, so contrary to the case of positive characteristics, we do not need to distinguish between the permutation-invariant and the divided power algebras.} provides an extra grading, its weight grading, and a commutative algebra structure compatible with this grading. This constitutes its graded cdga structure.
			\item The de Rham differential satisfies the conditions to be a mixed structure on this graded cdga.
		\end{itemize}
		Therefore we can canonically consider that:
		\[\DR(B/A) \in CAlg(\epsilon-dg^{gr}) = \epsilon-cdga^{gr}\]
	\end{ex}
	\separation
	This allows us to define relative graded mixed algebras in the following way.
	\begin{df} \label[df]{def_rel_eps-cdga}
		The $\infty$-category of $A$-linear graded mixed $B$-modules is defined as the $\infty$-category of $\DR(B/A)$-modules in $\epsilon-dg^{gr}$, and is denoted by $\epsilon-dg^{gr}(B/A)$, i.e.:
		\[\epsilon-dg^{gr}(B/A) := \DR(B/A)-Mod(\epsilon-dg^{gr})\]
		The $\infty$-category of $A$-linear graded mixed $B$-algebras is defined as the slice $\infty$-category $\DR(B/A)/\epsilon-cdga^{gr}$, and is denoted by $\epsilon-cdga^{gr}(B/A)$, i.e.:
		\begin{align*}
			\epsilon-cdga^{gr}(B/A) :=& \, \DR(B/A)-CAlg(\epsilon-dg^{gr}) \\
			\simeq& \, \DR(B/A)/\epsilon-cdga^{gr}
		\end{align*}
	\end{df}
	\begin{rmk}
		Note that by considering a graded mixed module or algebra $F$ over $\DR(B/A)$, since $B$ is the weight-0 part of the latter, we can induce an actual $B$-module/algebra structure on $F$, but also an $A$-module/algebra structure inherited from the structural map $\phi : A \to B$. 
		Moreover, since the structural action of $\DR(B/A)$ on $F$ has to be compatible with the graded mixed structures, and since the mixed differential in $\DR(B/A)$ cancels on $A$, this requires that the mixed differential in $F$ also cancels on $A$, i.e.\ that it is $A$-linear.
		This provides the justification for the name of the $\infty$-category, and also describes the forgetful functors $\epsilon-dg^{gr}(B/A) \to dg^{gr}(B)$ and $\epsilon-cdga^{gr}(B/A) \to cdga^{gr}(B)$.
	\end{rmk}
	\begin{df} \label[df]{def_forget_gr}
		Let us define the two forgetful functors from graded mixed objects to graded objects:
		\begin{itemize}
			\item The forgetful functor which maps graded mixed modules to their underlying graded modules will be denoted by:
			\[-^{gr} : \epsilon-dg^{gr}(B/A) \to dg^{gr}(B)\]
			\item By abuse of notation, the analogous functor on graded mixed algebras will also be denoted by:
			\[-^{gr} : \epsilon-cdga^{gr}(B/A) \to cdga^{gr}(B)\]
		\end{itemize}
	\end{df}
	In fact we can break down these forgetful functors in a composition of two forgetful steps.
	\begin{notation}
		Let us set the following notations:
		\begin{itemize}
			\item The $\infty$-category of graded modules over $\DR(B/A)$ viewed as a graded algebra (without taking into account its de Rham differential) is:
			\[\DR(B/A)-Mod(dg^{gr})\]
			\item The functor forgetting the mixed differential of a graded mixed module is:
			\[\cancel \epsilon : \epsilon-dg^{gr}(B/A) \to \DR(B/A)-Mod(dg^{gr})\]
			\item The functor inducing a $B$-action on a graded $\DR(B/A)$-module along the structural map $B \to \DR(B/A)$ is:
			\[-^B : \DR(B/A)-Mod(dg^{gr}) \to dg^{gr}(B)\]
		\end{itemize}
	\end{notation}
	\begin{rmk} \label[rmk]{diagram_algebraic_defs}
		The forgetful functor $-^{gr} : \epsilon-dg^{gr}(B/A) \to dg^{gr}(B)$ is the composite of the forgetful functors $\cancel \epsilon : \epsilon-dg^{gr}(B/A) \to \DR(B/A)-Mod(dg^{gr})$ and $-^B : \DR(B/A)-Mod(dg^{gr}) \to dg^{gr}(B)$.
		To summarize, so far we have defined the following $\infty$-categories and forgetful functors:
		\[\begin{tikzcd}
			{\epsilon-dg^{gr}(B/A)} & {\DR(B/A)-Mod(dg^{gr})} & {dg^{gr}(B)}
			\arrow["{\cancel \epsilon}", from=1-1, to=1-2]
			\arrow["{-^{gr}}"', shift right=4, from=1-1, to=1-3]
			\arrow["{-^B}", from=1-2, to=1-3]
		\end{tikzcd}\]
		Note that all of these functors are conservative.
	\end{rmk}
	\subsubsection{Change of Base}
	We will now inspect what happens when we change our relative base $B/A$.
	\begin{construction} \label[construction]{pull-back}
		Consider a commutative square of connective $k$-cdga's:
		\[\begin{tikzcd}
			B & {B'} \\
			A & {A'}
			\arrow["b", from=1-1, to=1-2]
			\arrow["f", from=2-1, to=1-1]
			\arrow["a"', from=2-1, to=2-2]
			\arrow["{f'}"', from=2-2, to=1-2]
		\end{tikzcd}\]
		Note that this data can equivalently be phrased as a map $\phi : f \to f'$ in $\Arr(\dRings)$, where in this case $\phi$ is the data of $a$, $b$ and the commutativity of the square. In this situation, $\phi$ induces canonical morphisms $\DR(B/A) \xrightarrow{\beta} b^* \DR(B/A) \xrightarrow{\theta} \DR(B'/A')$ of graded mixed algebras. We will denote by $\psi : \DR(B/A) \to \DR(B'/A')$ the composite $\theta \circ \beta$. We can use it to construct a push-forward/pull-back adjunction between $\epsilon-cdga^{gr}(B'/A')$ and $\epsilon-cdga^{gr}(B/A)$, mildly generalizing \cite[Prop.~1.2.3]{AGI}:
		\begin{itemize}
			\item The push-forward functor $\phi_* : \epsilon-cdga^{gr}(B'/A') \to \epsilon-cdga^{gr}(B/A)$ takes a graded mixed algebra $F$ equipped with a structural map $\alpha' : \DR(B'/A') \to F$ and associates to it the composite $\DR(B/A) \xrightarrow{\psi} \DR(B'/A') \xrightarrow{\alpha'} F$.
			\item The pull-back functor $\phi^* : \epsilon-cdga^{gr}(B/A) \to \epsilon-cdga^{gr}(B'/A')$ takes a graded mixed algebra $F$ equipped with a structural map $\alpha : \DR(B/A) \to F$, and associates to it the graded mixed algebra $\phi^* F := b^* F \overset{\bb L}{\otimes}_{b^* \DR(B/A)} \DR(B'/A')$, equipped with its structural map $\DR(B'/A') \to \phi^* F$, or more explicitly it is the map labeled $\phi^* \alpha$ in the following pushout in graded mixed algebras
			\[\begin{tikzcd}
				{b^* \DR(B/A)} & {b^* F} \\
				{\DR(B'/A')} & {\phi^* F}
				\arrow["{b^* \alpha}", from=1-1, to=1-2]
				\arrow["\theta"', from=1-1, to=2-1]
				\arrow[from=1-2, to=2-2]
				\arrow["{\phi^* \alpha}"', from=2-1, to=2-2]
				\arrow["\lrcorner"{anchor=center, pos=0.125, rotate=180}, draw=none, from=2-2, to=1-1]
			\end{tikzcd}\]
			\item The two functors above form an adjunction $\phi^* \dashv \phi_*$.
		\end{itemize}
	\end{construction}
	
	Some $A$-linear graded mixed $B$-algebras serve as models for derived foliations on $\Spec(B)$ relative to $\Spec(A)$. We will now discuss the properties defining them.
	\begin{df}
		An $A$-linear graded mixed $B$-algebra $F$ will be called \textup{quasi-free} if the natural map $\Sym_{F(0)} F(1) \to F$ is an equivalence of graded algebras.
	\end{df}
	We can show that this property is stable by push-forward.
	\begin{lem}
		Let $F' \in \epsilon-cdga^{gr}(B'/A')$, and let $F = \phi_*(F') \in \epsilon-cdga^{gr}(B/A)$ denote its push-forward. If $F'$ is quasi-free, then so is $F$.
	\end{lem}
	\begin{proof}
		We have to show that the natural map $\Sym_{F(0)} F(1) \to F$ is an equivalence, assuming $\Sym_{F'(0)} F'(1) \to F'$ already is. But $F$ and $F'$ are canonically isomorphic as graded algebras (forgetting the respective $B$ and $B'$-algebra structures), so in particular $F(0)$ and $F'(0)$ are canonically equivalent as algebras, and $F(1)$ and $F'(1)$ also are canonically equivalent as $F(0)$-modules. Therefore, we can construct this commutative square of graded algebras:
		\[\begin{tikzcd}
			{\Sym_{F(0)} F(1)} & F \\
			{\Sym_{F'(0)} F'(1)} & {F'}
			\arrow["\sim"', from=1-1, to=2-1]
			\arrow[from=1-1, to=1-2]
			\arrow["\sim", from=2-1, to=2-2]
			\arrow["\sim", from=1-2, to=2-2]
		\end{tikzcd}\]
		and thus the top map is also an equivalence of graded algebras.
	\end{proof}
	To conclude this subsection, we give the definition of a derived foliation over $\Spec(B)$ relative to $\Spec(A)$.
	\begin{df} \label[df]{def_foliation}
		An $A$-linear graded mixed $B$-algebra $F$ will be called a \textup{derived foliation on $B$ relative to $A$} if it satisfies the following conditions:
		\begin{enumerate}
			\item The structural map $B \to F(0)$ is an equivalence
			\item $F$ is quasi-free
			\item $F(1)$ is perfect over $B$
		\end{enumerate}
		We define the $\infty$-category $\mcal Fol(B/A) \subset \epsilon-cdga^{gr}(B/A)^\op$ as the full subcategory spanned by these objects (note the $\op$). Moreover, the $B$-module $F(1)[-1]$ will be denoted by $\bb L_F$ and is called the cotangent complex of $F$ which is required to be perfect by condition 3, while condition 1.\ and 2.\ are equivalent to saying that $F \simeq \Sym_B(\bb L_F[1])$ canonically as graded $B$-algebras.
	\end{df}
	\begin{rmk}
		Note that more generally we can define non-perfect derived foliations by removing condition 3, which is the point of view taken in \cite[Def.~2.1.1.2,2.1.1.3]{toen2025derivedfoliations}. However, in this paper, all derived foliations will be assumed to be perfect, with the sole exception of the next remark, in which we will explicitly mention that we are considering non-perfect derived foliations to emphasize that we remove condition 3.
		
		A conceptual motivation for requiring condition 3, i.e.\ that the cotangent of the foliation is perfect, is that we wish to view foliations as a way to encode sub-bundles of the tangent bundle of the ambient space (equipped with a Lie algebra structure). Therefore, we need to be able to dualize the cotangent to obtain this tangent structure. However, we do not claim that condition 3.\ is the most general condition under which a notion of dual of the cotangent might exist, as for instance requiring instead that the cotangent is almost-perfect also yields a notion of dual, but only as a Pro-object in this case. Nonetheless, condition 3.\ is general enough for the examples of this paper, and more general conditions will not be studied here.
	\end{rmk}
	\begin{rmk} \label[rmk]{formula_pull-back}
		Firstly, conditions 1.\ and 2.\ of \cref{def_foliation} are stable by the pull-back of graded mixed algebras $\phi^*$ defined in \cref{pull-back}, so that the pull-back of graded mixed algebras restricts to a pull-back of non-perfect derived foliations, as described in the non-relative case in \cite[Def.~2.1.3]{toen2025derivedfoliations}, and as we will show explicitly below.
		
		Secondly, note that condition 1.\ in \cref{def_foliation} is not stable by the push-forward of graded mixed algebras $\phi_*$ defined in \cref{pull-back}, since it does not change the underlying graded mixed algebra and quasi-free $B'$-algebras need not remain quasi-free when viewed as $B$-algebras. Therefore, in general we cannot construct a push-forward of derived foliations naively as a restriction of the push-forward of graded mixed algebras.
		Instead, we will construct a push-forward of (perfect) derived foliations in \cref{Push-forward} using a different approach, through a point of view which we will introduce in the next subsection and develop in \cref{Foliations_Geometric}.
		Moreover, note that all the conditions of \cref{def_foliation} only depend on the graded structure of $F$, and not on its mixed structure.
		
		Coming back to our first point, we will now show explicitly that conditions 1.\ and 2.\ in \cref{def_foliation} are stable by the pull-back as defined in \cref{pull-back}. Consider a commutative square
		\[\begin{tikzcd}
			B & {B'} \\
			A & {A'}
			\arrow["b", from=1-1, to=1-2]
			\arrow["f", from=2-1, to=1-1]
			\arrow["a"', from=2-1, to=2-2]
			\arrow["{f'}"', from=2-2, to=1-2]
		\end{tikzcd}\]
		and denote by $\phi : f \to f'$ the corresponding map in $\Arr(\dRings)$. Then $\phi^*$ has the following description on non-perfect derived foliations: for $F \in \epsilon-cdga^{gr}(B/A)$, by definition, there is a pushout square in graded mixed algebras
		\[\begin{tikzcd}
			{b^* \DR(B/A)} & {b^* F} \\
			{\DR(B'/A')} & {\phi^* F}
			\arrow["{b^* \alpha}", from=1-1, to=1-2]
			\arrow["\theta"', from=1-1, to=2-1]
			\arrow[from=1-2, to=2-2]
			\arrow["{\phi^* \alpha}"', from=2-1, to=2-2]
			\arrow["\lrcorner"{anchor=center, pos=0.125, rotate=180}, draw=none, from=2-2, to=1-1]
		\end{tikzcd}\]
		If $F \simeq \Sym_B(\bb L_F[1])$ is a non-perfect derived foliation, we can identify $b^* F \simeq \Sym_{B'}(b^* \bb L_{\mcal F}[1])$, and moreover $b^* \DR(B/A) \simeq \Sym_{B'}(b^* \bb L_{B/A}[1])$, so we can rewrite this pushout as
		\[\begin{tikzcd}
			{\Sym_{B'}(b^* \bb L_{B/A}[1])} & {\Sym_{B'}(b^* \bb L_{\mcal F}[1])} \\
			{\Sym_{B'}(\bb L_{B'/A'}[1])} & {\phi^* F}
			\arrow["{b^* \alpha}", from=1-1, to=1-2]
			\arrow["\theta"', from=1-1, to=2-1]
			\arrow[from=1-2, to=2-2]
			\arrow["{\phi^* \alpha}"', from=2-1, to=2-2]
			\arrow["\lrcorner"{anchor=center, pos=0.125, rotate=180}, draw=none, from=2-2, to=1-1]
		\end{tikzcd}\]
		Since $\Sym_{B'}$ is a left adjoint, it commutes with colimits, and thus we obtain that $\phi^* F \simeq \Sym_{B'}(\bb L_{\phi^* \mcal F}[1])$, where $\bb L_{\phi^* \mcal F}$ is the following pushout of $B'$-modules
		\[\begin{tikzcd}
			{b^* \bb L_{B/A}} & {b^* \bb L_{\mcal F}} \\
			{\bb L_{B'/A'}} & {\bb L_{\phi^* \mcal F}}
			\arrow[from=1-1, to=1-2]
			\arrow[from=1-1, to=2-1]
			\arrow[from=1-2, to=2-2]
			\arrow[from=2-1, to=2-2]
			\arrow["\lrcorner"{anchor=center, pos=0.125, rotate=180}, draw=none, from=2-2, to=1-1]
		\end{tikzcd}\]
		As a result, we obtain that $\phi^* F$ is a non-perfect derived foliation on $B'$ relative to $A'$, and that we have the following formula to compute its cotangent:
		\[\bb L_{\phi^* \mcal F} \simeq b^* \bb L_{\mcal F} \oplus_{b^* \bb L_{B/A}} \bb L_{B'/A'}\]
		In general, $\bb L_{\phi^* \mcal F}$ may fail to remain perfect even if $\bb L_{\mcal F}$ was perfect. However, note that if $\bb L_{B/A}$ and $\bb L_{B'/A'}$ are also perfect, e.g.\ if $f$ and $f'$ are of finite presentation, then $b^* \bb L_{B/A}$ is also perfect because $b^*$ preserves perfects, and finally $\bb L_{\phi^* \mcal F}$ is perfect since perfect $B'$-modules are stable by pushouts. Moreover, even without these assumptions, $\bb L_{\phi^* \mcal F}$ may happen to be perfect for other reasons, as will be the case in the proof of \cref{main_prop-restate}.
	\end{rmk}
	\subsection{A Geometric Reformulation}
	We will now review results from \cite{Mou} (resp.\ \cite{HKR}), which provide equivalences between graded modules (resp.\ graded mixed modules) over a derived stack $X$ and quasi-coherent sheaves on derived stacks associated to $X$. In the next section, this will allow us to provide equivalences between bundles over these stacks and graded algebras (resp.\ graded mixed algebras) satisfying extra regularity conditions, such as derived foliations.
	
	More precisely, in the previous subsection we have introduced the following diagram (see \cref{diagram_algebraic_defs}):
	\[\begin{tikzcd}
		{\epsilon-dg^{gr}(B/A)} & {\DR(B/A)-Mod(dg^{gr})} & {dg^{gr}(B)}
		\arrow["{\cancel \epsilon}", from=1-1, to=1-2]
		\arrow["{-^{gr}}"', shift right=4, from=1-1, to=1-3]
		\arrow["{-^B}", from=1-2, to=1-3]
	\end{tikzcd}\]
	In this subsection, our goal is to view each $\infty$-category of this diagram as the $\infty$-category of quasi-coherent sheaves over a given stack, and each functor between them as either a pull-back or a push-forward along a given map between those stacks.
	\begin{rmk} \label[rmk]{rmk_QCohBGm}
		Consider first absolute graded modules and algebras.
		As is well known in algebraic geometry, and has been shown in derived geometry in \cite[Thm.~4.1]{Mou}, there is a canonical equivalence of symmetric monoidal $\infty$-categories $dg^{gr} \simeq \QCoh(B \bb G_m)$, which promotes to an equivalence between their $\infty$-categories of commutative algebras.
		Moreover, by \cite[Thm.~4.5.4.7]{HA}, the $\infty$-category of $\infty$-commutative algebras in the $\infty$-categorical localization of $dg^{gr}_{st}$ is equivalent to the $\infty$-categorical localization of $cdga^{gr}_{st} = \CAlg(dg^{gr}_{st})$, the latter of which is by definition the $\infty$-category of graded mixed algebras. In other words, if $\mcal W$ (resp. $\mcal W'$) is the class of weak equivalences of $dg^{gr}_{st}$ (resp. $cdga^{gr}_{st}$), then $\infty-\CAlg(dg^{gr}_{st}[\mcal W^{-1}]) \simeq \CAlg(dg^{gr}_{st})[\mcal W'^{-1}]$.
		In conclusion, we have canonical equivalences of symmetric monoidal $\infty$-categories:
		\begin{align*}
			dg^{gr} \simeq&\, \QCoh(B \bb G_m) \\
			cdga^{gr} \simeq&\, \CAlg(\QCoh(B \bb G_m))
		\end{align*}
	\end{rmk}
	To obtain a similar equivalence for graded mixed objects, let us now briefly introduce the group scheme called $\mcal H$, by a definition which is only valid since we are working in characteristic zero, as introduced in \cite[Rmk.~1.5]{PTVV}.\footnote{In positive characteristic, two groups $\mcal H$ and $\mcal H_{\pi}$ exist to encode respectively derived foliations and infinitesimal derived foliations, see \cite[\S1.2]{Hposchar}. The pull-back of these two groups over $\bb Q$ both coincide with the group $\mcal H$ as defined above, so the distinction disappears in our setting. We expect constructions similar to those of this paper to also apply for these two groups.}
	\begin{df} \label[df]{def_H}
		Let $\bb G_m$ and $\bb G_a$ denote respectively the multiplicative and additive group schemes. Recall that there is a canonical action of $\bb G_m$ on $\bb G_a$ by multiplication, which induces a canonical action of $\bb G_m$ on $B \bb G_a$ after delooping. Note that in derived algebraic geometry, all these group schemes including $B \bb G_a$ are affine and of finite type over $k$.
		We will call $\mcal H$ the semi-direct product $B \bb G_a \rtimes \bb G_m$, i.e.\ there is a canonical split short exact sequence of abelian group schemes:
		\[1 \to B \bb G_a \to \mcal H \leftrightarrows \bb G_m \to 1\]
		Moreover, note that $\mcal H$ is thus affine and of finite type over $k$ too.
	\end{df}
	\begin{rmk} \label[rmk]{QCohBH}
		Consider now absolute graded mixed modules and algebras.
		In this context, $\mcal H$-actions are equivalent to graded mixed objects.
		More precisely, $\mcal H$ is closely related to the object $(S^1_{Fil})^{gr}$ defined over $B \bb G_m$ in \cite[Def.~1.2.2,Def.~2.2.11]{HKR}, in the sense that the usual delooping $B \mcal H$ is equivalent to $B_{B \bb G_m} \left( (S^1_{Fil})^{gr} \right)$ the relative delooping of $(S^1_{Fil})^{gr}$ in the $\infty$-topos $\dSt/B \bb G_m$ (see \cite[Rmk.~2.3.65]{robalo2024applications}). Thus, by collecting the results \cite[Thm.~1.2.1(ii),Thm.~1.2.1(iii),Prop.~4.2.3(ii)]{HKR} there is a canonical equivalence of monoidal $\infty$-categories $\epsilon-dg^{gr} \simeq \QCoh(B \mcal H)$. Note that this equivalence has also been proven directly in characteristic zero in \cite[Prop.~2.7]{thmBH}.
		In conclusion, by a similar reasoning as for graded objects, we obtain the following equivalences of symmetric monoidal $\infty$-categories:
		\begin{align*}
			\epsilon-dg^{gr} \simeq&\, \QCoh(B \mcal H) \\
			\epsilon-cdga^{gr} \simeq&\, \CAlg(\QCoh(B \mcal H))
		\end{align*}
	\end{rmk}
	\separation
	We will now refine these equivalences to the relative cases, i.e.\ $cdga^{gr}(B)$ and $\epsilon-cdga^{gr}(B/A)$.
	\begin{rmk} \label[rmk]{equiv_gr_BGm}
		For the graded case, consider first the trivial $\bb G_m$-action on $\Spec(B)$: its quotient is $[\Spec(B)/\bb G_m] = \Spec(B) \times B \bb G_m$, and similarly as before we have symmetric monoidal equivalences:
		\begin{align*}
			dg^{gr}(B) \simeq&\, \QCoh(\Spec(B) \times B \bb G_m) \\
			cdga^{gr}(B) \simeq&\, \CAlg(\QCoh(\Spec(B) \times B \bb G_m))
		\end{align*}
	\end{rmk}
	For the graded mixed case, let us first introduce an auxiliary derived stack.
	\begin{df} \label[df]{def_Lgr}
		For any derived stack $Z$ over a derived stack $X$, we will call graded loops on $Z$ relative to $X$ the following mapping stack, named after the object $(L_{Fil} X)^{gr}$ in \cite[Thm.~5.1.3(2)]{HKR}:
		\[\mcal L^{gr}(Z/X) = \underline{\Map}_{X}(B \bb G_a \times X,Z)\]
		It comes equipped with a canonical $\mcal H$-action by acting on the source (more precisely, $B \bb G_a$ acts on the source, and $\bb G_m$ acts on $B \bb G_a$, which induces the $\mcal H$-action), and a canonical map $\mcal L^{gr}(Z/X) \to Z$ by evaluation at the base point of $B \bb G_a \times X$ over $X$.
		Moreover, when both $Z$ and $X$ are equipped with the trivial $\mcal H$-actions, we have that:
		\begin{enumerate}
			\item The induced map $\mcal L^{gr}(Z/X) \to X$ is $\mcal H$-equivariant.
			\item The canonical map $\mcal L^{gr}(Z/X) \to Z$ is not $\mcal H$-equivariant, but it is $\bb G_m$-equivariant for the induced actions.
		\end{enumerate}
		Moreover, when $Z = \Spec(B)$ and $X = \Spec(A)$ are affines, by \cite[Thm.~5.1.3(2)]{HKR}, applied over the base ring $A$, which is also in characteristic zero (i.e. a $\bb Q$-algebra) since $k$ is, this last stack is equivalent to $\Spec_B(\DR(B/A))$, where $\DR(B/A)$ is the (graded mixed) $B$-algebra introduced in \cref{deRham_graded_mixed_alg}.
	\end{df}
	\begin{rmk} \label[rmk]{equiv_eps_gr_BH}
		Note that since $\mcal L^{gr}(B/A)$ is relatively affine over the affine $\Spec(B)$, it is relatively affine over the point $* = \Spec(k)$\footnote{More concretely $\mcal L^{gr}(B/A) \simeq \Spec_k(\DR(B/A))$ as derived stacks (over $\Spec(k)$). In \cite[\S2.2]{AffineStacks}, this derived stack is simply denoted as $\Spec(\DR(B/A))$. However, this notation can be confusing here since $\DR(B/A)$ is not a derived ring in the sense of this paper, as it is not connective in general. Thus, the object $\Spec(\DR(B/A))$ as defined in the above sense is not a derived affine in our sense.}, and moreover the unique map to the point is also $\mcal H$-equivariant since it was already the case over $\Spec(A)$, and $\Spec(A) \to *$ is $\mcal H$-equivariant for the trivial actions.
		It follows that the quotient map $[\mcal L^{gr}(B/A)/\mcal H] \to B \mcal H$ is relatively affine, and in fact $[\mcal L^{gr}(B/A)/\mcal H] \simeq \Spec_{B \mcal H}(\DR(B/A))$ as derived stacks over $B \mcal H$, where $\DR(B/A)$ is viewed as an object of $\CAlg(\QCoh(B \mcal H))$ under the equivalence of \cref{QCohBH}. Therefore, by \cite[Chap.~I.3,Prop.~3.3.3]{gaitsgoryRozenblyumI} there is a canonical symmetric monoidal equivalence:
		\[\QCoh([\mcal L^{gr}(B/A)/\mcal H]) \simeq \DR(B/A)-\Mod(\QCoh(B \mcal H))\]
		which promotes to the following equivalences by \cref{QCohBH} and \cref{def_rel_eps-cdga}:
		\begin{align*}
			\epsilon-dg^{gr}(B/A) \simeq&\, \QCoh([\mcal L^{gr}(B/A)/\mcal H]) \\
			\epsilon-cdga^{gr}(B/A) \simeq&\, \CAlg(\QCoh([\mcal L^{gr}(B/A)/\mcal H]))
		\end{align*}
	\end{rmk}
	\begin{rmk}
		We can give a similar geometric description to the forgetful functor $-^{gr} : \epsilon-dg^{gr}(B/A) \to dg^{gr}(B)$.
		On the algebraic side, recall by \cref{diagram_algebraic_defs} that we have defined the following $\infty$-categories and forgetful functors:
		\[\begin{tikzcd}
			{\epsilon-dg^{gr}(B/A)} & {\DR(B/A)-Mod(dg^{gr})} & {dg^{gr}(B)}
			\arrow["{\cancel \epsilon}", from=1-1, to=1-2]
			\arrow["{-^{gr}}"', shift right=4, from=1-1, to=1-3]
			\arrow["{-^B}", from=1-2, to=1-3]
		\end{tikzcd}\]
		On the geometric side, note that the canonical map $\mcal L^{gr}(B/A) \to \Spec(B)$ is $\bb G_m$-equivariant as stated in \cref{def_Lgr}, and that $\bb G_m$ embeds as a subgroup of $\mcal H$ by the splitting of the short exact sequence in \cref{def_H}, so we can form the following span of quotient stacks:
		\[\begin{tikzcd}[cramped]
			{[\mcal L^{gr}(B/A)/\mcal H]} & {[\mcal L^{gr}(B/A)/\bb G_m]} & {\Spec(B) \times \bb G_m}
			\arrow["\pi"', from=1-2, to=1-1]
			\arrow["p", from=1-2, to=1-3]
		\end{tikzcd}\]
		Recall that we also have the following equivalences:
		\begin{itemize}
			\item $dg^{gr}(B) \simeq \QCoh(\Spec(B) \times B \bb G_m)$ (see \cref{equiv_gr_BGm})
			\item $\epsilon-dg^{gr}(B/A) \simeq \QCoh([\mcal L^{gr}(B/A)/\mcal H])$ (see \cref{equiv_eps_gr_BH})
			\item $\DR(B/A)-Mod(dg^{gr}) \simeq \QCoh([\mcal L^{gr}(B/A)/\bb G_m])$ (similarly)
		\end{itemize}
		Finally, note that under those equivalences, the functor $\cancel \epsilon$ corresponds to the pull-back $\pi^* : \QCoh([\mcal L^{gr}(B/A)/\mcal H]) \to \QCoh([\mcal L^{gr}(B/A)/\bb G_m])$, and that the functor $-^B$ corresponds to the push-forward $p_* : \QCoh([\mcal L^{gr}(B/A)/\bb G_m]) \to \QCoh(\Spec(B) \times \bb G_m)$.
	\end{rmk}
	\separation
	To conclude this subsection, we can summarize it in the following way.
	\begin{conclusion} \label[conclusion]{defs_conclusion}
		Starting from the following canonical span of quotient stacks:
		\[\begin{tikzcd}[cramped]
			{[\mcal L^{gr}(B/A)/\mcal H]} & {[\mcal L^{gr}(B/A)/\bb G_m]} & {\Spec(B) \times \bb G_m}
			\arrow["\pi"', from=1-2, to=1-1]
			\arrow["p", from=1-2, to=1-3]
		\end{tikzcd},\]
		we can form the following commutative diagram:
		\[\begin{tikzcd}[ampersand replacement=\&,cramped,column sep=scriptsize]
			{\QCoh([\mcal L^{gr}(B/A)/\mcal H])} \& {\QCoh([\mcal L^{gr}(B/A)/\bb G_m])} \& {\QCoh(\Spec(B) \times B \bb G_m)} \\
			{\epsilon-dg^{gr}(B/A)} \& {\DR(B/A)-Mod(dg^{gr})} \& {dg^{gr}(B)}
			\arrow["{\pi^*}", from=1-1, to=1-2]
			\arrow["\simeq"{marking, allow upside down}, draw=none, from=1-1, to=2-1]
			\arrow["{p_*}", from=1-2, to=1-3]
			\arrow["\simeq"{marking, allow upside down}, draw=none, from=1-2, to=2-2]
			\arrow["\simeq"{marking, allow upside down}, draw=none, from=1-3, to=2-3]
			\arrow["{\cancel \epsilon}", from=2-1, to=2-2]
			\arrow["{-^{gr}}"', shift right=4, from=2-1, to=2-3]
			\arrow["{-^B}", from=2-2, to=2-3]
		\end{tikzcd}\]
	\end{conclusion}
	\section{Derived Foliations as Equivariant Stacks} \label[section]{Foliations_Geometric}
	Compared to the previous section, in the present one we will let $Z = \Spec(B)$ be an affine derived stack over an affine base $X = \Spec(A)$, and replace all occurrences of $B$ (resp.\ $A$) by $Z$ (resp.\ $X$) in the notations of the previous section.
	
	Our goal in this section is to describe the $\infty$-category of $X$-linear derived foliations over $Z$ as a full subcategory of $\mcal H$-equivariant derived stacks over $\mcal L^{gr}(Z/X)$, the graded loop stack of $Z$ over $X$. We will then globalize this result to the general case where $Z \to X$ is locally of finite presentation, with a refinement in the relative derived Deligne-Mumford case.
	In order to prove it, we will follow a strategy similar to \cite{Mon}: we will construct an adjunction between $X$-linear graded-mixed algebras over $Z$ and the latter $\infty$-category, show that it is fully faithful, and describe its essential image.
	In \cref{Push-forward}, we will use the latter description of foliations to define their push-forward, as a restriction of the push-forward of stacks along a map $f : X \to Y$ satisfying extra conditions, after showing that it does preserve foliations.
	\subsection{Categorical Overview of the Setting}
	To make the goal of the present section more precise, this first subsection is dedicated to introducing all the $\infty$-categories of stacks, algebras of quasi-coherent sheaves, and all their relationships, half of which were introduced in \cref{Defs}. As a result, this subsection is heavy in definitions and notations, the point of which is to arrive at the diagram presented at the end in \cref{main_diagram}.
	\separation
	To begin, we will use the conclusion of \cref{Defs}.
	\begin{reminder} \label[reminder]{reminder_defs_conclusion}
		In \cref{defs_conclusion}, we arrived at the following setting.
		\begin{itemize}
			\item There is a canonical span of quotient stacks:
			\[\begin{tikzcd}[cramped]
				{[\mcal L^{gr}(Z/X)/\mcal H]} & {[\mcal L^{gr}(Z/X)/\bb G_m]} & {Z \times B \bb G_m}
				\arrow["\pi"', from=1-2, to=1-1]
				\arrow["p", from=1-2, to=1-3]
			\end{tikzcd}\]
			\item The following diagram is commutative:
			\[\begin{tikzcd}[cramped]
				{\QCoh([\mcal L^{gr}(Z/X)/\mcal H])} & {\QCoh([\mcal L^{gr}(Z/X)/\bb G_m])} & {\QCoh(Z \times B \bb G_m)} \\
				{\epsilon-dg^{gr}(Z/X)} & {\DR(Z/X)-Mod(dg^{gr})} & {dg^{gr}(Z)}
				\arrow["{\pi^*}", from=1-1, to=1-2]
				\arrow["\simeq"{marking, allow upside down}, draw=none, from=1-1, to=2-1]
				\arrow["{p_*}", from=1-2, to=1-3]
				\arrow["\simeq"{marking, allow upside down}, draw=none, from=1-2, to=2-2]
				\arrow["\simeq"{marking, allow upside down}, draw=none, from=1-3, to=2-3]
				\arrow["{\cancel \epsilon}", from=2-1, to=2-2]
				\arrow["{-^{gr}}"', shift right=4, from=2-1, to=2-3]
				\arrow["{-^B}", from=2-2, to=2-3]
			\end{tikzcd}\]
		\end{itemize}
	\end{reminder}
	This provides a description of the $\infty$-categories of quasi-coherent sheaves of those stacks in terms of the $\infty$-categories of graded and graded mixed modules introduced in \cref{alg_defs}.
	\begin{rmk} \label[rmk]{span_dst_desc}
		In parallel, we can inspect the $\infty$-categories of derived stacks over these stacks:
		\begin{itemize}
			\item By \cite[Prop.~1.3.5.1]{toenVezzosiHAGII}, the $\infty$-category $\dSt/[\mcal L^{gr}(Z/X)/\mcal H]$ is canonically equivalent to the $\infty$-category whose objects are derived stacks $T$ over $\mcal L^{gr}(Z/X)$, equipped with $\mcal H$-actions, whose structural maps $T \to \mcal L^{gr}(Z/X)$ are compatible with the actions of $\mcal H$ on both sides, i.e. 
			$\mcal H$-equivariant, and maps between such objects $T \to T'$ are also required to be compatible with the $\mcal H$-actions on both sides.
			\item The latter $\infty$-category is usually called the $\infty$-category of $\mcal H$-equivariant derived stacks over $\mcal L^{gr}(Z/X)$, and denoted by $\dSt/\mcal L^{gr}(Z/X) - \mcal H$.
			\item More concretely, the equivalence $\dSt/[\mcal L^{gr}(Z/X)/\mcal H] \xrightarrow{\sim} \dSt/\mcal L^{gr}(Z/X) - \mcal H$ of \cite[Prop.~1.3.5.1]{toenVezzosiHAGII} is given by taking an object $T \to [\mcal L^{gr}(Z/X)/\mcal H]$, and pulling it back along the canonical map $\mcal L^{gr}(Z/X) \to [\mcal L^{gr}(Z/X)/\mcal H]$. In the converse direction, it takes an $\mcal H$-equivariant derived stack $T' \to \mcal L^{gr}(Z/X)$, and associates to it the quotient map $[T'/\mcal H] \to [\mcal L^{gr}(Z/X)/\mcal H]$.
			\item Similar remarks can be made about $\dSt/[\mcal L^{gr}(Z/X)/\bb G_m]$ and $\dSt/(Z \times B \bb G_m)$.
		\end{itemize}
		In conclusion, we have:
		\begin{align*}
			\dSt/[\mcal L^{gr}(Z/X)/\mcal H] \simeq&\, \dSt/\mcal L^{gr}(Z/X) - \mcal H \\
			\dSt/[\mcal L^{gr}(Z/X)/\bb G_m] \simeq&\, \dSt/\mcal L^{gr}(Z/X) - \bb G_m \\
			\dSt/(Z \times B \bb G_m) \simeq&\, \dSt/Z - \bb G_m
		\end{align*}
	\end{rmk}
	Our current goal is to relate these $\infty$-categories with the $\infty$-categories of quasi-coherent sheaves described in \cref{reminder_defs_conclusion}.
	To proceed, let us now also recall a standard duality between derived stacks defined over a base stack $S$ and commutative algebras of quasi-coherent sheaves over $S$.
	\begin{reminder} \label[reminder]{qcoh_dst_adj}
		Let $S$ be any derived stack.
		\begin{itemize}
			\item Let $p : T \to S \in \dSt/S$ be a map of derived stacks.
			We define $\mcal O_{/S}(p) = p_*(\mcal O_T) \in \CAlg(\QCoh(S))$, where $\mcal O_T$ is the structure sheaf of $T$, and where $p_* : \CAlg(\QCoh(T)) \to \CAlg(\QCoh(S))$ is the push-forward functor of quasi-coherent sheaves.
			This canonically promotes to a contravariant functor, called the $S$-relative structure sheaf, and denoted by:
			\[\mcal O_{/S} : \dSt/S \to \CAlg(\QCoh(S))^{\op}\]
			\item Conversely, let $A \in \CAlg(\QCoh(S))$ be a quasi-coherent algebra.
			We define $\Spec_S(A) \in \dSt/S$ to be the derived stack on $S$ defined on the generating site $\dAff/S$ by mapping $q : \Spec(B) \to S$ to $\Map_{B-Alg}(q^*A,B)$.
			This canonically promotes to a contravariant functor, called the $S$-relative spectrum, and denoted by:
			\[\Spec_S : \CAlg(\QCoh(S))^{\op} \to \dSt/S\]
			\item These form the canonical adjoint pair:
			\[\begin{tikzcd}[cramped]
				{\CAlg(\QCoh(S))^{\op}} & {\dSt/S}
				\arrow[""{name=0, anchor=center, inner sep=0}, "{\Spec_S}"', shift right, from=1-1, to=1-2]
				\arrow[""{name=1, anchor=center, inner sep=0}, "{\mcal O_{/S}}"', shift right, from=1-2, to=1-1]
				\arrow["\bot"{description}, draw=none, from=1, to=0]
			\end{tikzcd}\]
		\end{itemize}
	\end{reminder}
	Moreover, we can show that this construction itself is functorial in $S$, in the sense of the following remark.
	\begin{reminder} \label[reminder]{base_chg_adj}
		Let $f : S' \to S$ be any map of derived stacks.
		\begin{itemize}
			\item There is a forgetful/pull-back/push-forward adjunction triplet $f_! \dashv f^* \dashv f_*$ between their $\infty$-categories of relative derived stacks:
			\[\begin{tikzcd}[cramped]
				{\dSt/S'} && {\dSt/S}
				\arrow[""{name=0, anchor=center, inner sep=0}, "{f_!}", shift left=5, from=1-1, to=1-3]
				\arrow[""{name=1, anchor=center, inner sep=0}, "{f_*}"', shift right=5, from=1-1, to=1-3]
				\arrow[""{name=2, anchor=center, inner sep=0}, "{f^*}"{description}, from=1-3, to=1-1]
				\arrow["\dashv"{marking, allow upside down, pos=0.4}, draw=none, from=0, to=2]
				\arrow["\dashv"{marking, allow upside down, pos=0.6}, draw=none, from=2, to=1]
			\end{tikzcd}\]
			where, given $U \to S \in \dSt/S$, we define $f^*(U)$ to be its pull-back along $f$ equipped with its canonical projection to $S'$:
			\[f^*(U) := U \times_{S} S' \to S' \in \dSt/S'\]
			its left adjoint $f_! : \dSt/S' \to \dSt/S$ is given simply by composition along $f$, in the sense that if $U \to S'$ is an object of $\dSt/S'$, then we have:
			\[f_!(U') := U \to S' \xrightarrow{f} S \in \dSt/S\]
			Moreover, $f^*$ also admits a right adjoint $f_* : \dSt/S' \to \dSt/S$, called the push-forward, direct image, or Weil restriction along $f$, given on an object $U' \to S'$ by the following functor of points:
			\[\forall \Spec(A) \to S \in \dAff/S, \; f_*(U')(A) := \Map_{\dSt/S'}(\Spec(A) \times_S S', U')\]
			More generally, this extends to arbitrary stacks over $S$, so $f_*(U')$ satisfies the following universal property in $\dSt/S$:
			\[\forall T \to S, \; \Map_S(T,f_*(U')) \simeq \Map_{S'}(T \times_S S',U')\]
			which recovers exactly the statement that $f_*$ is the right adjoint of $f^*$.
			\item There is a pull-back/push-forward adjunction between their $\infty$-categories of algebras of quasi-coherent sheaves:
			\[\begin{tikzcd}[cramped]
				{\CAlg(\QCoh(S'))^\op} \\
				{\CAlg(\QCoh(S))^\op}
				\arrow[""{name=0, anchor=center, inner sep=0}, "{f_*}"', shift right, from=1-1, to=2-1]
				\arrow[""{name=1, anchor=center, inner sep=0}, "{f^*}"', shift right, from=2-1, to=1-1]
				\arrow["\dashv"{description}, draw=none, from=0, to=1]
			\end{tikzcd}\]
			where, given $A \in \CAlg(\QCoh(S))$, we define:
			\[f^*(A) := A \otimes_{\mcal O_S} \mcal O_{S'} \in \CAlg(\QCoh(S'))\]
			Note that usually $f^*$ is considered as the left-adjoint to $f_*$, but in this diagram it is denoted as the right-adjoint to $f_*$ instead, because we are taking opposites of $\infty$-categories on both sides.
			\item The following diagrams commute:
			\[\begin{tikzcd}[cramped]
				{\CAlg(\QCoh(S'))^\op} & {\dSt/S'} && {\CAlg(\QCoh(S'))^\op} & {\dSt/S'} \\
				{\CAlg(\QCoh(S))^\op} & {\dSt/S} && {\CAlg(\QCoh(S))^\op} & {\dSt/S}
				\arrow["{f_*}"', from=1-1, to=2-1]
				\arrow["{\mcal O_{/S'}}"', from=1-2, to=1-1]
				\arrow["{f_!}", from=1-2, to=2-2]
				\arrow["{\Spec_{S'}}"', from=1-4, to=1-5]
				\arrow["{\mcal O_{/S}}"', from=2-2, to=2-1]
				\arrow["{f^*}"', shift left, from=2-4, to=1-4]
				\arrow["{\Spec_S}"', from=2-4, to=2-5]
				\arrow["{f^*}"', from=2-5, to=1-5]
			\end{tikzcd}\]
		\end{itemize}
		On the left-hand side, commutativity is simply the statement that for any $p : T \to S'$, we have $(f \circ p)_*(\mcal O_T) \simeq f_*(p_* \mcal O_T)$. Moreover, the commutativity of the left-hand side implies the commutativity of the right-hand side, because the functors on this side are right-adjoints to the functors on the left-hand side, and right-adjoints are unique.
	\end{reminder}
	\begin{rmk}
		In particular, we can apply these constructions to our span from \cref{defs_conclusion}, i.e.\ apply it to:
		\[\begin{tikzcd}[cramped]
			{[\mcal L^{gr}(Z/X)/\mcal H]} & {[\mcal L^{gr}(Z/X)/\bb G_m]} & {Z \times B \bb G_m}
			\arrow["\pi"', from=1-2, to=1-1]
			\arrow["p", from=1-2, to=1-3]
		\end{tikzcd}\]
		Moreover, we can describe the corresponding $\infty$-categories of quasi-coherent sheaves by \cref{reminder_defs_conclusion}, and of derived stacks by \cref{span_dst_desc}.
		We thus obtain the following diagram, which a priori need not commute:
		\[\begin{tikzcd}[cramped]
			{[\mcal L^{gr}(Z/X)/\mcal H]} & {\epsilon-cdga^{gr}(Z/X)^\op} &&& {\dSt/\mcal L^{gr}(Z/X) - \mcal H} \\
			\\
			{[\mcal L^{gr}(Z/X)/\bb G_m]} & {cdga^{gr}(\DR(Z/X))^\op} &&& {\dSt/\mcal L^{gr}(Z/X) - \bb G_m} \\
			\\
			{[Z/\bb G_m]} & {cdga^{gr}(Z)^\op} &&& {\dSt/Z - \bb G_m}
			\arrow[""{name=0, anchor=center, inner sep=0}, "{\Spec_{[\mcal L^{gr}(Z/X)/\mcal H]}}"', shift right, from=1-2, to=1-5]
			\arrow["{\pi^*}", from=1-2, to=3-2]
			\arrow[""{name=1, anchor=center, inner sep=0}, "{\mcal O_{[\mcal L^{gr}(Z/X)/\mcal H]}}"', shift right, from=1-5, to=1-2]
			\arrow["{\pi^*}", from=1-5, to=3-5]
			\arrow["\pi"', from=3-1, to=1-1]
			\arrow[squiggly, from=3-1, to=3-2]
			\arrow["p", from=3-1, to=5-1]
			\arrow[""{name=2, anchor=center, inner sep=0}, "{\Spec_{[\mcal L^{gr}(Z/X)/\bb G_m]}}"', shift right, from=3-2, to=3-5]
			\arrow["{p_*}", from=3-2, to=5-2]
			\arrow[""{name=3, anchor=center, inner sep=0}, "{\mcal O_{[\mcal L^{gr}(Z/X)/\bb G_m]}}"', shift right, from=3-5, to=3-2]
			\arrow["{p_!}", from=3-5, to=5-5]
			\arrow[""{name=4, anchor=center, inner sep=0}, "{\Spec_{[Z/\bb G_m]}}"', shift right, from=5-2, to=5-5]
			\arrow[""{name=5, anchor=center, inner sep=0}, "{\mcal O_{[Z/\bb G_m]}}"', shift right, from=5-5, to=5-2]
			\arrow["\bot"{description}, draw=none, from=1, to=0]
			\arrow["\bot"{description}, draw=none, from=3, to=2]
			\arrow["\bot"{description}, draw=none, from=5, to=4]
		\end{tikzcd}\]
		Note that so far, by the last point of \cref{base_chg_adj}, we only know that the following diagrams are commutative:
		\[\begin{tikzcd}[cramped]
			{\epsilon-cdga^{gr}(Z/X)^\op} && {\dSt/\mcal L^{gr}(Z/X) - \mcal H} \\
			{cdga^{gr}(\DR(Z/X))^\op} && {\dSt/\mcal L^{gr}(Z/X) - \bb G_m} \\
			\\
			{cdga^{gr}(\DR(Z/X))^\op} && {\dSt/\mcal L^{gr}(Z/X) - \bb G_m} \\
			{cdga^{gr}(Z)^\op} && {\dSt/Z - \bb G_m}
			\arrow["{\Spec_{[\mcal L^{gr}(Z/X)/\mcal H]}}", shift right, from=1-1, to=1-3]
			\arrow["{\pi^*}"', from=1-1, to=2-1]
			\arrow["{\pi^*}", from=1-3, to=2-3]
			\arrow["{\Spec_{[\mcal L^{gr}(Z/X)/\bb G_m]}}"', from=2-1, to=2-3]
			\arrow["{p_*}"', from=4-1, to=5-1]
			\arrow["{\mcal O_{[\mcal L^{gr}(Z/X)/\bb G_m]}}"', from=4-3, to=4-1]
			\arrow["{p_!}", from=4-3, to=5-3]
			\arrow["{\mcal O_{[Z/\bb G_m]}}", from=5-3, to=5-1]
		\end{tikzcd}\]
	\end{rmk}
	In the above diagrams, let us now give special notations to the functors involved.
	\begin{notation}
		For the functors obtained from \cref{qcoh_dst_adj} which are of the form $\dSt/S : \mcal O_{/S} \dashv \Spec_S : \CAlg(\QCoh(S))^\op$, we introduce notations in the following cases:
		\begin{itemize}
			\item When $S = [\mcal L^{gr}(Z/X)/\mcal H]$, we denote it by: $\mcal O_{\epsilon-gr} \dashv \Spec^{\epsilon-gr}$.
			\item When $S = [\mcal L^{gr}(Z/X)/\bb G_m]$, we denote it by: $\mcal O_{\DR-gr} \dashv \Spec^{\DR-gr}$.
			\item When $S = [Z/\bb G_m]$, we denote it by: $\mcal O_{gr} \dashv \Spec^{gr}$.
			Note that $\mcal O_{gr}$ was introduced and studied in \cite[Def.~1.11]{Mon}.
		\end{itemize}
	\end{notation}
	\begin{notation}
		For the functors obtained from \cref{base_chg_adj} which are of the form $\CAlg(\QCoh(S))^\op : f_* \dashv f^* : \CAlg(\QCoh(S'))^\op$ and $\dSt/S : f_! \dashv f^* : \dSt/S'$, we introduce notations in the following cases:
		\begin{itemize}
			\item When $f = \pi : [\mcal L^{gr}(Z/X)/\mcal H] \to [\mcal L^{gr}(Z/X)/\bb G_m]$, we introduce:\\
			$\cancel \epsilon : \epsilon-cdga^{gr}(Z/X)^\op \to cdga^{gr}(\DR(Z/X))^\op := \pi^* \textup{as in \cref{defs_conclusion}}$\\
			$-^{\bb G_m} : \dSt/\mcal L^{gr}(Z/X) - \mcal H \to \dSt/\mcal L^{gr}(Z/X) - \bb G_m := \pi^*$
			\item When $f = p : [\mcal L^{gr}(Z/X)/\bb G_m] \to [Z/\bb G_m]$, we introduce:\\
			$-^Z : cdga^{gr}(\DR(Z/X))^\op \to cdga^{gr}(Z)^\op := p_* \textup{as in \cref{defs_conclusion}}$\\
			$-^Z : \dSt/\mcal L^{gr}(Z/X) - \bb G_m \to \dSt/Z - \bb G_m := p_!$
		\end{itemize}
		Finally, we also introduce notations for the following compositions:\\
		$-^{gr} : \epsilon-cdga^{gr}(Z/X)^\op \to cdga^{gr}(Z)^\op := -^Z \circ \cancel \epsilon, \textup{ following \cref{diagram_algebraic_defs}}$\\
		$-^{Z-\bb G_m} : \dSt/\mcal L^{gr}(Z/X) - \mcal H \to \dSt/Z - \bb G_m := -^Z \circ -^{\bb G_m}$
	\end{notation}
	\begin{rmk} \label[rmk]{first_adj_diagram}
		So far, we have the following diagram, which a priori need not commute:
		\[\begin{tikzcd}
			{} & {\epsilon-cdga^{gr}(Z/X)^\op} & {\dSt/\mcal L^{gr}(Z/X) - \mcal H} & {} \\
			& {cdga^{gr}(\DR(Z/X))^\op} & {\dSt/\mcal L^{gr}(Z/X) - \bb G_m} \\
			{} & {cdga^{gr}(Z)^\op} & {\dSt/Z - \bb G_m} & {}
			\arrow["{-^{gr}}"{description}, shift left=5, from=1-1, to=3-1]
			\arrow[""{name=0, anchor=center, inner sep=0}, "{\Spec^{\epsilon-gr}}"', shift right, from=1-2, to=1-3]
			\arrow["{\cancel \epsilon}"', from=1-2, to=2-2]
			\arrow[""{name=1, anchor=center, inner sep=0}, "{\mcal O_{\epsilon-gr}}"', shift right, from=1-3, to=1-2]
			\arrow["{-^{\bb G_m}}", from=1-3, to=2-3]
			\arrow["{-^{Z-\bb G_m}}"{description}, shift right=5, from=1-4, to=3-4]
			\arrow[""{name=2, anchor=center, inner sep=0}, "{\Spec^{\DR-gr}}"', shift right, from=2-2, to=2-3]
			\arrow["{-^Z}"', from=2-2, to=3-2]
			\arrow[""{name=3, anchor=center, inner sep=0}, "{\mcal O_{\DR-gr}}"', shift right, from=2-3, to=2-2]
			\arrow["{-^Z}", from=2-3, to=3-3]
			\arrow[""{name=4, anchor=center, inner sep=0}, "{\Spec^{gr}}"', shift right, from=3-2, to=3-3]
			\arrow[""{name=5, anchor=center, inner sep=0}, "{\mcal O_{gr}}"', shift right, from=3-3, to=3-2]
			\arrow["\bot"{description}, draw=none, from=1, to=0]
			\arrow["\bot"{description}, draw=none, from=3, to=2]
			\arrow["\bot"{description}, draw=none, from=5, to=4]
		\end{tikzcd}\]
		
		As above, so far we know only of the following commutativity:
		\begin{align*}
			\Spec^{\DR-gr} \circ \; \cancel \epsilon &\simeq -^{\bb G_m} \circ \Spec^{\epsilon-gr}\\
			\mcal O_{gr} \circ -^Z &\simeq -^Z \circ \mcal O_{\DR-gr}
		\end{align*}
	\end{rmk}
	To complete the above diagram, let us recall and introduce additional $\infty$-categories and functors involved in our setting.
	\begin{reminder}
		We recall following objects:
		\begin{itemize}
			\item The $\infty$-category $\mcal Fol(Z/X) \subset \epsilon-cdga^{gr}(Z/X)^\op$ is the full subcategory spanned by graded mixed algebras satisfying the conditions of \cref{def_foliation}.
			\item We recall the following from \cite[Prop.~2.1]{Mon}:\\
			- There is a forgetful functor:
			\[-(1) : cdga^{gr}(Z)^\op \to \QCoh(Z)^\op,\]
			which takes a graded $Z$-algebra and maps it to its part sitting in weight 1.\\
			- The functor $-(1)$ is left-adjoint to the functor:
			\[\Sym_{\mcal O_Z} : \QCoh(Z)^\op \to cdga^{gr}(Z)^\op\]
			which comes equipped with the canonical grading $\Sym_{\mcal O_Z} = \underset{n \in \bb N}{\oplus} \Sym^i_{\mcal O_Z}$.\\
			- The composition $\Spec^{gr} \circ \Sym_{\mcal O_Z}$ is denoted by $\bb V : \QCoh(Z)^\op \to \dSt/Z - \bb G_m$.\\
			- The essential image of $\bb V$ is called the $\infty$-category of linear stacks over $Z$.
		\end{itemize}
	\end{reminder}
	Let us now introduce the following objects and notations.
	\begin{df} \label[df]{df-perflinst}
		Let $\PerfLinSt_Z \subset \dSt/Z - \bb G_m$ be the $\infty$-category of perfect linear stacks over $Z$, defined as the essential image of $\bb V_{|\Perf(Z)^\op} : \Perf(Z)^\op \to \dSt/Z - \bb G_m$.
	\end{df}
	\begin{df}
		Let $\PerfLinSt_Z/\mcal L^{gr}(Z/X) - \mcal H$ be defined as the full subcategory of $\dSt/\mcal L^{gr}(Z/X) - \mcal H$ spanned by stacks whose image under the forgetful functor $-^{Z - \bb G_m} : \dSt/\mcal L^{gr}(Z/X) - \mcal H \to \dSt/Z - \bb G_m$ belongs to $\PerfLinSt_Z$.
		Equivalently, we can define $\PerfLinSt_Z/\mcal L^{gr}(Z/X) - \mcal H$ using the following pull-back square of $\infty$-categories:
		\[\begin{tikzcd}
			{\PerfLinSt_Z/\mcal L^{gr}(Z/X) - \mcal H} & {\PerfLinSt_Z} \\
			{\dSt/\mcal L^{gr}(Z/X) - \mcal H} & {\dSt/Z - \bb G_m}
			\arrow[hook, from=1-2, to=2-2]
			\arrow[hook, from=1-1, to=2-1]
			\arrow["{-^{Z - \bb G_m}}", from=2-1, to=2-2]
			\arrow[from=1-1, to=1-2]
			\arrow["\lrcorner"{anchor=center, pos=0.125, rotate=0}, draw=none, from=1-1, to=2-2]
		\end{tikzcd}\]
	\end{df}
	\separation
	To conclude this subsection, let us combine all its data.
	\begin{conclusion} \label[conclusion]{main_diagram}
		We have constructed the following diagram, which a priori need not commute:
		\[\begin{tikzcd}
			& {\mcal Fol(Z/X)} & {\PerfLinSt_Z/\mcal L^{gr}(Z/X) - \mcal H} \\
			{} & {\epsilon-cdga^{gr}(Z/X)^\op} & {\dSt/\mcal L^{gr}(Z/X) - \mcal H} & {} \\
			& {cdga^{gr}(\DR(Z/X))^\op} & {\dSt/\mcal L^{gr}(Z/X) - \bb G_m} \\
			{} & {cdga^{gr}(Z)^\op} & {\dSt/Z - \bb G_m} & {} \\
			& {\QCoh(Z)^\op}
			\arrow[hook, from=1-2, to=2-2]
			\arrow[hook, from=1-3, to=2-3]
			\arrow["{-^{gr}}"{description}, shift left=5, from=2-1, to=4-1]
			\arrow[""{name=0, anchor=center, inner sep=0}, "{\Spec^{\epsilon-gr}}"', shift right, from=2-2, to=2-3]
			\arrow["{\cancel \epsilon}"', from=2-2, to=3-2]
			\arrow[""{name=1, anchor=center, inner sep=0}, "{\mcal O_{\epsilon-gr}}"', shift right, from=2-3, to=2-2]
			\arrow["{-^{\bb G_m}}", from=2-3, to=3-3]
			\arrow["{-^{Z-\bb G_m}}"{description}, shift right=5, from=2-4, to=4-4]
			\arrow[""{name=2, anchor=center, inner sep=0}, "{\Spec^{\DR-gr}}"', shift right, from=3-2, to=3-3]
			\arrow["{-^Z}"', from=3-2, to=4-2]
			\arrow[""{name=3, anchor=center, inner sep=0}, "{\mcal O_{\DR-gr}}"', shift right, from=3-3, to=3-2]
			\arrow["{-^Z}", from=3-3, to=4-3]
			\arrow[""{name=4, anchor=center, inner sep=0}, "{\Spec^{gr}}"', shift right, from=4-2, to=4-3]
			\arrow[""{name=5, anchor=center, inner sep=0}, "{-(1)}"', shift right, from=4-2, to=5-2]
			\arrow[""{name=6, anchor=center, inner sep=0}, "{\mcal O_{gr}}"', shift right, from=4-3, to=4-2]
			\arrow[""{name=7, anchor=center, inner sep=0}, "{\Sym_{\mcal O_Z}}"', shift right, from=5-2, to=4-2]
			\arrow["{\bb V}"', from=5-2, to=4-3]
			\arrow["\bot"{description}, draw=none, from=1, to=0]
			\arrow["\bot"{description}, draw=none, from=3, to=2]
			\arrow["\dashv"{description}, draw=none, from=5, to=7]
			\arrow["\bot"{description}, draw=none, from=6, to=4]
		\end{tikzcd}\]
		In this diagram, so far we only know the following commutativity:
		\begin{align*}
			\Spec^{\DR-gr} \circ \; \cancel \epsilon &\simeq -^{\bb G_m} \circ \Spec^{\epsilon-gr}\\
			\mcal O_{gr} \circ -^Z &\simeq -^Z \circ \mcal O_{\DR-gr}\\
			\bb V &:= \Spec^{gr} \circ \Sym_{\mcal O_Z}
		\end{align*}
	\end{conclusion}
	\separation
	By the end of the next subsection, in \cref{equiv_fol_dst}, we will prove that the adjunction $\mcal O_{\epsilon-gr} \dashv \Spec^{\epsilon-gr}$ restricts to an adjoint equivalence of $\infty$-categories:
	\[\mcal Fol(Z/X) \simeq \PerfLinSt_Z/\mcal L^{gr}(Z/X) - \mcal H\]
	To prove it, we will show that the vertical forgetful functors are compatible with the horizontal adjunctions, in a sense we will discuss in the next subsection. Essentially, we will prove that the diagram is in fact commutative.
	In the final subsection, \cref{equiv_global}, we will globalize these results to non-affine cases, more specifically to the general case of derived stacks and with a refinement in the case of derived Deligne-Mumford stacks.
	\subsection{Equivalence in the Affine Case} \label[subsection]{equiv-affine-case}
	As stated above, the equivalence in \cref{equiv_fol_dst} will follow from the fact that the squares of the diagram in \cref{main_diagram} satisfy a compatibility condition. This condition is known as the Beck-Chevalley condition (note that squares satisfying this condition are also called adjointable squares in the terminology of \cite[Def.~4.7.4.13]{HA}). Details on Beck-Chevalley conditions are recalled in \cref{beck-chevalley-section}, where the general definitions and properties are recalled in \cref{beck-chevalley-defs}, while their appearance in derived algebraic geometry are recalled in \cref{beck-chevalley-dag}. Intuitively, this condition means that the square formed by the left adjoints, and the square formed by the right adjoints, are both commutative, and that the vertical functors are compatible with the units and co-units of these adjunctions.
	In this subsection, this condition will allow us to prove on the one hand that the adjunction $\mcal O_{\epsilon-gr} \dashv \Spec^{\epsilon-gr}$ does restrict to an adjunction between derived foliations and perfect linear stacks, and on the other hand that checking that the unit and co-unit of this adjunction are equivalences reduces to checking that it is the case for the unit and co-unit of $\mcal O_{gr} \dashv \Spec^{gr}$, which in turn will follow from results of \cite[Cor.~2.9,Prop.~2.11]{Mon}. We will start this subsection by applying the results of \cref{beck-chevalley-section}.
	\begin{rmk} \label[rmk]{past_bc} \label[rmk]{p_pi_bc}
		Consider the maps $\pi$ and $p$ in the canonical span of \cref{defs_conclusion}:
		\[\begin{tikzcd}
			{[\mcal L^{gr}(Z/X)/\mcal H]} & {[\mcal L^{gr}(Z/X)/\bb G_m]} & {Z \times B \bb G_m}
			\arrow["\pi"', from=1-2, to=1-1]
			\arrow["p", from=1-2, to=1-3]
		\end{tikzcd}\]
		which form the following two squares, which a priori need not commute:
		\[\begin{tikzcd}
			{\epsilon-cdga^{gr}(Z/X)^\op} & {\dSt/\mcal L^{gr}(Z/X) - \mcal H} \\
			{\DR(Z/X)/cdga^{gr}(Z)^\op} & {\dSt/\mcal L^{gr}(Z/X)-\bb G_m} \\
			{cdga^{gr}(Z)^\op} & {\dSt/Z - \bb G_m}
			\arrow[""{name=0, anchor=center, inner sep=0}, "{\Spec^{\epsilon-gr}}"', shift right, from=1-1, to=1-2]
			\arrow["{\pi^*}", from=1-1, to=2-1]
			\arrow[""{name=1, anchor=center, inner sep=0}, "{\mcal O_{\epsilon-gr}}"', shift right, from=1-2, to=1-1]
			\arrow["{\pi^*}", from=1-2, to=2-2]
			\arrow[""{name=2, anchor=center, inner sep=0}, "{\Spec_{\DR-gr}}"', shift right, from=2-1, to=2-2]
			\arrow["{p_*}", from=2-1, to=3-1]
			\arrow[""{name=3, anchor=center, inner sep=0}, "{\mcal O_{\DR-gr}}"', shift right, from=2-2, to=2-1]
			\arrow["{p_!}", from=2-2, to=3-2]
			\arrow[""{name=4, anchor=center, inner sep=0}, "{\Spec^{gr}}"', shift right, from=3-1, to=3-2]
			\arrow[""{name=5, anchor=center, inner sep=0}, "{\mcal O_{gr}}"', shift right, from=3-2, to=3-1]
			\arrow["\bot"{description}, draw=none, from=1, to=0]
			\arrow["\bot"{description}, draw=none, from=3, to=2]
			\arrow["\bot"{description}, draw=none, from=5, to=4]
		\end{tikzcd}\]
		In these squares, by the last point of \cref{base_chg_adj}, so far we only know that the following squares commute:
		\[\begin{tikzcd}[cramped]
			{\epsilon-cdga^{gr}(Z/X)^\op} && {\dSt/\mcal L^{gr}(Z/X) - \mcal H} \\
			{cdga^{gr}(\DR(Z/X))^\op} && {\dSt/\mcal L^{gr}(Z/X) - \bb G_m} \\
			\\
			{cdga^{gr}(\DR(Z/X))^\op} && {\dSt/\mcal L^{gr}(Z/X) - \bb G_m} \\
			{cdga^{gr}(Z)^\op} && {\dSt/Z - \bb G_m}
			\arrow["{\Spec_{[\mcal L^{gr}(Z/X)/\mcal H]}}", shift right, from=1-1, to=1-3]
			\arrow["{\pi^*}"', from=1-1, to=2-1]
			\arrow["{\pi^*}", from=1-3, to=2-3]
			\arrow["{\Spec_{[\mcal L^{gr}(Z/X)/\bb G_m]}}"', from=2-1, to=2-3]
			\arrow["{p_*}"', from=4-1, to=5-1]
			\arrow["{\mcal O_{[\mcal L^{gr}(Z/X)/\bb G_m]}}"', from=4-3, to=4-1]
			\arrow["{p_!}", from=4-3, to=5-3]
			\arrow["{\mcal O_{[Z/\bb G_m]}}", from=5-3, to=5-1]
		\end{tikzcd}\]
		We will apply \cref{conclusion-bc} to obtain that the restriction of these squares to relative affines on the right-hand side are Beck-Chevalley in the sense recalled in \cref{def_beck_chevalley}. As such, we will first prove that $p$ and $\pi$ satisfy the conditions of this property.
	\end{rmk}
	\begin{lem} \label[lem]{props_p_pi}
		The maps $p$ and $\pi$ have the following properties:
		\begin{enumerate}
			\item $p$ is a quotient of a map between derived stacks relatively affine over $\Spec(k)$ by an affine group scheme.
			\item $\pi$ is a finite map.
		\end{enumerate}
	\end{lem}
	\begin{proof}
		\begin{enumerate}
			\item By definition, $p$ is the quotient of the following map:
			\[\mcal L^{gr}(Z/X) = \Spec_Z(\DR(Z/X)) \to Z\]
			by the affine group scheme $\bb G_m$. Moreover, $Z = \Spec(B)$ is affine and thus $\Spec_Z(\DR(Z/X))$ is relatively affine over $\Spec(k)$ too (see \cref{equiv_eps_gr_BH}).
			\item Firstly, recall from \cref{def_H} that $\mcal H$ is defined by the following split short exact sequence of abelian group schemes:
			\[1 \to B \bb G_a \to \mcal H \leftrightarrows \bb G_m \to 1\]
			We will call $\iota : \bb G_m \to \mcal H$ the section appearing in this diagram, which is a group morphism and thus admits a delooping $B(\iota) : B \bb G_m \to B \mcal H$. 
			Also note that the projection $\pi : [\mcal L^{gr}(Z/X)/\bb G_m] \to [\mcal L^{gr}(Z/X)/\mcal H]$ sits in the following pull-back square:
			\[\begin{tikzcd}[cramped]
				{[\mcal L^{gr}(Z/X)/\bb G_m]} & {B \bb G_m} \\
				{[\mcal L^{gr}(Z/X)/\mcal H]} & {B \mcal H}
				\arrow[from=1-1, to=1-2]
				\arrow["\pi"', from=1-1, to=2-1]
				\arrow["\lrcorner"{anchor=center, pos=0.125, rotate=0}, draw=none, from=1-1, to=2-2]
				\arrow["{B(\iota)}", from=1-2, to=2-2]
				\arrow[from=2-1, to=2-2]
			\end{tikzcd}\]
			Therefore, since $\pi$ is the pull-back of $B(\iota) : B \bb G_m \to B \mcal H$, it suffices to show that $B(\iota)$ is finite.
			To see that $B(\iota) : B \bb G_m \to B \mcal H$ is finite, let us first look at the embedding $\iota : \bb G_m \to \mcal H$, of which it is the delooping.
			We will show that $\iota$ has fibers equivalent to $\bb G_a$. Since it is a group morphism, it suffices to check the fiber over the unit of $\mcal H$, since all fibers are equivalent in this case. Let us thus consider the following diagram:
			\[\begin{tikzcd}[cramped]
				{\bb G_a} & \bullet & \\
				{\bb G_m} & {\bb G_m \rtimes B \bb G_a} & {\mcal H} \\
				\bullet & {B \bb G_a}
				\arrow[from=1-1, to=1-2]
				\arrow[from=1-1, to=2-1]
				\arrow["\lrcorner"{anchor=center, pos=0.125}, draw=none, from=1-1, to=2-2]
				\arrow["{1_{\mcal H}}", from=1-2, to=2-2]
				\arrow["\iota", hook, from=2-1, to=2-2]
				\arrow[two heads, from=2-1, to=3-1]
				\arrow["\lrcorner"{anchor=center, pos=0.125}, draw=none, from=2-1, to=3-2]
				\arrow[two heads, from=2-2, to=3-2]
				\arrow["{=:}"{description}, draw=none, from=2-3, to=2-2]
				\arrow["{1_{B \bb G_a}}"', from=3-1, to=3-2]
			\end{tikzcd}\]
			Note that in this diagram not all maps are group morphisms, so this diagram is considered inside derived stacks only.
			Inside this diagram:\\
			-- The projection map $\mcal H \to B \bb G_a$, which is not a group morphism, comes from the identification $\mcal H \simeq \bb G_m \times B \bb G_a$ as derived stacks (but the group structure of this derived stack is twisted by the action of $\bb G_m$, as indicated by the direct product $\rtimes$ in the definition of $\mcal H$).\\
			-- The bottom square is a pull-back since the embedding $\iota : \bb G_m \to \mcal H$ and the projection $\mcal H \to B \bb G_a$ do form a fiber sequence whose image is exactly the unit element of $B \bb G_a$.\\
			-- The outer square is a pull-back by definition of $B \bb G_a$.\\
			Therefore, by pasting, the top square is a pull-back, i.e.\ the fiber of the embedding $\iota : \bb G_m \to \mcal H$ is equivalent to $\bb G_a$. Moreover, by commutativity of the top square, the map $\bb G_a \to \bb G_m$ is constant equal to the unit of $\bb G_m$, and is thus a group morphism.
			As such, we can apply delooping, and in conclusion all the fibers of the map $B \bb G_m \to B \mcal H$ are equivalent to $B \bb G_a$, which is relatively affine and finite over $\Spec(k)$ as it is the relative spectrum of $k[\eta] = k \oplus k[-1]$ in characteristic zero, as noted in \cite[\S1.1, (A)]{HKR}. This implies that the map $B \bb G_m \to B \mcal H$ is finite.
			
			Alternatively, while the above proof only used elementary facts about $\mcal H$, note that using \cite[Rmk.~2.3.65]{robalo2024applications}, there is a canonical identification $B \mcal H = [B^2 \bb G_a/\bb G_m]$, and thus we can form the following diagram:
			\[\begin{tikzcd}
				{B \bb G_a} & \bullet & {B \bb G_m} & {[\bullet/\bb G_m]} \\
				\bullet & {B^2 \bb G_a} & {B \mcal H} & {[B^2 \bb G_a/\bb G_m]}
				\arrow[from=1-1, to=1-2]
				\arrow[from=1-1, to=2-1]
				\arrow["\lrcorner"{anchor=center, pos=0.125}, draw=none, from=1-1, to=2-2]
				\arrow[from=1-2, to=1-3]
				\arrow[from=1-2, to=2-2]
				\arrow["\lrcorner"{anchor=center, pos=0.125}, draw=none, from=1-2, to=2-3]
				\arrow[equals, from=1-3, to=1-4]
				\arrow[from=1-3, to=2-3]
				\arrow[from=1-4, to=2-4]
				\arrow[from=2-1, to=2-2]
				\arrow[from=2-2, to=2-3]
				\arrow[equals, from=2-3, to=2-4]
			\end{tikzcd}\]
			which directly proves that fibers of $B \bb G_m \to B \mcal H$ are equivalent to $B \bb G_a$. \qedhere
		\end{enumerate}
	\end{proof}
	Applying \cref{conclusion-bc} to the previous lemma, we obtain the following property.
	\begin{cor} \label[cor]{main_bc_square}
		The two inner and thus the outer square in the following diagram are Beck-Chevalley:
		\[\begin{tikzcd}
			{} & {\epsilon-cdga^{gr}(Z/X)^\op} & {\textup{RelAff}/\mcal L^{gr}(Z/X) - \mcal H} & {} \\
			& {\DR(Z/X)/cdga^{gr}(Z)^\op} & {\textup{RelAff}/\mcal L^{gr}(Z/X)-\bb G_m} \\
			{} & {cdga^{gr}(Z)^\op} & {\textup{RelAff}/Z - \bb G_m} & {}
			\arrow["{-^{gr}}"{description}, shift left=5, from=1-1, to=3-1]
			\arrow[""{name=0, anchor=center, inner sep=0}, "{\Spec^{\epsilon-gr}}"', shift right, from=1-2, to=1-3]
			\arrow["{\pi^*}", from=1-2, to=2-2]
			\arrow[""{name=1, anchor=center, inner sep=0}, "{\mcal O_{\epsilon-gr}}"', shift right, from=1-3, to=1-2]
			\arrow["{\pi^*}", from=1-3, to=2-3]
			\arrow["{-^{Z-\bb G_m}}"{description}, shift right=5, from=1-4, to=3-4]
			\arrow[""{name=2, anchor=center, inner sep=0}, "{\Spec_{\DR-gr}}"', shift right, from=2-2, to=2-3]
			\arrow["{p_*}", from=2-2, to=3-2]
			\arrow[""{name=3, anchor=center, inner sep=0}, "{\mcal O_{\DR-gr}}"', shift right, from=2-3, to=2-2]
			\arrow["{p_!}", from=2-3, to=3-3]
			\arrow[""{name=4, anchor=center, inner sep=0}, "{\Spec^{gr}}"', shift right, from=3-2, to=3-3]
			\arrow[""{name=5, anchor=center, inner sep=0}, "{\mcal O_{gr}}"', shift right, from=3-3, to=3-2]
			\arrow["\bot"{description}, draw=none, from=1, to=0]
			\arrow["\bot"{description}, draw=none, from=3, to=2]
			\arrow["\bot"{description}, draw=none, from=5, to=4]
		\end{tikzcd}\]
		In particular, the two following diagrams commute:
		\begin{enumerate}
			\item
			\[\begin{tikzcd}[cramped]
				{\epsilon-cdga^{gr}(Z/X)^\op} & {\textup{RelAff}/\mcal L^{gr}(Z/X) - \mcal H} \\
				\\
				{cdga^{gr}(Z)^\op} & {\textup{RelAff}/Z - \bb G_m}
				\arrow["{\Spec^{\epsilon-gr}}"', from=1-1, to=1-2]
				\arrow["{-^{gr}}"', from=1-1, to=3-1]
				\arrow["{-^{Z-\bb G_m}}", from=1-2, to=3-2]
				\arrow["{\Spec^{gr}}"', from=3-1, to=3-2]
			\end{tikzcd}\]
			\item
			\[\begin{tikzcd}[cramped]
				{\epsilon-cdga^{gr}(Z/X)^\op} & {\textup{RelAff}/\mcal L^{gr}(Z/X) - \mcal H} \\
				\\
				{cdga^{gr}(Z)^\op} & {\textup{RelAff}/Z - \bb G_m}
				\arrow["{-^{gr}}"', from=1-1, to=3-1]
				\arrow["{\mcal O_{\epsilon-gr}}"', from=1-2, to=1-1]
				\arrow["{-^{Z-\bb G_m}}", from=1-2, to=3-2]
				\arrow["{\mcal O_{gr}}"', from=3-2, to=3-1]
			\end{tikzcd}\]
		\end{enumerate}
	\end{cor}
	Using this corollary, we will first show that the adjunction $\mcal O_{\epsilon-gr} \dashv \Spec^{\epsilon-gr}$ restricts to $\mcal Fol(Z/X)$ on the left-hand side and $\PerfLinSt_Z/\mcal L^{gr}(Z/X) - \mcal H$ on the right-hand side.
	\begin{lem}
		Let $F \in \epsilon-cdga^{gr}(Z/X)$, and $V := \Spec^{\epsilon-gr}(F) \in \textup{RelAff}/\mcal L^{gr}(Z/X)$ $- \mcal H$.\\
		If $F \in \mcal Fol(Z/X)$, then $V \in \PerfLinSt_Z/\mcal L^{gr}(Z/X)-\mcal H$.
	\end{lem}
	\begin{proof}
		By definition, it suffices to prove that $V$ is equivalent to $\bb V(L)$ for some $L \in \Perf(Z)$ when viewed as an object in $\dSt/Z - \bb G_m$ after applying $-^{Z-\bb G_m}$. By the commutativity of square 1. in \cref{main_bc_square}, we have that viewing $V$ in this $\infty$-category is equivalent to computing $\Spec^{gr}(F^{gr})$, where $F^{gr}$ is the underlying graded algebra of $F$. Since $F$ is a derived foliation by assumption, we have that $F^{gr} \simeq \Sym_{\mcal O_Z} L$ as graded algebras for some $L \in \Perf(Z)$, so that $V \simeq \bb V(L)$ as $\bb G_m$-equivariant derived stacks over $Z$ as required.
	\end{proof}
	\begin{lem}
		Let $V \in \dSt/\mcal L^{gr}(Z/X)-\mcal H$, and let $F := \mcal O_{\epsilon-gr} V \in \epsilon-cdga^{gr}(Z/X)$.\\
		If $V \in \PerfLinSt_Z/\mcal L^{gr}(Z/X)-\mcal H$, then $F \in \mcal Fol(Z/X)$.
	\end{lem}
	\begin{proof}
		Similarly, in this case $V^{Z-\bb G_m} \simeq \bb V(L)$ for some $L \in \Perf(Z)$ by assumption, and by the commutativity of square 2. in \cref{main_bc_square}, we have that the graded algebra underlying $F$ is equivalent as a graded algebra to $\mcal O_{gr} \bb V(L)$. Since $L$ is perfect, by \cite[Prop.~2.11]{Mon}, we obtain that $\mcal O_{gr} \bb V(L) \simeq \Sym_{\mcal O_Z} L$ as graded algebras. This implies that $F$ is quasi-free, and that its weight-1 part $F(1)$ as in \cref{def_foliation}, which is thus $L$ here, is perfect. Therefore, $F$ is a derived foliation as required.
	\end{proof}
	\begin{cor}
		The adjoint functors $\mcal O_{\epsilon-gr} \dashv \Spec^{\epsilon-gr}$ restrict to an adjunction:
		\[\begin{tikzcd}
			{\mcal Fol(Z/X)} & {\PerfLinSt_Z/\mcal L^{gr}(Z/X)- \mcal H}
			\arrow[""{name=0, anchor=center, inner sep=0}, "{\Spec^{\epsilon-gr}}"', shift right, from=1-1, to=1-2]
			\arrow[""{name=1, anchor=center, inner sep=0}, "{\mcal O_{\epsilon-gr}}"', shift right, from=1-2, to=1-1]
			\arrow["\dashv"{anchor=center, rotate=-90}, draw=none, from=1, to=0]
		\end{tikzcd}\]
	\end{cor}
	We will now show that it is an equivalence of $\infty$-categories.
	\begin{lem}
		The co-unit of $\mcal O_{\epsilon-gr} \dashv \Spec^{\epsilon-gr}$ is an equivalence on objects that belong to $\mcal Fol(Z/X)$.
	\end{lem}
	\begin{proof}
		Let $F \in \mcal Fol(Z/X)$. By \cref{bc_unit} applied to the outer square of \cref{main_bc_square}, the co-unit of the adjunction $\mcal O_{\epsilon-gr} \Spec^{\epsilon-gr}(F) \leftarrow F$ is mapped by the forgetful functor $-^{gr} : \epsilon-cdga^{gr}(Z/X) \to cdga^{gr}(Z)$ to the co-unit of the adjunction $\mcal O_{gr} \Spec^{gr}(F) \leftarrow F$. By \cite[Cor.~2.9]{Mon}, this map is an equivalence whenever $F$ is quasi-free and $F(1) \in \QCoh^-(Z)$, i.e.\ $F(1)$ is left-bounded in amplitude. Since here $F$ is a derived foliation, $F(1)$ is perfect and thus left-bounded, so this co-unit is an equivalence. Finally, since $-^{gr}$ is a conservative functor, the initial co-unit is an equivalence too.
	\end{proof}
	\begin{lem}
		The unit of $\mcal O_{\epsilon-gr} \dashv \Spec^{\epsilon-gr}$ is an equivalence on objects that belong to $\PerfLinSt_Z/\mcal L^{gr}(Z/X) - \mcal H$.
	\end{lem}
	\begin{proof}
		Similarly, for $V = \bb V(L) \in \PerfLinSt_Z/\mcal L^{gr}(Z/X) - \mcal H$, the unit $V \to \Spec^{\epsilon-gr}(\mcal O_{\epsilon-gr} V)$ is sent by the forgetful functor $\PerfLinSt_Z/\mcal L^{gr}(Z/X)-\mcal H \to \dSt/Z - \bb G_m$ to the unit $V \to \Spec^{gr}(\mcal O_{gr} V)$. By \cite[Prop.~2.11]{Mon}, since $L$ is perfect, this last map is an equivalence, but the forgetful functor is conservative, so the former unit also is an equivalence.
	\end{proof}
	\begin{thm} \label[thm]{equiv_fol_dst}
		The following adjunction is an equivalence of $\infty$-categories:
		\[\begin{tikzcd}
			{\mcal Fol(Z/X)} & {\PerfLinSt_Z/\mcal L^{gr}(Z/X)- \mcal H}
			\arrow[""{name=0, anchor=center, inner sep=0}, "{\Spec^{\epsilon-gr}}"', shift right, from=1-1, to=1-2]
			\arrow[""{name=1, anchor=center, inner sep=0}, "{\mcal O_{\epsilon-gr}}"', shift right, from=1-2, to=1-1]
			\arrow["\dashv"{anchor=center, rotate=-90}, draw=none, from=1, to=0]
		\end{tikzcd}\]
	\end{thm}
	\subsection{Equivalence in the Global Case} \label{equiv_global}
	To finish this section, after dealing with the affine case, we will now globalize these results by inspecting how \cref{equiv_fol_dst} behaves under gluing. This will yield a slightly different equivalence in the general case, but we will show that it reduces to the same equivalence when $Z \to X$ is a relative derived Deligne-Mumford stack.
	\subsubsection{Structures in the Global Case}
	Let us start by defining structures from the previous subsection in the global case. This will require first inspecting the functoriality of these categories on affines. Note that by \cref{formula_pull-back}, the $\infty$-category $\mcal Fol$ of perfect foliations is not functorial in general, but it is functorial when restricted to affines $\Spec(B) \to \Spec(A)$ which are of finite presentation. This motivates the following remark.
	\begin{rmk}
		Let $p : Z \to X$ be an arbitrary map of derived stacks. We denote by $\Arr(\dAff)/p$ the $\infty$-category of commutative squares:
		\[\begin{tikzcd}
			V & Z \\
			U & X
			\arrow[from=1-1, to=1-2]
			\arrow[from=1-1, to=2-1]
			\arrow[from=1-2, to=2-2]
			\arrow[from=2-1, to=2-2]
		\end{tikzcd}\]
		where $U$ and $V$ are derived affines. This is the canonical generating site of the $\infty$-topos $\Arr(\dSt)$ as in \cite[Def.~6.2.4.1, Prop.~6.2.4.2]{HTT}, and as such the colimit of the objects $V \to U$ in this category is always equivalent to $p$.
		
		We now define $\Arr(\dAff)^{fp}/p$ to be the full subcategory of $\Arr(\dAff)/p$ spanned by the objects as above such that the map $V \to U$ is of finite presentation. In general, it is no longer the case that colimits of such objects recover $p$, but we will show that it is the case when $p$ is locally of finite presentation.
	\end{rmk}
	\begin{prop} \label[prop]{daff_fp_cofinal}
		Let $p : Z \to X$ be a map of derived stacks which is locally of finite presentation in the sense of \cite[Def.~2.3.1]{lurieDAGXIV}. Then $\Arr(\dAff)^{fp}/p$ is cofinal in $\Arr(\dAff)/p$, and thus $f \simeq \underset{(V \to U) \in \Arr(\dAff)^{fp}/p}{\textup{colim}}(V \to U)$.
	\end{prop}
	\begin{proof}
		By the cofinality criterion \cite[Thm.~4.1.3.1]{HTT}, it suffices to show that, for any object $q \in \Arr(\dAff)/p$, the $\infty$-category $q/(\Arr(\dAff)/p)$ is weakly contractible, in the sense that its $\infty$-groupoid completion is a contractible space. Let us thus fix a square $q=(\Spec(B)\to\Spec(A))$ over $p$, and set
		\[
		\mcal C=q/(\Arr(\dAff)^{fp}/p)
		\]
		which we must show to be weakly contractible in this sense.
		
		By definition, an object of $\mcal C$ is the data of a factorization
		\[\begin{tikzcd}[ampersand replacement=\&,cramped]
			{\Spec(B)} \& {\Spec(B')} \& Z \\
			{\Spec(A)} \& {\Spec(A')} \& X
			\arrow[from=1-1, to=1-2]
			\arrow[from=1-1, to=2-1]
			\arrow[from=1-2, to=1-3]
			\arrow[from=1-2, to=2-2]
			\arrow[from=1-3, to=2-3]
			\arrow[from=2-1, to=2-2]
			\arrow[from=2-2, to=2-3]
		\end{tikzcd}\]
		where $\Spec(B') \to \Spec(A')$ is of finite presentation.
		Let $\mcal C_A \subset \mcal C$ denote the full subcategory of those factorizations such that $\Spec(A') = \Spec(A)$ and the map $\Spec(A) \to \Spec(A')$ is the identity. Its inclusion has a right adjoint, which maps $(\Spec(B') \to \Spec(A'))$ to $(\Spec(B' \otimes_{A'} A) \to \Spec(A))$. Thus $\mcal C_A$ and $\mcal C$ are weakly equivalent, and it suffices to show that $\mcal C_A$ is weakly contractible.
		
		Let now $Y = Z \times_X \Spec(A)$, and notice that the commutativity of the square $q$ induces a map $\Spec(B) \to Y$, i.e. an element $z \in Y(B)$.
		Now note that $B$ is canonically the colimit of the finitely presented $A$-algebras mapping to $B$, and that their diagram $I$ is a filtered diagram, so that we have $B \simeq \underset{i \in I}{\textup{colim}}\ B_i$.
		Since $p$ is locally of finite presentation, by definition its pull-back $Y$ is locally of finite presentation over $\Spec(A)$, and as such it preserves filtered
		colimits of $A$-algebras by definition, so we have:
		\[
		\underset{i \in I}{\textup{colim}}\ Y(B_i) \simeq Y(B)
		\]
		Now notice that by unfolding the definitions $\mcal C_A^\op$ is exactly the category of elements of the functor	$i \mapsto Y(B_i) \times_{Y(B)} \{z\}$.	By \cite[Cor.~3.3.4.6]{HTT}, its groupoid completion can therefore be computed by:
		\begin{align*}
			|\mcal C_A^\op| &\simeq\ \underset{i \in I}{\textup{colim}}\ \bigl(Y(B_i)\times_{Y(B)}\{z\}\bigr) \\
			&\simeq\ \left(\underset{i \in I}{\textup{colim}}\ Y(B_i)\right)\times_{Y(B)}\{z\} \\
			&\simeq\ Y(B) \times_{Y(B)} \{z\} \\
			&\simeq\ \{z\}
		\end{align*}
		where the second equivalence holds because filtered colimits of spaces commute with finite limits. Hence $\mcal C_A^\op$ and thus $\mcal C_A$ and $\mcal C$ are weakly contractible, as required.
	\end{proof}
	\begin{rmk}
		In fact, we expect the above proposition to be an equivalence, i.e. a map $p$ is locally of finite presentation if and only if $\Arr(\dAff)^{fp}/p$ is cofinal in $\Arr(\dAff)/p$, but we will not need this fact here.
	\end{rmk}
	\begin{df} \label[df]{def_global_structures}
		Let $p : Z \to X$ be an arbitrary map of derived stacks.
		\begin{enumerate}
			\item By \cref{pull-back}, we defined contravariant functors $\epsilon-dg^{gr},\epsilon-cdga^{gr} : \Arr(\dAff)^\op \to \infty-Cat$. By \cref{formula_pull-back}, we can restrict $\epsilon-cdga^{gr}$ to a contravariant functor $\mcal Fol : \Arr(\dAff)^{fp,\op} \to \infty-Cat$ on the finitely-presented objects. Similarly, by pull-back of graded algebras, we get contravariant functors $dg^{gr},cdga^{gr} : \Arr(\dAff)^\op \to \infty-Cat$ which map $\Spec(B) \to \Spec(A)$ to $dg^{gr}(B),cdga^{gr}(B)$.
			\item Analogously to \cite[Def.~2.1.3]{toen2025derivedfoliations}, when $p : Z \to X$ is moreover locally of finite presentation\footnote{We expect the definition of $\mcal Fol$ to extend beyond the locally finitely presented case.}, we set:
			\[\mcal Fol(Z/X) := \underset{(V \to U) \in \left(\Arr(\dAff)^{fp}/p\right)^\op}{\lim} \mcal Fol(V/U)\]
			\item Similarly, we set:
			\begin{align*}
				\epsilon-dg^{gr}(Z/X) &:= \underset{(V \to U) \in (Arr(\dAff)/p)^\op}{\lim} \epsilon-dg^{gr}(V/U)\\
				\epsilon-cdga^{gr}(Z/X) &:= \underset{(V \to U) \in (Arr(\dAff)/p)^\op}{\lim} \epsilon-cdga^{gr}(V/U)
			\end{align*}
			\item Finally, we set:
			\begin{align*}
				dg^{gr}(Z) &:= \underset{V \in (\dAff/Z)^\op}{\lim} dg^{gr}(V)\\
				cdga^{gr}(Z) &:= \underset{V \in (\dAff/Z)^\op}{\lim} cdga^{gr}(V)
			\end{align*}
		\end{enumerate}
	\end{df}
	\begin{rmk}
		Note that by \cref{daff_fp_cofinal}, in the above definitions, when $p : Z \to X$ is locally of finite presentation, all instances of $Arr(\dAff)/p$ can be replaced by $Arr(\dAff)^{fp}/p$, as they compute the same limits.
		
		Also note that if $Z$ and $X$ are affines, the above definitions do recover the definitions in the affine case, because in this case the above limits admit an initial object given by taking $(V \to U) = (Z \to X)$.
		
		Moreover, it is important to note that while objects of $dg^{gr},cdga^{gr}$ are sheaves of graded objects in the usual sense, it is not the case for $\epsilon-dg^{gr}$ and $\epsilon-cdga^{gr}$ because the functoriality of the limit defining them is not the usual one (the pull-back is given by \cref{pull-back}, which is more complicated than usual pull-backs of sheaves). This means that at this stage defining an object of $\mcal Fol(Z/X)$ is quite difficult. We will later on establish more simple descriptions of these categories.
	\end{rmk}
	\begin{rmk} \label[rmk]{rmk_global_structures}
		In this remark, let $p : Z \to X$ be an arbitrary map of derived stacks, not necessarily locally of finite presentation. Consider the definition of the above point 3. For $(V \to U) \in \Arr(\dAff)/p$, by \cref{equiv_eps_gr_BH} we have:
		\[\epsilon-dg^{gr}(V/U) \simeq \QCoh([\mcal L^{gr}(V/U)/\mcal H])\]
		As a consequence, since $\Arr(\dAff)/p$ covers $Z \to X$ by \cite[Def.~6.2.4.1, Prop.~6.2.4.2]{HTT}, and $\QCoh$ sends colimits to limits by \cite[Prop.~2.7.14]{lurieDAGVIII} and \cite[Prop.~5.5.4.20]{HTT}, we obtain that:
		\[\epsilon-dg^{gr}(Z/X) \simeq \QCoh(\underset{(V \to U) \in \Arr(\dAff)/p}{\textup{colim}} [\mcal L^{gr}(V/U)/\mcal H])\]
		Moreover, the derived stack $\underset{(V \to U) \in \Arr(\dAff)/p}{\textup{colim}} \mcal L^{gr}(V/U)$ inherits an $\mcal H$-action, and taking the quotient by $\mcal H$ commutes with colimits (as the quotient by $\mcal H$ is itself a colimit), so we obtain that:
		\[\epsilon-dg^{gr}(Z/X) \simeq \QCoh([\underset{(V \to U) \in \Arr(\dAff)/p}{\textup{colim}} \mcal L^{gr}(V/U)/\mcal H])\]
		This motivates the next definition.
		
		Also note that a similar computation in the graded case directly yields the known equivalence:
		\[dg^{gr}(Z) \simeq \QCoh(Z \times B \bb G_m)\]
	\end{rmk}
	\begin{df} \label[df]{df-lgrhat}
		Let $\widetilde{\mcal L^{gr}}(Z/X) := \underset{(V \to U) \in \Arr(\dAff)/p}{\textup{colim}} \mcal L^{gr}(V/U)$.
		It comes equipped with a canonical $\mcal H$-action, and a canonical $\mcal H$-equivariant map $\widetilde{\mcal L^{gr}}(Z/X) \to \mcal L^{gr}(Z/X)$ over $Z$, which we do not expect to be an equivalence in general.
	\end{df}
	\begin{rmk}
		In \cite[Def.~1.22]{loopspacesandconnections}, another object called $\widehat{\mcal L}$ is defined as the formal completion of the loop stack at the constant loops. Note that this is a distinct object, although we expect them to coincide at least in some cases, but we will not need this fact here.
	\end{rmk}
	\begin{cor} \label[cor]{global_structures}
		Let $Z \to X$ be an arbitrary map of derived stacks. By \cref{rmk_global_structures}, we have the following equivalences:
		\begin{align*}
			\epsilon-dg^{gr}(Z/X) &\simeq \QCoh([\widetilde{\mcal L^{gr}}(Z/X)/\mcal H])\\
			\epsilon-cdga^{gr}(Z/X) &\simeq \CAlg(\QCoh([\widetilde{\mcal L^{gr}}(Z/X)/\mcal H]))\\
			dg^{gr}(Z) &\simeq \QCoh(Z \times B \bb G_m)\\
			cdga^{gr}(Z) &\simeq \CAlg(\QCoh(Z \times B \bb G_m))
		\end{align*}
	\end{cor}
	\begin{rmk} \label[rmk]{gr-stack}
		Note that it is immediate to see that the functor $(V \to U) \mapsto V \times B \bb G_m : \Arr(\dAff) \to \dSt$ satisfies étale descent, and that its left Kan extension is simply given by $(Z \to X) \mapsto Z \times B \bb G_m$. Since $\QCoh$ and $\CAlg$ preserve étale descent, we also obtain that $dg^{gr},cdga^{gr} : \Arr(\dAff)^\op \to \infty-Cat$ are stacks for the étale topology. In particular, we obtain that $dg^{gr}$ and $cdga^{gr}$ from \cref{def_global_structures} send arbitrary colimits to limits.
		
		Similarly, to obtain that $\epsilon-dg^{gr},\epsilon-cdga^{gr} : \Arr(\dAff)^\op \to \infty-Cat$ are stacks for the étale topology, and that they send arbitrary colimits to limits, it suffices to prove that the functor $(V \to U) \mapsto \mcal L^{gr}(V/U) : \Arr(\dAff) \to \dSt$ satisfies étale codescent. This is the content of the following proposition.
	\end{rmk}
	\begin{prop} \label[prop]{rmk_lgrhat_stack}
		The functor $(V \to U) \mapsto \mcal L^{gr}(V/U) : \Arr(\dAff) \to \dSt$ satisfies étale codescent.
	\end{prop}
	\begin{proof}
		Let
		\[
		u:\Spec(B)\to \Spec(A)
		\]
		be an object of $\Arr(\dAff)$, and let
		\[
		u_\bullet:\Spec(B_\bullet)\to \Spec(A_\bullet)
		\]
		be an étale hypercover of $u$ in $\Arr(\dAff)$. Thus $u_\bullet$ is an augmented simplicial object over $u$, and the augmentations
		\[
		\Spec(B_\bullet)\to \Spec(B),
		\qquad
		\Spec(A_\bullet)\to \Spec(A)
		\]
		are étale hypercovers.
		
		In particular, for every $n\ge 0$, the morphisms
		\[
		\Spec(B_n)\to \Spec(B),
		\qquad
		\Spec(A_n)\to \Spec(A)
		\]
		are étale, $\Spec(B) \simeq \textup{colim} \; \Spec(B_\bullet)$ and $\Spec(A) \simeq \textup{colim} \; \Spec(A_\bullet)$. Let $n \in \bb N$. Since the maps $A \to A_n$ and $B \to B_n$ are étale, we have:
		\[
		\mathbb L_{B_n/B}\simeq 0,
		\qquad
		\mathbb L_{A_n/A}\simeq 0.
		\]
		Applying transitivity of the cotangent to the composable maps $A\to B\to B_n$, one gets a triangle
		\[
		B_n\otimes_B \mathbb L_{B/A} \to \mathbb L_{B_n/A} \to \bb L_{B_n/B}
		\]
		Since $\bb L_{B_n/B} \simeq 0$, this implies that $B_n\otimes_B \mathbb L_{B/A} \xrightarrow{\sim} \mathbb L_{B_n/A}$.
		Similarly, applying the same argument to $A\to A_n\to B_n$, one gets
		\[
		\mathbb L_{B_n/A}\xrightarrow{\sim}\mathbb L_{B_n/A_n}.
		\]
		Therefore
		\[
		\mathbb L_{B_n/A_n}\simeq B_n\otimes_B \mathbb L_{B/A}.
		\]
		It follows that
		\[
		\Sym_{B_n}\bigl(\mathbb L_{B_n/A_n}[1]\bigr)
		\simeq
		\Sym_{B_n}\bigl((B_n\otimes_B \mathbb L_{B/A})[1]\bigr)
		\simeq
		B_n\otimes_B \Sym_B\bigl(\mathbb L_{B/A}[1]\bigr).
		\]
		By compatibility of relative spectrum with base change, this yields
		\[
		\mathcal L^{gr}(u_n)
		\simeq
		\mathcal L^{gr}(u)\times_{\Spec(B)} \Spec(B_n)
		\]
		functorially in $n$. Passing to the colimit, we obtain:
		\[\textup{colim} \; \mcal L^{gr}(u_\bullet) \simeq \textup{colim} \; \mcal L^{gr}(u) \times_{\Spec(B)} \Spec(B_\bullet)\]
		Since colimits commute with pull-backs, we obtain:
		\[\textup{colim} \; \mcal L^{gr}(u_\bullet) \simeq \mcal L^{gr}(u) \times_{\Spec(B)} \Spec(B) \simeq \mcal L^{gr}(u)\]
		Thus $\mathcal L^{gr}$ satisfies étale codescent.
	\end{proof}
	\begin{rmk} \label[cor]{lgrhat_descent}
		Note that $\widetilde{\mcal L^{gr}} : \Arr(\dSt) \to \dSt$ is by definition the left Kan extension of $\mcal L^{gr} : \Arr(\dAff) \to \dSt$. As such, for affines $V \to U$, the canonical map $\widetilde{\mcal L^{gr}}(V/U) \to \mcal L^{gr}(V/U)$ is an equivalence. Over a fixed affine base $U$, this also extends immediately to a disjoint union of affines because mapping stacks commute with disjoint union in their target, therefore: $\mcal L^{gr}(\coprod_\bullet V_\bullet/U) = \underline{\Map}_U(B \bb G_a \times U,\coprod_\bullet V_\bullet) \simeq \coprod_\bullet  \underline{\Map}_U(B \bb G_a \times U,V_\bullet) \simeq \coprod_\bullet \mcal L^{gr}(V_\bullet/U)$. Also, since by construction this is a map of derived stacks over $V$, and it is $\mcal H$-equivariant, these derived stacks are equivalent as derived stacks over $V$ and as $\mcal H$-equivariant derived stacks over $U$.
		
		Additionally, as a left Kan extension of a functor satisfying étale codescent on affines, it commutes with all colimits in $\Arr(\dSt)$. Equivalently, for any diagram $(Z_i \to X_i)$ in $\Arr(\dSt)$ with colimit $Z \to X$, one has $\widetilde{\mcal L^{gr}}(Z/X) \simeq \underset{i}{\textup{colim}} \; \widetilde{\mcal L^{gr}}(Z_i/X_i)$.
	\end{rmk}
	\begin{rmk} 
		Also note, as in \cref{gr-stack}, that the above proposition also implies that $\epsilon-dg^{gr},\epsilon-cdga^{gr} : \Arr(\dAff)^\op \to \infty-Cat$ are stacks for the étale topology.
		
		Moreover, when $p : Z \to X$ is a map of derived stacks locally of finite presentation, the fact that $\Arr(\dAff)^{fp}/p$ is cofinal in $\Arr(\dAff)/p$ by \cref{daff_fp_cofinal} implies that $\epsilon-cdga^{gr}(Z/X) \simeq \underset{(V \to U) \in (\Arr(\dAff)^{fp}/p)^\op}{\textup{lim}} \epsilon-cdga^{gr}(V/U)$, and the fact that $\mcal Fol$ is stable by the functoriality of this diagram by \cref{formula_pull-back} also implies that $\mcal Fol : (\Arr(\dAff)^{fp}/p)^\op \to \infty-Cat$ satisfies descent for the étale topology restricted to this $\infty$-category.
	\end{rmk}
	\subsubsection{Equivalence for Maps Locally of Finite Presentation}
	To generalize the equivalence of the affine case, we will start by proving that it is stable under affine base-change.
	\begin{prop} \label[prop]{equiv_functorial}
		Let $Z \to X \in \Arr(\dAff)^{fp}$. The equivalence of $\infty$-categories $\mcal Fol(Z/X) \underset{\sim}{\xrightarrow{\Spec^{\epsilon-gr}}} \PerfLinSt_Z/\mcal L^{gr}(Z/X) - \mcal H$ provided by \cref{equiv_fol_dst} is contravariant in $Z \to X$, i.e.\ it defines a functor $\textup{Arr}(\dAff)^{fp,\op} \to \textup{Equiv}(\infty-Cat)$.
	\end{prop}
	\begin{proof}
		We have already described the functoriality of $\mcal Fol(Z/X)$ in \cref{def_global_structures}.  We will now describe the functoriality of $(Z \to X) \mapsto \PerfLinSt_Z/\mcal L^{gr}(Z/X) - \mcal H$.
		Consider a commutative square where each object is affine, and both $Z \to X$ and $Z' \to X'$ are of finite presentation:
		\[\begin{tikzcd}
			{Z'} & Z \\
			{X'} & X
			\arrow[from=1-2, to=2-2]
			\arrow[from=1-1, to=2-1]
			\arrow[from=1-1, to=1-2]
			\arrow[from=2-1, to=2-2]
		\end{tikzcd}\]
		By functoriality of $\mcal L^{gr}$, it induces a canonical map $p : \mcal L^{gr}(Z'/X') \to \mcal L^{gr}(Z/X)$ of derived stacks over $Z$.
		We will describe the pull-back functor $p^*$ on the side of derived stacks. Let $\bb V(E[1]) \in \PerfLinSt_Z/\mcal L^{gr}(Z/X) - \mcal H$. Since all the objects are affines, we can form the following pull-back diagram:
		\[\begin{tikzcd}[cramped,column sep=small]
			{\Spec(\DR(Z'/X') \otimes_{\DR(Z/X)} \Sym_{\mcal O_Z} E[1])} && {\Spec(\Sym_{\mcal O_Z} E[1])} \\
			{\mathcal L^{gr}(Z'/X') = \Spec(\DR(Z'/X'))} && {\textup{Spec}(\DR(Z/X)) = \mathcal L^{gr}(Z/X)} \\
			& Z
			\arrow[from=1-1, to=1-3]
			\arrow[from=1-1, to=2-1]
			\arrow["\lrcorner"{anchor=center, pos=0.125, rotate=0}, draw=none, from=1-1, to=2-3]
			\arrow[from=1-3, to=2-3]
			\arrow["p"', from=2-1, to=2-3]
			\arrow[from=2-1, to=3-2]
			\arrow[from=2-3, to=3-2]
		\end{tikzcd}\]
		Let $Z = \Spec(B)$, $Z' = \Spec(B')$, $X = \Spec(A)$ and $X' = \Spec(A')$, and note that $\DR(Z/X) = \Sym_{B}(\bb L_{B/A}[1])$ and $\DR(Z'/X') = \Sym_{B'}(\bb L_{B'/A'}[1])$. As such, the homotopy fiber product in the top left corner can be computed by derived tensor products, which yields $\Sym_{B'}\left(\left( \bb L_{B'/A'} \oplus_{\bb L_{B/A} \otimes_B B'} (E \otimes_B B') \right) [1]\right)$. In particular, since $\Spec(B') \to \Spec(A')$ and $\Spec(B) \to \Spec(A)$ are of finite presentation by assumption, then $\bb L_{B'/A'}$ and $\bb L_{B/A}$ are perfect. Moreover, $E$ is perfect by assumption, and $- \otimes_B B'$ preserves perfects. Therefore, we do obtain that $p^* \bb V(E[1])$ is a perfect linear stack over $Z'$, equipped with a structural map to $\mcal L^{gr}(Z'/X')$, and $p^*$ does induce an $\mcal H$-action on it, so $p^* \bb V(E[1]) \in \PerfLinSt_{Z'}/\mcal L^{gr}(Z'/X') - \mcal H$. Finally, by \cref{base_chg_adj}, we obtain that the following square is canonically commutative:
		\[\begin{tikzcd}
			{\mcal Fol(Z/X)} & {\PerfLinSt_{Z}/\mcal L^{gr}(Z/X) - \mcal H} \\
			{\mcal Fol(Z'/X')} & {\PerfLinSt_{Z'}/\mcal L^{gr}(Z'/X') - \mcal H}
			\arrow["{p^*}"', from=1-1, to=2-1]
			\arrow["{\Spec^{\epsilon-gr}_{Z'/X'}}"', from=2-1, to=2-2]
			\arrow["{p^*}", from=1-2, to=2-2]
			\arrow["{\Spec^{\epsilon-gr}_{Z/X}}"', from=1-1, to=1-2]
			\arrow["\simeq", from=1-1, to=1-2]
			\arrow["\simeq", from=2-1, to=2-2]
		\end{tikzcd}\]
		Which precisely means that the map:
		\[(Z \to X) \mapsto \left(\mcal Fol(Z/X) \simeq \PerfLinSt_Z/\mcal L^{gr}(Z/X) - \mcal H \right)\]
		is contravariant in $Z \to X$.
	\end{proof}
	\begin{rmk}
		Consider now $p : Z \to X$ a map of derived stacks which is locally of finite presentation.
		Starting from the equivalence of \cref{equiv_fol_dst}, and passing through the limit over $\Arr(\dAff)^{fp}/p$ using the functoriality described in the above corollary, we always obtain an equivalence of $\infty$-categories:
		\[\underset{(V \to U) \in (\Arr(\dAff)^{fp}/p)^\op}{\textup{lim}} \mcal Fol(V/U) \simeq \underset{(V \to U) \in (\Arr(\dAff)^{fp}/p)^\op}{\textup{lim}} \PerfLinSt_V/\mcal L^{gr}(V/U) - \mcal H\]
		Moreover, on the left-hand side, by definition 2.\ in \cref{def_global_structures}, it is always the case that:
		\[\mcal Fol(Z/X) := \underset{(V \to U) \in (\Arr(\dAff)^{fp}/p)^\op}{\textup{lim}} \mcal Fol(V/U)\]
		However, on the right-hand side, it is not always the case that:
		\[\PerfLinSt_Z/\mcal L^{gr}(Z/X) - \mcal H \overset{?}{=} \underset{(V \to U) \in (\Arr(\dAff)^{fp}/p)^\op}{\textup{lim}} \PerfLinSt_V/\mcal L^{gr}(V/U) - \mcal H\]
		Instead, we will describe the right-hand limit in the next definition.
	\end{rmk}
	\begin{df}
		Let $p : Z \to X$ be a map of derived stacks which is lfp. We define:
		\[\PerfLieSt/\widetilde{\mcal L^{gr}}(Z/X) - \mcal H \subset \dSt/\widetilde{\mcal L^{gr}}(Z/X) - \mcal H\]
		to be the full subcategory of objects $L \in \dSt/\widetilde{\mcal L^{gr}}(Z/X) - \mcal H$ such that, for all $V \to U$ affines over $Z \to X$ with $V \to U$ of finite presentation, the pull-back $L \times_{\widetilde{\mcal L^{gr}}(Z/X)} \mcal L^{gr}(V/U)$ is a perfect linear stack over $V$. We will call \emph{perfect Lie stacks} objects satisfying this condition.
	\end{df}
	\begin{rmk} \label[rmk]{perf_lie_st_rmk}
		In the above situation, let $L$ be a perfect Lie stack and consider a derived affine $\Spec(C) \to \widetilde{\mcal L^{gr}}(Z/X)$, which induces a commutative square $(\Spec(C) \to \Spec(C)) \to (Z \to X)$. Since we have a canonical identification $\mcal L^{gr}(C/C) \simeq \Spec(C)$ because $\bb L_{C/C} \simeq 0$, the fact that $L \times_{\widetilde{\mcal L^{gr}}(Z/X)} \mcal L^{gr}(C/C)$ is perfect linear over $\Spec(C)$ implies that $L$ is relatively affine over $\widetilde{\mcal L^{gr}}(Z/X)$. In the affine case, $\widetilde{\mcal L^{gr}}(Z/X) \simeq \mcal L^{gr}(Z/X) \simeq \Spec_Z(\DR(Z/X))$ is itself relatively affine over $Z$, and the definition of perfect Lie stacks simply recovers perfect linear stacks. The proof of this last point will be made in \cref{equiv_fol_dst_rel_dm}.
		
		In fact, the name perfect Lie stack has a more precise interpretation in the affine case, because then perfect Lie algebroids are equivalent to perfect derived foliations by the Chevalley-Eilenberg construction (see the main result of \cite{FuDuality} for the more general equivalence, and \cite[Lem.~3.2.19,3.2.26]{GrataloupThesis} for a description of the functoriality of this equivalence with respect to the base which recovers the one described in this paper). Under this equivalence, by our results in the affine case, perfect Lie algebras are covariantly equivalent to perfect Lie stacks (which in this case coincide with perfect linear stacks), by mapping a perfect Lie algebroid with tangent $\bb T_L$ to its total vector bundle $\bb V(\bb T_L^\vee)$ with extra structure, whose sections are exactly elements of $\bb T_L$ up to homotopy. Beyond the affine case, we expect this equivalence to hold at least in the derived Deligne-Mumford case by étale descent, and even more generally for arbitrary derived pre-stacks up to a choice of convention and assumptions on perfect Lie algebroid. In particular, we do not know if the global definition of Lie algebroids in the general case given by \cite[\S2.1]{gaitsgoryRozenblyumII} is equivalent to perfect Lie stacks in our setting, as this would require a formal integration theorem similar to the main result of \cite{brantner2025formalintegrationderivedfoliations} but with slightly different assumptions and definitions. However, we will not study this question in this paper.
	\end{rmk}
	\begin{prop}
		Let $p : Z \to X$ be a map of derived stacks which is locally of finite presentation. There is a canonical equivalence of $\infty$-categories:
		\[\underset{(V \to U) \in (\Arr(\dAff)^{fp}/p)^\op}{\textup{lim}} \PerfLinSt_V/\mcal L^{gr}(V/U) - \mcal H = \PerfLieSt/\widetilde{\mcal L^{gr}}(Z/X) - \mcal H\]
	\end{prop}
	\begin{proof}
		By definition of $\widetilde{\mcal L^{gr}}(Z/X)$ and its $\mcal H$-action, since $\Arr(\dAff)^{fp}/p$ is cofinal in $\Arr(\dAff)/p$ by \cref{daff_fp_cofinal}, and since $T \mapsto \dSt/T$ sends colimits to limits, we have:
		\[\underset{(V \to U) \in (\Arr(\dAff)^{fp}/p)^\op}{\textup{lim}} \dSt/\mcal L^{gr}(V/U) - \mcal H \simeq \dSt/\widetilde{\mcal L^{gr}}(Z/X) - \mcal H\]
		Moreover, since for each $(V \to U) \in \Aff(\dAff)^{fp}/p$, $\PerfLinSt_V/\mcal L^{gr}(V/U) - \mcal H$ is a full subcategory of $\dSt/\mcal L^{gr}(V/U) - \mcal H$, we also obtain that the induced functor:
		\[\underset{(V \to U) \in (\Arr(\dAff)^{fp}/p)^\op}{\textup{lim}} \PerfLinSt_V/\mcal L^{gr}(V/U) - \mcal H \to \dSt/\widetilde{\mcal L^{gr}}(Z/X) - \mcal H\]
		is fully faithful. Note that this functor is defined by mapping a family of objects on the left to its colimit in derived stacks, which lands on the right, and the resulting object $L \in \dSt/\widetilde{\mcal L^{gr}}(Z/X) - \mcal H$ satisfies that for any $(V \to U) \in \Arr(\dAff)^{fp}/p$, the pull-back of $L$ over $\mcal L^{gr}(V/U) \to \widetilde{\mcal L^{gr}}(Z/X)$ is in $\PerfLinSt_V/\mcal L^{gr}(V/U) - \mcal H$. As a result, by definition, an object $L \in \dSt/\widetilde{\mcal L^{gr}}(Z/X) - \mcal H$ belongs to the essential image of the above functor if and only if it satisfies the condition of being a perfect Lie stack as defined above.
	\end{proof}
	Therefore, one always has the following equivalence.
	\begin{cor} \label[cor]{general_foliations}
		For $Z \to X$ a map of derived stacks which is locally of finite presentation, there is a canonical equivalence of $\infty$-categories:
		\[\mcal Fol(Z/X) \simeq \PerfLieSt/\widetilde{\mcal L^{gr}}(Z/X) - \mcal H\]
		Explicitly, this functor maps the graded-mixed algebra $\DR(\mcal F)$ underlying a derived foliation $\mcal F$, viewed as a commutative algebra in $\QCoh([\widetilde{\mcal L^{gr}}(Z/X)/\mcal H])$, to its relative spectrum over $\widetilde{\mcal L^{gr}}(Z/X)$.
	\end{cor}
	\begin{rmk}
		In the above situation, let $q : \widetilde{\mcal L^{gr}}(Z/X) \to Z$ be the canonical map, and note that it is $\bb G_m$-equivariant for the trivial $\bb G_m$-action on $Z$ and the induced $\bb G_m$-action of the structural $\mcal H$-action on $\widetilde{\mcal L^{gr}}(Z/X)$. Moreover, given a perfect Lie stack $\ell : L \to \widetilde{\mcal L^{gr}}(Z/X)$, the different structures of the derived foliation associated to $L$ can be recovered in the following ways:
		\begin{enumerate}
			\item Its underlying graded-mixed algebra is:
			\begin{align*}
				\DR(L)^{\epsilon-gr} := [\ell/\mcal H]_*(\mcal O_{[L/\mcal H]}) &\in \CAlg(\QCoh([\widetilde{\mcal L^{gr}}(Z/X)/\mcal H]))\\ 
				&= \epsilon-cdga^{gr}(Z/X)
			\end{align*}
			Note that since $L$ is a perfect Lie stack, by \cref{perf_lie_st_rmk} it is relatively affine over $\widetilde{\mcal L^{gr}}(Z/X)$ and thus the relative spectrum of $\DR(L)^{\epsilon-gr}$ reconstructs $L$ as a perfect Lie stack with its $\mcal H$-action and its structural $\mcal H$-equivariant map to $\widetilde{\mcal L^{gr}}(Z/X)$.
			\item Its underlying graded-algebra is:
			\[\DR(L)^{gr} := [\ell \circ q/\bb G_m]_*(\mcal O_{[L/\bb G_m]}) \in \CAlg(\QCoh(Z \times B \bb G_m)) = cdga^{gr}(Z)\]
			\item Its cotangent complex $\bb L_L \in \QCoh(Z)$ is obtained from the previous item by applying the forgetful functor $-(1) : cdga^{gr}(Z) \to \QCoh(Z)$. Note that if we denote by $\pi : Z \times B \bb G_m \to Z$ the projection, this functor can also be described geometrically as:
			\[\bb L_L := \pi_*(\DR(L)^{gr} \otimes \mcal O(-1)) \in \QCoh(Z)\]
			where $\mcal O(-1)$ is inverse of the canonical line bundle on $B \bb G_m$, concentrated in weight -1.
		\end{enumerate}
		Also note that $\mcal Fol(Z/X)$ always admits a terminal derived foliation $\star_{Z/X}$, which as a perfect Lie stack is simply the identity map $\widetilde{\mcal L^{gr}}(Z/X) \to \widetilde{\mcal L^{gr}}(Z/X)$. Moreover, when $p : Z \to X$ admits a relative cotangent, we can also obtain an explicit description of $\DR(\widetilde{\mcal L^{gr}}(Z/X))^{gr}$.
	\end{rmk}
	\begin{prop} \label[prop]{Ogr_L_gr}
		Let $p : Z \to X$ be any map of derived stacks which is locally of finite presentation and assume that it admits a relative cotangent $\bb L_{Z/X}$, which is thus perfect by local finite presentation. Then, there is a canonical equivalence $\DR(\star_{Z/X})^{gr} \simeq \DR(Z/X)$ as graded algebras over $Z$, where $\DR(Z/X) := \Sym_{\mcal O_Z}(\bb L_{Z/X}[1])$.
	\end{prop}
	\begin{proof}
		Consider $C := q_* \mcal O_{\widetilde{\mcal L^{gr}}(Z/X)}$, where $q : \widetilde{\mcal L^{gr}}(Z/X) \to Z$ is the structural map. In fact, we can write this as $\mcal O_{/Z}(\widetilde{\mcal L^{gr}}(Z/X))$ in the sense of \cref{qcoh_dst_adj}, and since the functor $\mcal O_{/Z}$ is left adjoint to $\Spec_Z$, it sends colimits to limits. Therefore, in general, we obtain that
		\[C \simeq \underset{\phi : (\Spec(B) \to \Spec(A)) \to p}{\textup{lim}} v(\phi)_* \DR(B/A)\]
		where $v(\phi) : \Spec(B) \to Z$ denotes the top arrow of each square $\phi \in \Arr(\dAff)/p$.
		
		We will now prove that $C \simeq \DR(Z/X)$, which will be similar to the proof of \cite[Thm.~2.7]{CalaqueSafronov}. Let $u : \Spec(A) \to X$, and consider the following diagram:
		\[\begin{tikzcd}[cramped]
			{Z_u} & Z \\
			{\Spec(A)} & X
			\arrow["{u'}", from=1-1, to=1-2]
			\arrow["{p'}"', from=1-1, to=2-1]
			\arrow["\lrcorner"{anchor=center, pos=0.125}, draw=none, from=1-1, to=2-2]
			\arrow["p", from=1-2, to=2-2]
			\arrow["u"', from=2-1, to=2-2]
		\end{tikzcd}\]
		Then, by \cite[Cor.~B.5.6]{CHS2025}, we obtain that:
		\begin{align*}
			\DR(Z/X) &\simeq \underset{u : \Spec(A) \to X}{\textup{lim}} u'_* u'^* \DR(Z/X) \\
			&\simeq \underset{u : \Spec(A) \to X}{\textup{lim}} u'_* \DR(Z_u/\Spec(A))
		\end{align*}
		where the identification $u'^* \DR(Z/X) \simeq \DR(Z_u/\Spec(A))$ comes from the fact that the square is a pull-back. More precisely, as a pull-back of $p$, $p'$ is locally of finite presentation and admits a relative cotangent which is equivalent to $u'^* \bb L_{Z/X}$. We will now show that:
		\[\DR(Z_u/\Spec(A)) \simeq \underset{v : \Spec(B) \to Z_u}{\textup{lim}} v_* \DR(B/A)\]
		Assuming this is the case, we can end the proof by the following computation:
		\begin{align*}
			\DR(Z/X) &\simeq \underset{u : \Spec(A) \to X}{\textup{lim}} u'_* \DR(Z_u/\Spec(A)) \\
			&\simeq \underset{u : \Spec(A) \to X}{\textup{lim}} u'_* \underset{v : \Spec(B) \to Z_u}{\textup{lim}} v_* \DR(B/A) \\
			&\simeq \underset{u : \Spec(A) \to X}{\textup{lim}} \underset{v : \Spec(B) \to Z_u}{\textup{lim}} u'_* v_* \DR(B/A) \\
			&\simeq \underset{\phi : (\Spec(B) \to \Spec(A)) \to p}{\textup{lim}} v(\phi)_* \DR(B/A) \\
			&\simeq C
		\end{align*}
		Let us denote $Z' := Z_u$ and $X' := \Spec(A)$. In this case, since $X'$ is a derived affine scheme, we can apply \cite[Cor.~1.29]{CalaqueSafronov} to obtain in our notations that $\bb V(\bb L_{Z'/X'}[1]) \simeq \underset{v : \Spec(B) \to Z'}{\textup{colim}} \bb V(\bb L_{B/A}[1])$ as $\bb G_m$-equivariant stacks over $Z$. Therefore, taking graded functions over $Z$, denoted by $\mcal O_{gr}$, and using the fact that it sends colimits to limits, we obtain that $\mcal O_{gr}(\bb V(\bb L_{Z'/X'}[1])) \simeq \underset{v : \Spec(B) \to Z'}{\textup{lim}} \mcal O_{gr}(\bb V(\bb L_{B/A}[1]))$. Using the fact that $p : Z \to X$ is locally of finite presentation, we obtain that $p' : Z' \to X'$ is also locally of finite presentation. As a result, we obtain that $\bb L_{Z'/X'}$ is perfect, and that we can compute the same limit when restricting to affines $v : \Spec(B) \to Z$ such that the composite $\Spec(B) \to Z' \to \Spec(A)$ is of finite presentation, and hence we can assume that all the $\bb L_{B/A}$ are perfect $B$-modules as well. As such, by \cite[Thm.~2.5]{Mou}, we can identify $\mcal O_{gr}(\bb V(\bb L_{Z'/X'}[1])) \simeq \DR(Z'/X')$ as graded $Z'$-algebras, and for any $v : \Spec(B) \to Z'$ satisfying the above condition we have $\mcal O_{gr}(\bb V(\bb L_{B/A}[1])) \simeq v_* \DR(B/A)$ as graded $Z'$-algebras. Combining the above equivalences, we obtain as required that $\DR(Z'/X') \simeq \underset{v' : \Spec(B) \to Z'}{\textup{lim}} v'_* \DR(B/A)$ as graded $Z'$-algebras. This implies by the above computation that $\mcal O_{/Z}(\widetilde{\mcal L^{gr}}(Z/X)) \simeq \DR(Z/X)$ as graded $Z$-algebras.
	\end{proof}
	\begin{notation}
		By abuse of notation, we will also denote by $\DR(Z/X)$ the terminal foliation in the above case.
	\end{notation}
	\begin{rmk}
		Note that $\DR(Z/X)$ usually denotes a sheaf of commutative algebras over $Z$, so it is the data of a commutative $B$-algebra for each map $v : \Spec(B) \to Z$, with transition maps between their usual pull-backs as described in \cite[Not.~2.7.11]{lurieDAGVIII}. On the other hand, the data of a graded mixed algebra on $Z$ relative to $X$ is the data of an $A$-linear graded mixed $B$-algebra for each commutative square $(\Spec(B) \to \Spec(A)) \to p$, where $\Spec(B) \to \Spec(A)$ can be assumed to be of finite presentation as in \cref{def_global_structures} since $p : Z \to X$ is locally of finite presentation, and with transition maps between their pull-backs as graded mixed algebras, which are very different from the usual pull-backs as explained in \cref{pull-back}. Therefore, it is not immediate to see how to associate the second kind of data from the first, and viewing $\DR(Z/X)$ as a graded mixed algebra on $Z$ relative to $X$ does require a clarification, which explains the need to prove the above results.
		
		Moreover, note that even when $p : Z \to X$ does not admit a genuine relative cotangent, for formal reasons there always exists a terminal derived foliation on $Z$ relative to $X$, that is because $\mcal Fol(Z/X)$ is a limit of $\infty$-categories with terminal objects and the transition functors preserve these objects, so the family of these objects is a terminal object in the limit. However, in general outside the above case we expect its associated graded $Z$-algebra to be related to $\Sym_{\mcal O_Z}(\bb L_{Z/X}^{big,qcoh}[1])$, where $\bb L_{Z/X}^{big}$ the big $\mcal O_Z$-module of derivations as defined in \cite[\S A.5]{toen2025derivedfoliations}. We will thus avoid denoting it as $\DR(Z/X)$ outside of the above case, but we suspect the notation might still make sense when for instance the cotangent only exists as a pro-object.
	\end{rmk}
	\subsubsection{Equivalence for Relative Derived Deligne-Mumford Stacks}
	We will now compare $\widetilde{\mcal L^{gr}}$ to $\mcal L^{gr}$, with the final goal of proving that they coincide when $Z \to X$ is a relative derived Deligne-Mumford stack locally of finite presentation in \cref{lgrhat_is_lgr_for_dm}, which combined with \cref{general_foliations} provides a more explicit description of $\mcal Fol(Z/X)$ in this case. To begin with, we will inspect properties of $\mcal L^{gr}$ and $\widetilde{\mcal L^{gr}}$.
	
	More precisely, when $Z \to X$ is a relative derived Deligne-Mumford stack, we will want to show that the canonical map $\widetilde{\mcal L^{gr}}(Z/X) \to \mcal L^{gr}(Z/X)$ is an equivalence by reducing to the case where $X$ is affine and $Z$ is a relative derived Deligne-Mumford stack over $X$. The next two propositions will essentially enable us to do so.
	\begin{prop} \label[prop]{pullback_lgr}
		Let $Z \to X$ be any map of derived stacks, and $u : U \to X$ any map of derived stacks. Then $u^* \mcal L^{gr}(Z/X) \simeq \mcal L^{gr}(u^* Z/U)$ as derived stacks over $U$.
	\end{prop}
	\begin{proof}
		Let $T \to U \in \dSt/U$, we have the following chain of equivalences:
		\begin{align*}
			\Map_U(T,u^* \mcal L^{gr}(Z/X)) &= \Map_U(T,u^* \underline \Map_X(B \bb G_a \times X, Z)) \\
			&\simeq \Map_X(u_! T,\underline \Map_X(B \bb G_a \times X, Z)) \\
			&\simeq \Map_X((u_! T) \times_X (B \bb G_a \times X), Z) \\
			&\simeq \Map_X(u_! (T \times_U (B \bb G_a \times U)), Z) \\
			&\simeq \Map_U(T \times_U (B \bb G_a \times U),u^* Z) \\
			&\simeq \Map_U(T,\underline \Map_U(B \bb G_a \times U,u^* Z)) \\
			&= \Map_U(T,\mcal L^{gr}(u^* Z/U))
		\end{align*}
		where the equivalence $(u_! T) \times_X (B \bb G_a \times X) \simeq u_! (T \times_U (B \bb G_a \times U))$ can be seen in the following diagram:
		\[\begin{tikzcd}[cramped]
			&& Z & \\
			T & U & X & \bullet \\
			{T \times_X (X \times B \bb G_a)} & {U \times_X (X \times B \bb G_a)} & {X \times B \bb G_a} & {B \bb G_a} \\
			{T \times_U (U \times B \bb G_a)} & {U \times B \bb G_a}
			\arrow[from=1-3, to=2-3]
			\arrow[from=2-1, to=2-2]
			\arrow[from=2-1, to=3-1]
			\arrow["u", from=2-2, to=2-3]
			\arrow[from=2-3, to=2-4]
			\arrow["\lrcorner"{anchor=center, pos=0.125, rotate=90}, draw=none, from=3-1, to=2-2]
			\arrow[from=3-1, to=3-2]
			\arrow[from=3-2, to=2-2]
			\arrow["\lrcorner"{anchor=center, pos=0.125, rotate=90}, draw=none, from=3-2, to=2-3]
			\arrow[from=3-2, to=3-3]
			\arrow["\sim", no head, from=3-2, to=4-2]
			\arrow[from=3-3, to=2-3]
			\arrow["\lrcorner"{anchor=center, pos=0.125, rotate=90}, draw=none, from=3-3, to=2-4]
			\arrow[from=3-3, to=3-4]
			\arrow[from=3-4, to=2-4]
			\arrow["\sim", no head, from=4-1, to=3-1]
			\arrow["\lrcorner"{anchor=center, pos=0.125, rotate=90}, draw=none, from=4-1, to=3-2]
			\arrow[from=4-1, to=4-2]
		\end{tikzcd}\]
		Being natural in $T$, this equivalence proves that $u^* \mcal L^{gr}(Z/X) \simeq \mcal L^{gr}(u^* Z/U)$. Note that we can also prove that this equivalence is also $\mcal H$-equivariant, although we will not need this fact. The proof would use a method similar to the proof of \cref{f_*_lgr}.
	\end{proof}
	A similar statement for $\widetilde{\mcal L^{gr}}$ can be obtained, compatibly with the above statement:
	\begin{prop} \label[prop]{pullback_lgrhat}
		Let $p : Z \to X$ be a map of arbitrary derived stacks, and $u : U \to X$ be a map of arbitrary derived stacks. Then, $u^* \widetilde{\mcal L^{gr}}(Z/X) \simeq \widetilde{\mcal L^{gr}}(u^* Z/U)$. Moreover, the following diagram commutes:
		\[\begin{tikzcd}
			{ \widetilde{\mcal L^{gr}}(u^* Z/U)} & {\mcal L^{gr}(u^* Z/U)} \\
			{u^*  \widetilde{\mcal L^{gr}}(Z/X)} & {u^* \mcal L^{gr}(Z/X)}
			\arrow[from=1-1, to=1-2]
			\arrow["\simeq"{description}, no head, from=2-1, to=1-1]
			\arrow[from=2-1, to=2-2]
			\arrow["\simeq"{description}, no head, from=2-2, to=1-2]
		\end{tikzcd}\]
	\end{prop}
	\begin{proof}
		Consider the category of commutative diagrams of the following form:
		\[\begin{tikzcd}
			{U'} & {X'} & {Z'} \\
			U & X & Z
			\arrow[from=1-1, to=1-2]
			\arrow[from=1-1, to=2-1]
			\arrow[from=1-2, to=2-2]
			\arrow[from=1-3, to=1-2]
			\arrow[from=1-3, to=2-3]
			\arrow[from=2-1, to=2-2]
			\arrow[from=2-3, to=2-2]
		\end{tikzcd}\]
		where $U'$, $X'$, and $Z'$ are affines. Let $\dAff/Cospan(U,X,Z)$ denote this category. For simplicity, we will denote its objects as triples $(U',X',Z')$ and omit the maps and their commutativity in the notation.
		
		We know that the colimit of $U' \to X' \leftarrow Z'$ over this category is equivalent to $U \to X \leftarrow Z$. This is essentially because such objects generate the topos $Cospan(\dSt)/(U \to X \leftarrow Z)$, which in turn is implied by \cite[Prop.~6.3.4.9]{HTT}.
		Consider an object of this category, we can form the following commutative diagram:
		\[\begin{tikzcd}
			& {u'^* Z'} && {Z'} \\
			\\
			& {U'} && {X'} \\
			{u^* Z} && Z \\
			\\
			U && X
			\arrow[from=1-2, to=1-4]
			\arrow[dashed, from=1-2, to=3-2]
			\arrow["\lrcorner"{anchor=center, pos=0.125}, draw=none, from=1-2, to=3-4]
			\arrow[from=1-2, to=4-1]
			\arrow[from=1-4, to=3-4]
			\arrow[from=1-4, to=4-3]
			\arrow["{u'}"{description}, dashed, from=3-2, to=3-4]
			\arrow[dashed, from=3-2, to=6-1]
			\arrow[from=3-4, to=6-3]
			\arrow[from=4-1, to=4-3]
			\arrow[from=4-1, to=6-1]
			\arrow["\lrcorner"{anchor=center, pos=0.125}, draw=none, from=4-1, to=6-3]
			\arrow[from=4-3, to=6-3]
			\arrow["u"{description}, from=6-1, to=6-3]
		\end{tikzcd}\]
		Note that $u'^* Z'$ is also affine since $U'$, $X'$ and $Z'$ are affines. We also know that the colimit of the pull-back square formed by $u'^* Z'$, $Z'$, $U'$ and $X'$ is the pull-back square formed by $u^* Z$, $Z$, $U$ and $X$. By \cref{lgrhat_descent}, $\widetilde{\mcal L^{gr}}$ commutes with colimits, which implies that:
		\[\widetilde{\mcal L^{gr}}(u^*Z/U) \simeq \underset{(U',X',Z) \in \dAff/Cospan(U,X,Z)}{\textup{colim}} \widetilde{\mcal L^{gr}}(u'^* Z'/U')\]
		and also:
		\[\widetilde{\mcal L^{gr}}(Z/X) \simeq \underset{(U',X',Z) \in \dAff/Cospan(U,X,Z)}{\textup{colim}} \widetilde{\mcal L^{gr}}(Z'/X')\]
		Moreover, for $(U',X',Z')$ such an object, since $u'^* Z'$ is also affine, we also know that $\widetilde{\mcal L^{gr}}(u'^* Z'/U') \simeq \mcal L^{gr}(u'^* Z'/U')$ by \cref{lgrhat_descent}, and that $\mcal L^{gr}(u'^* Z'/U') \simeq u'^* \mcal L^{gr}(Z'/X')$ by \cref{pullback_lgr}, and finally $u'^* \mcal L^{gr}(Z'/X') \simeq u'^* \widetilde{\mcal L^{gr}}(Z'/X')$ again by \cref{lgrhat_descent}. Therefore, we can compute:
		\begin{align*}
			\widetilde{\mcal L^{gr}}(u^*Z/U) &\simeq \underset{(U',X',Z) \in \dAff/Cospan(U,X,Z)}{\textup{colim}} \widetilde{\mcal L^{gr}}(u'^* Z'/U') \\
			&\simeq \underset{(U',X',Z) \in \dAff/Cospan(U,X,Z)}{\textup{colim}} u'^* \widetilde{\mcal L^{gr}}(Z'/X') \\
			&= \underset{(U',X',Z) \in \dAff/Cospan(U,X,Z)}{\textup{colim}} U' \times_{X'} \widetilde{\mcal L^{gr}}(Z'/X')
		\end{align*}
		Now, since colimits commute with pull-backs, we know that this last colimit is equivalent to the pull-back of the following objects:
		\begin{align*}
			U &\simeq \underset{(U',X',Z) \in \dAff/Cospan(U,X,Z)}{\textup{colim}} U' \\
			X &\simeq \underset{(U',X',Z) \in \dAff/Cospan(U,X,Z)}{\textup{colim}} X' \\
			\widetilde{\mcal L^{gr}}(Z/X) &\simeq \underset{(U',X',Z) \in \dAff/Cospan(U,X,Z)}{\textup{colim}} \widetilde{\mcal L^{gr}}(Z'/X')
		\end{align*}
		where the last equivalence holds because by \cref{lgrhat_descent} since $Z \to X$ is the colimit of the $Z' \to X'$ over this category. Therefore, we obtain that:
		\begin{align*}
			\widetilde{\mcal L^{gr}}(u^*Z/U) &\simeq \underset{(U',X',Z) \in \dAff/Cospan(U,X,Z)}{\textup{colim}} U' \times_{X'} \widetilde{\mcal L^{gr}}(Z'/X') \\
			&\simeq U \times_X \widetilde{\mcal L^{gr}}(Z/X) \\
			&= u^* \widetilde{\mcal L^{gr}}(Z/X)
		\end{align*}
		So we obtain the desired equivalence.
		\separation
		Moreover, by functoriality of $\mcal L^{gr}$, we know that for each object $(U',X',Z') \in \dAff/Cospan(U,X,Z)$, the following square is commutative:
		\[\begin{tikzcd}
			{\mcal L^{gr}(Z'/X')} & {\mcal L^{gr}(Z/X)} \\
			{\mcal L^{gr}(u'^* Z'/U')} & {\mcal L^{gr}(u^*Z/U)} \\
			{u'^* \mcal L^{gr}(Z'/U')} & {u^* \mcal L^{gr}(u^*Z/U)}
			\arrow[from=1-1, to=1-2]
			\arrow[from=2-1, to=1-1]
			\arrow[from=2-1, to=2-2]
			\arrow["\sim"', no head, from=2-1, to=3-1]
			\arrow[from=2-2, to=1-2]
			\arrow["\sim", no head, from=2-2, to=3-2]
			\arrow[from=3-1, to=3-2]
		\end{tikzcd}\]
		Passing to the colimit over $\dAff/Cospan(U,X,Z)$, we obtain that the following square is commutative:
		\[\begin{tikzcd}
			{\widetilde{\mcal L^{gr}}(Z/X)} & {\mcal L^{gr}(Z/X)} \\
			{\widetilde{\mcal L^{gr}}(u^* Z/X)} & {\mcal L^{gr}(u^*Z/U)} \\
			{u^* \widetilde{\mcal L^{gr}}(Z'/U')} & {u^* \mcal L^{gr}(u^*Z/U)}
			\arrow[from=1-1, to=1-2]
			\arrow[from=2-1, to=1-1]
			\arrow[from=2-1, to=2-2]
			\arrow["\sim"', no head, from=2-1, to=3-1]
			\arrow[from=2-2, to=1-2]
			\arrow["\sim", no head, from=2-2, to=3-2]
			\arrow[from=3-1, to=3-2]
		\end{tikzcd}\]
		This proves the second statement of this proposition, and concludes the proof.
	\end{proof}
	As noted above, using these properties we can show the following theorem by restricting to the case where $X$ is affine and $Z$ is a derived Deligne-Mumford stack.
	\begin{thm} \label[thm]{lgrhat_is_lgr_for_dm}
		Assuming $p : Z \to X$ is a relative derived Deligne-Mumford stack, the canonical map:
		\[\widetilde{\mcal L^{gr}}(Z/X) \to \mcal L^{gr}(Z/X)\]
		is an equivalence.
		
		Moreover, since this canonical map is $\mcal H$-equivariant \cref{df-lgrhat}, and is a map of derived stacks over $Z$, these stacks are equivalent as $\mcal H$-equivariant stacks over $X$ and as stacks $\bb G_m$-equivariant derived stacks over $Z$.
	\end{thm}
	\begin{proof}
		To check that $\widetilde{\mcal L^{gr}}(Z/X) \to \mcal L^{gr}(Z/X)$ is an equivalence, it suffices to check that the pull-back of this map along any map $u : U \to X$ is an equivalence, where $U$ is affine.
		Let $u : U \to X$ be such a map, by \cref{pullback_lgrhat}, we know that the pull-back of the above canonical map is exactly the canonical map $\widetilde{\mcal L^{gr}}(u^*Z/U) \to \mcal L^{gr}(u^* Z/U)$, and therefore in particular it is $\mcal H$-equivariant. Notice that by assumption, since $Z \to X$ is a relative derived Deligne-Mumford stack, $u^* Z$ is now a derived Deligne-Mumford stack, and $U$ is affine. Therefore, we are reduced to proving the statement in the case where $X$ is affine and $Z$ is a derived Deligne-Mumford stack. We will now assume that this is the case.
		
		By \cref{def_Lgr}, we have that:
		\[\mcal L^{gr}(Z/X) := \underline{\Map}_X(B \bb G_a \times X,Z)\]
		Note that, as derived stacks over $X$, there is a canonical equivalence:
		\[\underline{\Map}_X(B \bb G_a \times X,Z) \simeq \underline{\Map}(B \bb G_a,Z) \times_{\underline{\Map}(B \bb G_a,X)} X\]
		where the map $u:\Map(B\mathbb G_a,Z)\to \Map(B\mathbb G_a,X)$ is induced by $p : Z \to X$ and the map $s:X\to \Map(B\mathbb G_a,X)$ is the transpose of the projection $B \bb G_a \times X \to X$.
		To see this, it suffices to show that for every \(y:Y\to X\) in \(\dSt/X\), there is a natural equivalence
		\[
		\Map_{\dSt/X}\!\bigl(Y,\underline{\Map}_X(B\mathbb G_a\times X,Z)\bigr)
		\simeq
		\Map_{\dSt/X}\!\bigl(Y,\underline{\Map}(B\mathbb G_a,Z)\times_{\underline{\Map}(B\mathbb G_a,X)}X\bigr).
		\]
		Let $y : Y \to X$, we have:
		\begin{align*}
			\Map_{\dSt/X}\!\bigl(Y,\underline{\Map}_X(B\mathbb G_a\times X,Z)\bigr) &\simeq \Map_{\dSt/X}\!\bigl(Y\times_X(B\mathbb G_a\times X),Z\bigr) \\
			&\simeq \Map_{\dSt/X}(Y\times B\mathbb G_a,Z)
		\end{align*}
		since $Y \times_X (B \bb G_a \times X) \simeq Y \times B \bb G_a$ as derived stacks over $X$, where the structural map $Y \times B \bb G_a \to X$ is the composite of the projection $pr_1 : Y \times B \bb G_a \to Y$ with $y : Y \to X$.
		Now, a morphism \(f:Y\times B\mathbb G_a\to Z\) over \(X\) is precisely the data of a morphism $f : Y \times B \bb G_a \to Z$ in \(\dSt\) together with a homotopy $p\circ f \simeq y\circ pr_1$.
		By adjunction, the data of \(f\) in $\dSt$ corresponds to a morphism $\widetilde f:Y\to \underline{\Map}(B\mathbb G_a,Z)$, and under this correspondence the data of an equivalence \(p\circ f \simeq y\circ pr_1\) becomes the data of an equivalence $u \circ \widetilde f \simeq s \circ y :Y\to \underline{\Map}(B\mathbb G_a,X)$, where as above $s:X\to \underline{\Map}(B\mathbb G_a,X)$ is the adjoint of the projection $B \bb G_a \times X \to X$, and $u:\underline{\Map}(B \bb G_a,Z) \to \underline{\Map}(B \bb G_a,X)$ is induced by $p$. Since $X$ is the terminal object of $\dSt/X$, the only map $Y \to X$ in $\dSt/X$ is $y$, and thus we obtain that:
		\[
		\Map_{\dSt/X}(Y\times B\mathbb G_a,Z)
		\simeq
		\Map_{\dSt/X}\!\bigl(Y,\underline{\Map}(B\mathbb G_a,Z)\times_{\underline{\Map}(B\mathbb G_a,X)}X\bigr).
		\]
		
		Combining the previous equivalences, in total we obtain:
		\[
		\Map_{\dSt/X}\!\bigl(Y,\underline{\Map}_X(B\mathbb G_a\times X,Z)\bigr)
		\simeq
		\Map_{\dSt/X}\!\bigl(Y,\underline{\Map}(B\mathbb G_a,Z)\times_{\underline{\Map}(B\mathbb G_a,X)}X\bigr).
		\]
		and this equivalence is natural in $y : Y \to X$. As a result, we obtain the desired equivalence of derived stacks over $X$:
		\[
		\mcal L^{gr}(Z/X) := \underline{\Map}_X(B\mathbb G_a\times X,Z) \simeq \underline{\Map}(B\mathbb G_a,Z)\times_{\underline{\Map}(B\mathbb G_a,X)}X
		\]
		
		Since $Z$ is now a derived Deligne-Mumford stack, there exists an étale hypercovering $v_\bullet : V_\bullet \to Z$ where the $V_\bullet$ are disjoint union of affines. In particular, since the colimit of the maps $(V_\bullet \to X)$ in $\Arr(\dSt)/p$ is the map $(Z \to X)$, by \cref{lgrhat_descent} we know that:
		\[\widetilde{\mcal L^{gr}}(Z/X) \simeq \textup{colim} \; \widetilde{\mcal L^{gr}}(V_\bullet/X) \simeq \textup{colim} \; \mcal L^{gr}(V_\bullet/X)\]
		where the last equivalence also holds by \cref{lgrhat_descent} since the $V_\bullet$ are disjoint union of affines and $X$ is affine. Note that by functoriality of $\mcal L^{gr}$, there are canonical maps $\mcal L^{gr}(V_\bullet/X) \to \mcal L^{gr}(Z/X)$, which induces a canonical map $\textup{colim} \; \mcal L^{gr}(V_\bullet/X) \to \mcal L^{gr}(Z/X)$. Moreover, the following triangle commutes by definition:
		\[\begin{tikzcd}[cramped]
			{\widetilde{\mcal L^{gr}}(Z/X)} && {\mcal L^{gr}(Z/X)} \\
			\\
			& {\textup{colim} \; \mcal L^{gr}(V_\bullet/X)}
			\arrow[from=1-1, to=1-3]
			\arrow["\sim", from=3-2, to=1-1]
			\arrow[from=3-2, to=1-3]
		\end{tikzcd}\]
		Therefore, it suffices to show that $\textup{colim} \; \mcal L^{gr}(V_\bullet/X) \to \mcal L^{gr}(Z/X)$ is an equivalence to finish the proof. By \cite[Thm.~2.11]{fu2025hochschild}, we know that $\textup{colim} \; \underline{\Map}(B \bb G_a,V_\bullet) \xrightarrow{\sim} \underline{\Map}(B \bb G_a,Z)$ is an equivalence. Therefore, by pull-back, the map:
		\[(\textup{colim} \; \underline{\Map}(B\mathbb G_a,V_\bullet))\times_{\underline{\Map}(B\mathbb G_a,X)}X \xrightarrow{\sim} \underline{\Map}(B\mathbb G_a,Z)\times_{\underline{\Map}(B\mathbb G_a,X)}X\]
		is an equivalence. Moreover, since colimits commute with pull-backs, we obtain that:
		\[\textup{colim} \; \bigl( \underline{\Map}(B\mathbb G_a,V_\bullet)\times_{\underline{\Map}(B\mathbb G_a,X)}X \bigr) \xrightarrow{\sim} \underline{\Map}(B\mathbb G_a,Z)\times_{\underline{\Map}(B\mathbb G_a,X)}X\]
		is also an equivalence. Finally, as identified above, this is exactly the fact that the map:
		\[\textup{colim} \; \mcal L^{gr}(V_\bullet/X) \xrightarrow{\sim} \mcal L^{gr}(Z/X)\]
		is an equivalence. This concludes the proof.
	\end{proof}
	To finish this section, note that combining this result with \cref{general_foliations}, we obtain the following result.
	\begin{thm} \label[thm]{equiv_fol_dst_rel_dm}
		If $Z \to X$ is a relative derived Deligne-Mumford stack locally of finite presentation, then there is a canonical equivalence of $\infty$-categories:
		\[\mcal Fol(Z/X) \simeq \PerfLinSt_Z/\mcal L^{gr}(Z/X) - \mcal H\]
	\end{thm}
	\begin{proof}
		By \cref{general_foliations}, we know that:
		\[\mcal Fol(Z/X) \simeq \PerfLieSt/\widetilde{\mcal L^{gr}}(Z/X) - \mcal H\]
		By the above result, we obtain that:
		\[\widetilde{\mcal L^{gr}}(Z/X) \simeq \mcal L^{gr}(Z/X)\]
		as $\mcal H$-equivariant derived stacks over $X$ and as $\bb G_m$-equivariant derived stacks over $Z$. Moreover, in this case we will now show that $\widetilde{\mcal L^{gr}}(Z/X)$ is relatively affine over $Z$. Notice that by \cref{Ogr_L_gr}, this implies that $\widetilde{\mcal L^{gr}}(Z/X) \simeq \Spec_Z(\DR(Z/X)) = \bb V(\bb L_{Z/X}[1])$ as $\bb G_m$-equivariant derived stacks over $Z$. Using this fact, we will later on prove that perfect Lie stacks coincide with perfect linear stacks.
		
		To prove that $\widetilde{\mcal L^{gr}}(Z/X)$ is relatively affine over $Z$, notice that it suffices to prove that for any map $u : \Spec(A) \to X$, the pull-back $u^* \widetilde{\mcal L^{gr}}(Z/X)$ is relatively affine over $Z_u := Z \times_X \Spec(A)$. Indeed, assuming this is the case, for any $v : \Spec(B) \to Z$ we can consider the diagram:
		\[\begin{tikzcd}[ampersand replacement=\&,cramped]
			{v^* \widetilde{\mcal L^{gr}}(Z/X)} \& {u^* \widetilde{\mcal L^{gr}}(Z/X)} \&\& {\widetilde{\mcal L^{gr}}(Z/X)} \\
			{\Spec(B)} \& {Z_u} \&\& Z \\
			\& {\Spec(B)} \& Z \& X
			\arrow[from=1-1, to=1-2]
			\arrow[from=1-1, to=2-1]
			\arrow["\lrcorner"{anchor=center, pos=0.125}, draw=none, from=1-1, to=2-2]
			\arrow[from=1-2, to=1-4]
			\arrow[from=1-2, to=2-2]
			\arrow["\lrcorner"{anchor=center, pos=0.125}, draw=none, from=1-2, to=2-4]
			\arrow[from=1-4, to=2-4]
			\arrow[from=2-1, to=2-2]
			\arrow[from=2-2, to=2-4]
			\arrow[from=2-2, to=3-2]
			\arrow["\lrcorner"{anchor=center, pos=0.125}, draw=none, from=2-2, to=3-4]
			\arrow[from=2-4, to=3-4]
			\arrow[from=3-2, to=3-3]
			\arrow["u"', shift right=3, from=3-2, to=3-4]
			\arrow[from=3-3, to=3-4]
		\end{tikzcd}\]
		which presents $v^* \widetilde{\mcal L^{gr}}(Z/X) \to \Spec(B)$ as the pull-back of $u^* \widetilde{\mcal L^{gr}}(Z/X) \to Z_u$ and is thus relatively affine by assumption, so $\widetilde{\mcal L^{gr}}(Z/X)$ is relatively affine. Moreover, since $u^* \widetilde{\mcal L^{gr}}(Z/X) \simeq \widetilde{\mcal L^{gr}}(Z_u/\Spec(A))$ by \cref{pullback_lgrhat}, we are thus reduced to proving the statement when $X = \Spec(A)$ is a derived affine.
		
		In this case, $Z$ is a derived Deligne-Mumford stack, so it admits an étale affine hypercover $V_\bullet \to Z$. For each $n \in \bb N$, we consider $\mcal L_n = \mcal L^{gr}(V_n/\Spec(A)) \simeq \Spec_{V_n}(\DR(V_n/\Spec(A)))$. Since the transition maps between the $V_\bullet$ are themselves étale, we get in particular that $\mcal L_n \times_{V_n} V_m \simeq \mcal L_m$. As such, this family of objects defines an object in $\textup{lim} \dSt/V_\bullet$. Moreover, since the colimit of the $V_\bullet$ is $Z$, and the functor $T \mapsto \dSt/T$ sends colimits to limits, we obtain that $\textup{lim} \dSt/V_\bullet \simeq \dSt/Z$, where the functor from left to right is given by taking the colimit, and the converse functor is the pull-back over $V_\bullet \to Z$. As a consequence, the above family $\mcal L_\bullet$ defines a derived stack $\mcal L \to Z$ such that for any $n$, $\mcal L \times_Z V_n \simeq \mcal L_n$, and in particular this pull-back is relatively affine over $V_n$.
		
		Moreover, note that the family of arrows $V_n \to \Spec(A)$ has colimit $Z \to \Spec(A)$ in $\Arr(\dSt)$, so by \cref{lgrhat_descent} we know that $\widetilde{\mcal L^{gr}}(Z/\Spec(A)) \simeq \textup{colim} \mcal L_\bullet = \mcal L$. Finally, also consider $\Spec_Z(\DR(Z/X)) \to Z$, and note that by the transitivity of the cotangent over $v_n : V_n \to Z$ we get canonical maps $v_n^* \bb L_Z \to \bb L_{V_n}$ which promote to maps $\mcal L_n \to \Spec_Z(\DR(Z/X))$ of derived stacks over $Z$ and finally produce a map $\mcal L \to \Spec_Z(\DR(Z/X))$. To show that this map is an equivalence, it suffices to show that its pull-back to each $V_n$ is an equivalence since they cover $Z$, but this pull-back is exactly the map induced by $v_n^* \bb L_Z \to \bb L_{V_n}$ by construction, which is indeed an equivalence because $v_n : V_n \to Z$ is étale by assumption.
		
		As a result, we obtain even when $X$ is not affine that $\widetilde{\mcal L^{gr}}(Z/X) \to Z$ is relatively affine. Moreover, by \cref{Ogr_L_gr} and by the previous theorem, we deduce in fact that in the general relative derived Deligne-Mumford case we have $\mcal L^{gr}(Z/X) \simeq \widetilde{\mcal L^{gr}}(Z/X) \simeq \Spec_Z(\DR(Z/X))$. Using this fact, we can now prove that:
		\[\PerfLieSt/\Spec_Z(\DR(Z/X)) - \mcal H \simeq \PerfLinSt_Z/\Spec_Z(\DR(Z/X)) - \mcal H\]
		
		Let $L \in \PerfLieSt/\widetilde{\mcal L^{gr}}(Z/X) - \mcal H$, and let us show that $L \in \PerfLinSt_Z$. Since $L$ is relatively affine over $\Spec_Z(\DR(Z/X))$, let us write $L \simeq \Spec_{\Spec_Z(\DR(Z/X))}(C)$, where $C \in \CAlg(\QCoh(\Spec_Z(\DR(Z/X)))) \simeq \DR(Z/X)-\CAlg(\QCoh(Z))$. Note that since the map $L \to \Spec_Z(\DR(Z/X))$ is $\bb G_m$-equivariant, $C$ comes equipped with a canonical grading, and the structural map $\DR(Z/X) \to C$ is a graded map. Consider the canonical graded map $\Sym_{\mcal O_Z}(C(1)) \to C$, and note that the structural map $\DR(Z/X) \to C$ factors through $\Sym_{\mcal O_Z}(C(1)) \to C$ so this latter map is a map of graded $\DR(Z/X)$-algebras. Our goal is to show that it is an equivalence, and that $C(1)$ is perfect.
		
		Consider now $L' = \Spec_{\Spec_Z(\DR(Z/X))}(\Sym_{\mcal O_Z}(C(1)))$, and notice that there is a canonical $\bb G_m$-action on $L'$ corresponding to the canonical grading of its algebra, and a canonical $\bb G_m$-equivariant map $L \to L'$ over $\Spec_Z(\DR(Z/X))$ corresponding to the above graded map of $\DR(Z/X)$-algebras, and moreover this map is an equivalence if and only if $L \to L'$ is an equivalence. Since $\Spec_Z(\DR(Z/X)) \simeq \widetilde{\mcal L^{gr}}(Z/X)$, to show that $L \to L'$ is an equivalence it suffices to show that for any commutative square $(V \to U) \to (Z \to X) \in \Arr(\dAff)^{fp}/p$, the pull-back map $L \times_{\Spec_Z(\DR(Z/X))} \Spec_V(\DR(V/U)) \to L' \times_{\Spec_Z(\DR(Z/X))} \Spec_V(\DR(V/U))$ is an equivalence.
		
		Let $(V \to U) \to (Z \to X)$ be such a square, and consider $v : V \to Z$ the induced map. We compute directly that $L \times_{\Spec_Z(\DR(Z/X))} \Spec_V(\DR(V/U)) \simeq \Spec_V(v^* C \otimes_{v^* \DR(Z/X)} \DR(V/U))$ and $L' \times_{\Spec_Z(\DR(Z/X))} \Spec_V(\DR(V/U)) \simeq \Spec_V(\Sym_V(v^* C(1)) \otimes_{v^* \DR(Z/X)} \DR(V/U))$, and the corresponding map between the two is the induced canonical graded map $\Sym_V(v^* C(1)) \otimes_{v^* \DR(Z/X)} \DR(V/U) \to v^* C \otimes_{v^* \DR(Z/X)} \DR(V/U)$. Note that the first graded algebra is equivalent to $\Sym_V(v^* C(1) \oplus_{v^* \bb L_{Z/X}} \bb L_{V/U})$. Since $L$ is a perfect Lie stack, by assumption $v^* C \otimes_{v^* \DR(Z/X)} \DR(V/U)$ is a quasi-free graded algebra. Moreover, we compute that its weight-1 part is exactly $v^* C(1) \oplus_{v^* \bb L_{Z/X}} \bb L_{V/U}$ which coincides with the previous computation, so the above map is indeed an equivalence, and we immediately obtain that the pull-back of $L \to L'$ is an equivalence and thus that $L \to L'$ is an equivalence and finally that $\Sym_{\mcal O_Z}(C(1)) \to C$ is an equivalence of graded algebras in $\QCoh(Z)$.
		
		It remains to check that $C(1)$ is perfect. Consider now $V \to Z$ any derived affine over $Z$, and the induced commutative square $(V = V) \to (Z \to X)$ and map $\mcal L^{gr}(V/V) \simeq V \to \Spec_Z(\DR(Z/X))$ as in \cref{perf_lie_st_rmk}. In this case, the above computation over the identity map $V = V$ wheere $\DR(V/V) \simeq V$ concentrated in weight 0 yields that $v^* C \otimes_{v^* \DR(Z/X)} V \simeq \Sym_V(v^* C(1) \oplus_{v^* \bb L_{Z/X}} 0)$ and we set $E_V := v^* C(1) \oplus_{v^* \bb L_{Z/X}} 0 = \cofib(v^* \bb L_{Z/X} \to v^* C(1))$, which is perfect since $L$ is a perfect Lie stack. As a consequence, since $v^* C(1)$ is the fiber of $v^* \bb L_{Z/X} \to E_V$ and both are perfect, $v^* C(1)$ is perfect. Since this is true for all $v : V \to Z$, we obtain that $C(1)$ is perfect. Finally, as a result, we obtain that $L$ is a perfect linear stack over $Z$.
		
		Conversely, let $L = \Spec_Z(\Sym_{\mcal O_Z}(E))$ be a perfect linear stack over $Z$ with a structural $\mcal H$-equivariant map to $\widetilde{\mcal L^{gr}}(Z/X) \simeq \Spec_Z(\DR(Z/X))$. In this case, for any commutative square $(V \to U) \to (Z \to X)$ where $U$ and $V$ are derived affines with $V \to U$ of finite presentation and induced map $v : V \to Z$, we can compute directly that $L \times_{\Spec_Z(\DR(Z/X))} \Spec_V(\DR(V/U)) \simeq \Spec_V(\Sym_V(v^* E \oplus_{v^* \bb L_{Z/X}} \bb L_{V/U}))$, which is a perfect linear stack over $Z$ since $v^* E \oplus_{v^* \bb L_{Z/X}} \bb L_{V/U}$ is perfect because $v^*$ preserves perfects, $E$ is perfect by assumption, and both $\bb L_{Z/X}$ and $\bb L_{V/U}$ are perfects because $Z \to X$ and $V \to U$ are locally of finite presentation. As such, $L$ is a perfect Lie stack. This concludes the proof.
	\end{proof}
	\begin{rmk}
		We expect the equivalence between perfect Lie stacks and perfect linear stacks to hold precisely when $\widetilde{\mcal L^{gr}}(Z/X) \to Z$ is relatively affine, and we only expect this to hold in the relative derived Deligne-Mumford case. In general, we expect these notions to be genuinely distinct.
	\end{rmk}
	\section{Pushforward of Derived Foliations} \label[section]{Push-forward}
	Let $Z \to X$ and $f : X \to Y$ be two maps of derived stacks. Our goal in this section is to construct a functor $f_{*,\mcal Fol} : \mcal Fol(Z/X) \to \mcal Fol(f_* Z/Y)$ under some assumptions on $Z \to X$ and $f$, and use it to obtain in particular a derived foliation on some mapping stacks when the target comes equipped with one. Our plan to construct $f_{*,\mcal Fol} : \mcal Fol(Z/X) \to \mcal Fol(f_* Z/Y)$ is to use \cref{equiv_fol_dst_rel_dm}, so we will later restrict to the case where $Z \to X$ is a relative derived Deligne-Mumford stack locally of finite presentation. We will consider that an object of $\mcal Fol(Z/X)$ is an $\mcal H$-equivariant stack $\bb V(E) \to \mcal L^{gr}(Z/X)$ where $E \in \Perf(Z)$, construct its Weil restriction $f_* \bb V(E)$ as in \cref{base_chg_adj} and as recalled below, and identify it with a perfect linear stack on $f_* Z$, equipped with an $\mcal H$-action and a canonical $\mcal H$-equivariant map to $\mcal L^{gr}(f_* Z/Y)$, and finally identify this data as a derived foliation on $f_* Z$ relative to $Y$.
	
	More concretely, recall from \cref{base_chg_adj} that there is a functor
	\[f_* : \dSt/X \to \dSt/Y\]
	called the push-forward, direct image, or Weil restriction along $f$. By definition, it is the right adjoint to the usual pull-back $f^* : \dSt/Y \to \dSt/X$, and therefore given $U \to X$, the object $f_*(U)$ in $\dSt/Y$ satisfies the following universal property for all $V \to Y$:
	\[\Map_Y(V,f_*(U)) := \Map_X(V \times_Y X,U)\]
	Now, notice that, by functoriality of $f_*$, the data $\bb V(E) \to \mcal L^{gr}(Z/X) \to Z$ is mapped to $f_* \bb V(E) \to f_* \mcal L^{gr}(Z/X) \to f_* Z$.
	
	In the first subsection of this section, with no assumption on $Z \to X$, we will prove that $f_* \mcal L^{gr}(Z/X) \simeq \mcal L^{gr}(f_* Z/Y)$ as $\mcal H$-equivariant derived stacks over $f_* Z$, and that $f_* \bb V(E)$ is perfect linear over $f_* Z$, when assuming that $f$ is proper-schematic and local complete intersection. This consitutes a pushforward of perfect linear stacks, as established in \cref{pushforward_preserves_fols}.
	
	In the second subsection of this section, assuming that $f_* Z \to Y$ is a relative derived Artin stack, we will note that $f_* Z \to Y$ will moreover be a derived Deligne-Mumford stack locally of finite presentation when adding that $f$ is moreover flat, thus by \cref{equiv_fol_dst_rel_dm} we have that $f_* \bb V(E)$ corresponds to a derived foliation over $f_* Z$ relative to $Y$, which yields the functor $f_{*,\mcal Fol}$ in \cref{thm_pushforward_fol-restate} (\ref{thm_pushforward_fol}). Finally, applying this theorem to \cref{ex_main_thm} will yield the main \cref{main_theorem-restate} (\ref{main_theorem}) on mapping stacks.
	
	\subsection{Pushforward of Perfect Linear Stacks}
	Let us now follow the steps described above to construct the push-forward of perfect linear stacks, starting with two remarks on the push-forward, i.e.\ the Weil restriction along $f$ of stacks over $Z$, as introduced in \cref{base_chg_adj}. In all of this subsection, $Z \to X$ will be an arbitrary map of derived stacks.
	\begin{df} \label[df]{def_f_*_Z}
		Let $Z \to X$ be an arbitrary map of derived stacks, and $f : X \to Y$ an arbitrary map of derived stacks. We will denote by:
		\[f^Z_* : \dSt/Z \to \dSt/f_* Z\]
		the functor mapping a derived stack $T \to Z$, viewed as a map in $\dSt/X$, to the map $f_* T \to f_* Z$ in $\dSt/Y$, providing an object in $\dSt/f_* Z$.
	\end{df}
	\begin{rmk}
		The preceding functor is right adjoint to the functor mapping a derived stack $t' : T' \to f_* Z$ to the composite $f^* T' \xrightarrow{f^* t'} f^* f_* Z \xrightarrow{ev} Z$, where the map we will call $ev : f^* f_* Z \to Z$ is the co-unit of the adjunction $f^* \dashv f_* : \dSt/Y \rightleftarrows \dSt/X$.
		This left adjoint will be denoted by:
		\[f_Z^* : \dSt/f_* Z \to \dSt/Z\]
	\end{rmk}
	\begin{cor} \label[cor]{pushforward_monoidal}
		As a right adjoint, the functor $f^Z_* : \dSt/Z \to \dSt/f_* Z$ commutes with limits and thus has a canonical monoidal structure with respect to the cartesian product.
	\end{cor}
	After establishing these results on $f^Z_*$, the next step in our construction consists in providing a canonical equivalence between $f^Z_* \mcal L^{gr}(Z/X)$ and $\mcal L^{gr}(f_* Z/Y)$. Note that graded loop stacks come equipped with canonical $\mcal H$-actions, and since $f^Z_*$ is monoidal by the above corollary, $f^Z_* \mcal L^{gr}(Z/X)$ also inherits an $\mcal H$-action. Therefore, it makes sense to further prove that this canonical map is an equivalence of $\mcal H$-equivariant stacks over $f_* Z$, which we will need later on for our construction.
	\begin{lem} \label[lem]{f_*_lgr}
		There is a canonical equivalence $f^Z_* \mcal L^{gr}(Z/X) \simeq \mcal L^{gr}(f_* Z/Y)$ as $\mcal H$-equivariant derived stacks over $Y$ and as $\bb G_m$-equivariant derived stacks over $f_* Z$.
	\end{lem}
	\begin{proof}
		Note that, by \cite[Prop.~3.2.77]{SatiUrs26} we obtain that:
		\begin{align*}
			[\mcal L^{gr}(Z/X)/\mcal H] &= [\underline{\Map}_X(X \times B \bb G_a,Z)/\mcal H] \\
			&\simeq \underline{\Map}_{[X / \mcal H]}([X \times B \bb G_a / \mcal H],[Z / \mcal H]) \\
			&\simeq \underline{\Map}_{X \times B \mcal H}(X \times [B \bb G_a / \mcal H],Z \times B \mcal H)
		\end{align*}
		Moreover, in \cref{props_p_pi}, we proved that we have the following pull-back square:
		\[\begin{tikzcd}[ampersand replacement=\&,cramped]
			{B \bb G_a} \& {B \bb G_m} \\
			\bullet \& {B \mcal H}
			\arrow[from=1-1, to=1-2]
			\arrow[from=1-1, to=2-1]
			\arrow["\lrcorner"{anchor=center, pos=0.125}, draw=none, from=1-1, to=2-2]
			\arrow[from=1-2, to=2-2]
			\arrow[from=2-1, to=2-2]
		\end{tikzcd}\]
		This implies that $[B \bb G_a/\mcal H] \simeq B \bb G_m$. Therefore, in conclusion we have:
		\[[\mcal L^{gr}(Z/X)/\mcal H] \simeq \underline{\Map}_{X \times B \mcal H}(X \times B \bb G_m,Z \times B \mcal H)\]
		Moreover, let us denote by $\overline{f} : X \times B \mcal H \to Y \times B \mcal H$ the quotient of the map $f : X \to Y$ by the trivial $\mcal H$-action on both sides. Then, by \cite[Lem.~2.1]{FreyRasekh23}, we have:
		\[\overline{f}_* [\mcal L^{gr}(Z/X)/\mcal H] \simeq [f_* \mcal L^{gr}(Z/X)/\mcal H]\]
		We can thus compare $f^Z_* \mcal L^{gr}(Z/X)$ and $\mcal L^{gr}(f_* Z/Y)$ as $\mcal H$-equivariant derived stacks over $Y$ directly. Let $T \in \dSt/Y - \mcal H$, we have:
		\begin{align*}
			&\Map_{\dSt/Y - \mcal H}(T,f_* \mcal L^{gr}(Z/X)) \\
			\simeq\ &\Map_{\dSt/(Y \times B \mcal H)}([T/\mcal H], \overline{f}_* [\mcal L^{gr}(Z/X)/\mcal H]) \\
			\simeq\ &\Map_{\dSt/(X \times B \mcal H)}(\overline{f}^* [T/\mcal H],\underline{\Map}_{X \times B \mcal H}(X \times B \bb G_m,Z \times B \mcal H)) \\
			\simeq\ &\Map_{\dSt/(X \times B \mcal H)}(\overline{f}^* [T/\mcal H] \times_{X \times B \mcal H} (X \times B \bb G_m),Z \times B \mcal H) \\
			\simeq\ &\Map_{\dSt/(X \times B \mcal H)}(\overline{f}^* [T/\mcal H] \times_{X \times B \mcal H} \overline{f}^*(Y \times B \bb G_m),Z \times B \mcal H) \\
			\simeq\ &\Map_{\dSt/(X \times B \mcal H)}(\overline{f}^* ([T/\mcal H] \times_{Y \times B \mcal H} (Y \times B \bb G_m)),Z \times B \mcal H) \\
			\simeq\ &\Map_{\dSt/(Y \times B \mcal H)}([T/\mcal H] \times_{Y \times B \mcal H} (Y \times B \bb G_m),\overline{f}_*(Z \times B \mcal H)) \\
			\simeq\ &\Map_{\dSt/(Y \times B \mcal H)}([T/\mcal H],\underline{\Map}_{Y \times B \mcal H}(Y \times B \bb G_m,\overline{f}_*(Z \times B \mcal H))) \\
			\simeq\ &\Map_{\dSt/(Y \times B \mcal H)}([T/\mcal H],\underline{\Map}_{Y \times B \mcal H}(Y \times B \bb G_m,(f_* Z) \times B \mcal H)) \\
			\simeq\ &\Map_{\dSt/(Y \times B \mcal H)}([T/\mcal H],[\mcal L^{gr}(f_* Z/Y)/\mcal H]) \\
			\simeq\ &\Map_{\dSt/Y - \mcal H}(T,\mcal L^{gr}(f_* Z/Y))
		\end{align*}
		Moreover, this chain of equivalences is canonically functorial in $T$. This provides an equivalence $f^Z_* \mcal L^{gr}(Z/X) \simeq \mcal L^{gr}(f_* Z/Y)$ as $\mcal H$-equivariant derived stacks over $Y$. It remains to check that this equivalence is also an equivalence of $\bb G_m$-equivariant stacks over $f_* Z$. Recall that the map $ev_X : \mcal L^{gr}(Z/X) \to Z$ comes from the evaluation at the base point $e_X : X \times \bullet \to X \times B \bb G_a$. Equivalently, thanks to the identification $\underline{\Map}_X(X,Z) \simeq Z$ as derived stacks over $X$, this is obtained by contravariant functoriality of $\underline{\Map}_X(-,Z)$ applied to $e_X$ in the first argument, so we will denote $ev_X$ by $e_X^*$, and similarly for $ev_Y : \mcal L^{gr}(f_* Z/Y) \to f_* Z$. Therefore, we have to prove that the following square commutes:
		\[\begin{tikzcd}[ampersand replacement=\&,cramped]
			{f_* \underline{\Map}_X(X \times B \bb G_a,Z)} \& {\underline{\Map}_Y(Y \times B \bb G_a,Z)} \\
			{f_* \underline{\Map}_X(X,Z)} \& {\underline{\Map}_Y(Y,f_*Z)}
			\arrow["{f_*(e_X^*)}"', from=1-1, to=2-1]
			\arrow["\sim"', from=1-2, to=1-1]
			\arrow["{e_Y^*}", from=1-2, to=2-2]
			\arrow["{Id_{f_* Z}}"', equals, from=2-1, to=2-2]
		\end{tikzcd}\]
		By \cite[Prop.~3.2.63]{SatiUrs26}, under the adjunction $f^* \dashv f_*$, this is equivalent to the commutativity of the following square:
		\[\begin{tikzcd}[ampersand replacement=\&,cramped,column sep=scriptsize]
			{\underline{\Map}_X(X \times B \bb G_a,Z)} \& {\underline{\Map}_X(X \times B \bb G_a,f^* f_* Z)} \& {f^* \underline{\Map}_Y(Y \times B \bb G_a,f_* Z)} \\
			{\underline{\Map}_X(X,Z)} \& {f^* f_* Z} \& {f^* \underline{\Map}_Y(Y,f_*Z)}
			\arrow["{e_X^*}"', from=1-1, to=2-1]
			\arrow["{\epsilon_Z \circ -}"', from=1-2, to=1-1]
			\arrow[equals, from=1-2, to=1-3]
			\arrow["{e^*_X}", from=1-2, to=2-2]
			\arrow["{f^*(e_Y^*)}", from=1-3, to=2-3]
			\arrow["{\epsilon_Z}", from=2-2, to=2-1]
			\arrow[equals, from=2-2, to=2-3]
		\end{tikzcd}\]
		where $\epsilon_Z : f^* f_* Z \to Z$ is the co-unit of the adjunction. By the naturality of the same proposition \cite[Prop.~3.2.63]{SatiUrs26} in both arguments, this diagram commutes, and thus $f^Z_* \mcal L^{gr}(Z/X) \simeq \mcal L^{gr}(f_* Z/Y)$ as derived stacks over $f_* Z$. To see that the equivalence is $\bb G_m$-equivariant over $f_* Z$, it only remains to notice that $f_*(e_X^*)$ is $\bb G_m$-equivariant because $e_X^*$ is $\bb G_m$-equivariant and $f_*$ is monoidal. Since $e_Y^*$ is also $\bb G_m$-equivariant, the above commutative square promotes to a commutative square of $\bb G_m$-equivariant stacks, which concludes the proof.
	\end{proof}
	\begin{rmk}
		In the above context, note that the fact that $e_X^*$ is $\bb G_m$-equivariant can also be seen by the fact that the map $\bullet \to B \bb G_a$ is $\bb G_m$-equivariant, or equivalently that the action of $\bb G_m$ on $B \bb G_a$ preserves the canonical point, which is the case because this action is the delooping of the action of $\bb G_m$ on $\bb G_a$.
	\end{rmk}
	The next step in our construction now consists in proving that, assuming at some point that $f$ is proper-schematic and local complete intersection, and given $E \in \Perf(Z)$, the push-forward of $\bb V(E) \to Z \to X$ along $f$ is still perfect linear over $f_* Z$.
	The main difficulty in this situation is that the universal property of $f_*$ applies after forgetting the structural map $f_* \bb V(E) \to f_* Z$ along $f_* Z \to Y$, while the universal property of linear stacks depends on viewing them over $f_* Z$. Therefore, in order to prove our statement, we will first prove that $f^Z_*$ actually coincides with another functor having a universal property over $f_* Z$. Later on, we will apply this description to perfect linear stacks while assuming that $f$ is proper-schematic and local complete intersection to compute the functor of points of $f_* \bb V(E)$ over $f_* Z$ and show that it also corresponds to a perfect linear stack.
	More precisely, consider the following commutative diagram:
	\[\begin{tikzcd}
		Z & {f^* f_* Z} & {f_*Z} \\
		& X & Y
		\arrow["f"', from=2-2, to=2-3]
		\arrow[from=1-1, to=2-2]
		\arrow[from=1-3, to=2-3]
		\arrow["ev"', from=1-2, to=1-1]
		\arrow["\pi", from=1-2, to=1-3]
		\arrow[from=1-2, to=2-2]
		\arrow["\lrcorner"{anchor=center, pos=0.125}, draw=none, from=1-2, to=2-3]
	\end{tikzcd}\]
	where $ev : f^* f_* Z \to Z$ is the co-unit of the adjunction $f^* \dashv f_*$.
	\begin{lem} \label[lem]{desc_f_*}
		In the above situation, the functor $f^Z_* : \dSt/Z \to \dSt/f_* Z$ coincides with the composition $\pi_* \circ ev^*$, where $ev^* : \dSt/Z \to \dSt/f^* f_* Z$ is the pull-back and $\pi_* : \dSt/f^* f_* Z \to \dSt/f_* Z$ is the Weil restriction, see \cref{base_chg_adj}.
	\end{lem}
	\begin{proof}
		Note that $f^Z_*$ is right adjoint to $f_Z^*$, while $\pi_* \circ ev^*$ is right adjoint to $ev_! \circ \pi^*$, so it suffices to check that these left adjoints coincide to conclude. Inspecting the definition of $f_Z^*$ from \cref{def_f_*_Z}, we can immediately rewrite it as the composite $ev_! \circ f_{f^* f_* Z}^*$, where $f_{f^* f_* Z}^* : \dSt/f_* Z \to \dSt/f^* f_* Z$ is the functor mapping a stack $T \to f_* Z$ to $f^* T \to f^* f_* Z$ using the functoriality of $f^*$. Therefore, it actually only suffices to check that $f_{f^* f_* Z}^*$ coincides with $\pi^*$.
		Let thus $T \in \dSt/f_* Z$. Consider the following diagram:
		\[\begin{tikzcd}
			{\pi^* T} & T \\
			{f^* f_* Z} & {f_*Z} \\
			X & Y
			\arrow["f"', from=3-1, to=3-2]
			\arrow[from=2-2, to=3-2]
			\arrow["\pi", from=2-1, to=2-2]
			\arrow[from=2-1, to=3-1]
			\arrow["\lrcorner"{anchor=center, pos=0.125}, draw=none, from=2-1, to=3-2]
			\arrow[from=1-1, to=1-2]
			\arrow[from=1-2, to=2-2]
			\arrow[from=1-1, to=2-1]
			\arrow["\lrcorner"{anchor=center, pos=0.125}, draw=none, from=1-1, to=2-2]
		\end{tikzcd}\]
		As a pasting of two cartesian squares, the exterior square is also cartesian, which precisely means that $\pi^* T \simeq f^* T$ as derived stacks over $f^* f_* Z$, as required.
	\end{proof}
	Let us now introduce the following definition and lemmas.
	\begin{df} \label[df]{df_phi_+}
		Let $\phi : U \to V$ be a proper-schematic and local complete intersection map of derived stacks, and $E \in \Perf(U)$.\footnote{Note that by proper-schematic we mean a map $U \to V$ such that the pull-back along any derived affine $\Spec(A) \to V$ is a proper derived scheme $U_A \to \Spec(A)$. In particular, $\phi$ is a assumed to be a relative derived scheme.}
		Then, by \cite[Thm.~2.1]{LCI}, $\phi_* : \QCoh(U) \to \QCoh(V)$ preserves perfect objects, and thus for $E \in \Perf(U)$ we can set:
		\[\phi_+ E = \phi_*(E^\vee)^\vee \in \Perf(V),\]
		where $-^\vee$ denotes the dual of a perfect sheaf.
	\end{df}
	\begin{rmk}
		Note that $\phi_+$ canonically promotes to a covariant functor $\Perf(U) \to \Perf(V)$, which is moreover left adjoint to $\phi^* : \Perf(V) \to \Perf(U)$. As in \cite[Prop.~1.4]{LCI}, also note that $\phi_*$ satisfies the base-change formula for pull-back of affines $\Spec(A) \to V$ (more generally for quasi-compact quasi-separated derived schemes $X \to V$).
	\end{rmk}
	\begin{lem} \label[lem]{lem_phi_*} \label[lem]{phi_*_perf_lin_st}
		In the above setting, there is a canonical equivalence $\phi_* \bb V(E) \simeq \bb V(\phi_+ E)$ as $\bb G_m$-equivariant derived stacks over $V$.
	\end{lem}
	\begin{proof}
		We will compare their functors of points as $\bb G_m$-equivariant derived stacks over $V$, where the $\bb G_m$-action on $V$ is trivial. Let $v : T \to V$ be any map of $\bb G_m$-equivariant derived stacks, where $T$ is affine. We can form the following pull-back diagram:
		\[\begin{tikzcd}
			{\phi^* T} & T \\
			U & V
			\arrow["\phi"', from=2-1, to=2-2]
			\arrow["v", from=1-2, to=2-2]
			\arrow["p", from=1-1, to=1-2]
			\arrow["u"', from=1-1, to=2-1]
			\arrow["\lrcorner"{anchor=center, pos=0.125}, draw=none, from=1-1, to=2-2]
		\end{tikzcd}\]
		Notice that $\phi^* T$ inherits a canonical $\bb G_m$-action, and $p$ is $\bb G_m$-equivariant, because $\phi : U \to V$ is $\bb G_m$-equivariant for the trivial actions on both $U$ and $V$. This means that $\mcal O_T$ and $\mcal O_{\phi^* T}$ both are equipped with a weight grading and they are compatible under $p$, i.e. $p^* \mcal O_T \simeq \mcal O_{\phi^* T}$ as graded algebras in $\QCoh(\phi^* T)$. Freely using the fact that pull-back functors of quasi-coherent sheaves are monoidal and thus commute with $-^\vee$, we then have:
		\begin{align*}
			&\Map_{\dSt/V - \bb G_m}(T,\phi_* \bb V(E)) \\
			\simeq\ &\Map_{\dSt/U - \bb G_m}(\phi^*T,\bb V(E)) \\
			=\ &\Map_{\QCoh(\phi^* T)}(u^* E,\mcal O_{\phi^* T}(1)) \\
			\simeq\ &\Map_{\QCoh(\phi^* T)}(\mcal O_{\phi^* T}(1)^\vee,u^*(E^\vee)) \\
			\simeq\ &\Map_{\QCoh(\phi^* T)}(p^* \mcal O_T(1)^\vee,u^*(E^\vee)) \\
			\simeq\ &\Map_{\QCoh(T)}(\mcal O_T(1)^\vee,p_* u^*(E^\vee)) \\
			\simeq\ &\Map_{\QCoh(T)}(\mcal O_T(1)^\vee,v^* \phi_*(E^\vee)),\textup{ using the base-change formula} \\
			\simeq\ &\Map_{\QCoh(T)}(v^*(\phi_+ E),\mcal O_T(1)) \\
			=\ &\Map_{\dSt/V - \bb G_m}(T,\bb V(\phi_+ E))
		\end{align*}
		Note that the base-change formula applies since $T$ is affine and $\phi$ is proper-lci. Moreover, this equivalence is natural in $T \to V$, which concludes the proof.
	\end{proof}
	\begin{rmk}
		This statement can be interpreted in the context of \cref{rmk_f_beck_chevalley}, where conditions on a map were compared to certain canonical adjunction squares to be Beck-Chevalley. In this remark, on the side of derived stacks, only the adjunction $f_! \dashv f^*$ was considered, while on this side there is always a canonical triple $f_! \dashv f^* \dashv f_*$, as in \cref{base_chg_adj}. This is because, in general, there is no corresponding triple on the side of quasi-coherent sheaves, and only $f^* \dashv f_*$ can be considered on that side. However, in the above situation, the functor $\phi_+$ does fit in a similar triple $\phi_+ \dashv \phi^* \dashv \phi_*$, and the above statement can be interpreted as a proof that the appropriate extra canonical square that exists in this situation is in fact Beck-Chevalley too, at least when restricted to symmetric algebras. We will neither formally state nor prove this fact in this paper as the above result suffices.
	\end{rmk}
	Collecting what was proven so far, we obtain the following corollary.
	\begin{cor} \label[cor]{pushforward_cotangent_formula}
		For $f : X \to Y$ a proper-schematic and local complete intersection map between derived stacks, $Z \to X$ an arbitrary map of derived stack, and $E \in \Perf(Z)$, there is a canonical equivalence of $\bb G_m$-equivariant stacks over $f_*Z$:
		\[f^Z_* \bb V(E) \simeq \bb V(\pi_+ ev^* E),\]
		where $ev : f^* f_* Z \to Z$ is the co-unit of $f^* \dashv f_*$ at $Z$, and $\pi : f^* f_* Z \to f_* Z$ is the projection.
		In particular, $f^Z_*$ sends perfect linear stacks over $Z$ to perfect linear stacks over $f_* Z$.
	\end{cor}
	\begin{proof}
		As a consequence of \cref{desc_f_*}, $f^Z_* \bb V(E) \simeq \pi_* ev^* \bb V(E)$ as $\bb G_m$-equivariant stacks over $f_* Z$.
		Moreover, $ev^* \bb V(E) \simeq \bb V(ev^* E)$ as $\bb G_m$-equivariant stacks over $f^* f_* Z$, as in \cref{base_chg_adj}.
		Since $E$ is perfect and $ev^*$ is monoidal, $ev^* E$ is also perfect, and since $f$ is proper-schematic and local complete intersection and $\pi$ is its pull-back, $\pi$ is also proper-schematic and local complete intersection. We can thus apply \cref{phi_*_perf_lin_st}, and obtain an equivalence $\pi_* \bb V(ev^* E) \simeq \bb V(\pi_+ ev^* E)$ as $\bb G_m$-equivariant stacks over $f_* Z$.
		Combining these three equivalences, we obtain as required an equivalence $f^Z_* \bb V(E) \simeq \bb V(\pi_+ ev^* E)$ as $\bb G_m$-equivariant stacks over $f_* Z$.
	\end{proof}
	Finally, using the equivalence $f^Z_* \mcal L^{gr}(Z/X) \simeq \mcal L^{gr}(f_* Z/Y)$ of $\mcal H$-equivariant stacks over $Y$ and $\bb G_m$-equivariant stacks over $f_* Z$ provided by \cref{f_*_lgr}, we obtain our main construction.
	\begin{cor} \label[cor]{pushforward_preserves_fols}
		For $f : X \to Y$ a proper-schematic and local complete intersection map between derived stacks, and $Z \to X$ an arbitrary map of derived stack, the direct image functor $f_* : \dSt/X \to \dSt/Y$ induces a canonical functor:
		\[f_* : \PerfLinSt_Z/\mcal L^{gr}(Z/X) - \mcal H \to \PerfLinSt_{f_* Z}/\mcal L^{gr}(f_* Z/Y) - \mcal H\]
	\end{cor}
	\begin{proof}
		Let $F \in \PerfLinSt_Z/\mcal L^{gr}(Z/X) - \mcal H$. Firstly, we will equip $f_* F$ with an $\mcal H$-action over $f_* Z$, an $\mcal H$-equivariant map to $\mcal L^{gr}(f_* Z/Y)$, and show that it is perfect linear over $f_* Z$.
		Using the monoidal structure of $f^Z_*$ in the sense of \cref{pushforward_monoidal}, we can equip $f_* F$ with an $\mcal H$-action over $f_* Z$, inherited from the structural one on $F$ over $Z$.
		Similarly, the structural $\mcal H$-equivariant map $F \to \mcal L^{gr}(Z/X)$ is mapped to an $\mcal H$-equivariant map $f^Z_* F \to f^Z_* \mcal L^{gr}(Z/X)$. Moreover, by \cref{f_*_lgr}, there is a canonical equivalence $f^Z_* \mcal L^{gr}(Z/X) \simeq \mcal L^{gr}(f_* Z/Y)$ as $\mcal H$-equivariant stacks, which composed with the previous map produces an $\mcal H$-equivariant map $f^Z_* F \to \mcal L^{gr}(f_* Z/Y)$.
		Finally, by \cref{pushforward_cotangent_formula}, $f^Z_* F$ is perfect linear over $f_* Z$.
		As a consequence we have constructed a functor:
		\[f_* : \PerfLinSt_Z/\mcal L^{gr}(Z/X) - \mcal H \to \PerfLinSt_{f_* Z}/\mcal L^{gr}(f_* Z/Y) - \mcal H\]
	\end{proof}
	\subsection{Pushforward of Derived Foliations}
	At this point we want to apply \cref{equiv_fol_dst_rel_dm} to replace both sides of the above functor with $\mcal Fol(Z/X)$ and $\mcal Fol(f_*Z/Y)$ respectively to conclude. We will thus assume that $Z \to X$ is relative derived Deligne-Mumford stack locally of finite presentation to apply it on the source, but we will need the following lemma to show that the target also satisfies this condition.
	\begin{lem} \label[lem]{f_*_rel_dm}
		Let $Z \to X$ be a relative derived Deligne-Mumford stack locally of finite presentation, and $f : X \to Y$ a proper-schematic, flat, and local complete intersection map of derived stacks. Assume that $f_* Z \to Y$ is a relative derived Artin stack locally of finite presentation.
		Then, $f_* Z \to Y$ is a relative derived Deligne-Mumford stack locally of finite presentation.
	\end{lem}
	\begin{proof}
		Since $f_* Z \to Y$ is assumed to be a relative derived Artin stack and locally of finite presentation, for it to be a relative derived Deligne-Mumford stack locally of finite presentation, it suffices to show that its truncation $t_0(f_* Z) \to t_0(Y)$ is a relative (non-derived) Deligne-Mumford stack.
		Let us now consider the truncation:
		\[\begin{tikzcd}
			{t_0(Z)} & {t_0(f_*Z)} \\
			{t_0(X)} & {t_0(Y)}
			\arrow[from=1-1, to=2-1]
			\arrow[from=1-2, to=2-2]
			\arrow["{t_0(f)}"', from=2-1, to=2-2]
		\end{tikzcd}\]
		Since $Z \to X$ is a relative derived Deligne-Mumford stack, $t_0(Z) \to t_0(X)$ is a relative (non-derived) Deligne-Mumford stack, and similarly $t_0(f)$ is a proper-schematic, flat, and local complete intersection map of (non-derived) stacks, so we can apply \cite[Thm.~1.5]{olsson2006hom}, and obtain that $t_0(f)_* t_0(Z) \to t_0(Y)$ is a relative (non-derived) Deligne-Mumford stack. Finally, since $f$ is schematic and flat, let us now show that $t_0(f)_* t_0(Z)$ is equivalent to $t_ 0(f_* Z)$, which will prove that $f_* Z \to Z$ is a relative derived Deligne-Mumford stack.
		Let $\Spec(A)$ be a non-derived affine scheme, we have the following chain of equivalences, natural in $\Spec(A)$:
		\begin{align*}
			& \Map_{\textup{St}}(\Spec(A),t_0(f_* Z)) \\
			&\simeq \Map_{\dSt}(\Spec(A),f_* Z) \\
			&\simeq \Map_{\dSt}(f^* \Spec(A),Z) \\
			&\simeq \Map_{\textup{St}}(t_0(f)^* \Spec(A),t_0(Z))\ \ \ \ \ \ (*) \\
			&\simeq \Map_{\textup{St}}(\Spec(A),t_0(f)_* t_0(Z)) \\
		\end{align*}
		where the equivalence $(*)$ holds since, as $f$ is schematic, $f^* \Spec(A)$ is a derived scheme; moreover, as $f$ is flat, this derived scheme is in fact non-derived and coincides with $t_0(f)^* \Spec(A)$; and finally maps from non-derived schemes to derived stacks factor through their truncations.
	\end{proof}
	Combining these results, we obtain our push-forward of derived foliations.
	\restate{thm}{thm_pushforward_fol}
	\begin{proof}
		Applying \cref{pushforward_preserves_fols}, we have a functor:
		\[f_* : \PerfLinSt_Z/\mcal L^{gr}(Z/X) - \mcal H \to \PerfLinSt_{f_* Z}/\mcal L^{gr}(f_* Z/Y) - \mcal H\]
		Moreover, by \cref{equiv_fol_dst_rel_dm}, we have that $\mcal Fol(Z/X) \simeq \PerfLinSt_Z/\mcal L^{gr}(Z/X) - \mcal H$ since $Z \to X$ is a relative derived Deligne-Mumford stack locally of finite presentation. We also have that $\mcal Fol(f_* Z/Y) \simeq \PerfLinSt_{f_* Z}/\mcal L^{gr}(f_* Z/Y) - \mcal H$ since $f_* Z \to X$ is also a relative derived Deligne-Mumford stack locally of finite presentation by \cref{f_*_rel_dm}.
		Therefore, using these equivalences, we obtain the push-forward functor:
		\[f_{*,\mcal Fol} : \mcal Fol(Z/X) \to \mcal Fol(f_* Z/Y)\]
	\end{proof}
	\begin{rmk}
		Being a push-forward, we can wonder whether this functor is a right adjoint. More precisely, remember that this functor is essentially the application of $f_*^Z : \dSt/Z \to \dSt/f_* Z$, which is right adjoint to the functor $f^*_Z : \dSt/f_*Z \to \dSt/Z$, so we can ask if this left adjoint also preserves perfect linear stacks and $\mcal H$-actions and defines a left adjoint to $f_{*,\mcal Fol}$. Note that the pull-back of derived foliations as defined in \cref{pull-back} does not apply in this context and this would thus be a different construction, since there is no canonical map $Z \to f_* Z$ making a commutative square with the existing maps $Z \to X$, $X \to Y$ and $f_* Z \to Y$. However, we do not expect $f^*_Z$ to preserve neither perfect linear stacks nor their $\mcal H$-actions in general or under reasonable conditions. This is because, by \cref{desc_f_*} it can be rewritten as $ev_! \circ \pi^*$, where $\pi : f^* f_* Z \to f_* Z$ is the projection, and $ev : f^* f_* Z \to Z$ is the evaluation map (the co-unit of the adjunction $f^* \dashv f_*$). While $\pi^*$ always preserves this structure, we do not expect $ev_!$ to preserve it under reasonsable conditions.
	\end{rmk}
	As a final application, we will now prove the main \cref{main_theorem}.
	\restate{thm}{main_theorem}
	\begin{proof}
		Let us denote by $f : X \to S$ the structural map of $X$ over $S$, which by assumption is proper-schematic, flat, and local complete intersection. Consider the derived stack $Z = X \times_S Y$ over $X$, and note that since $Y$ is a relative derived Deligne-Mumford stack locally of finite presentation over $S$, $Z$ is a relative derived Deligne-Mumford stack locally of finite presentation over $X$. Now recall that $f_* Z \simeq \underline{\Map}_S(X,Y)$ as derived stacks over $S$. To see this, let us compute its functor of points. For $\Spec(A) \to S$, we have the following equivalences, natural in $\Spec(A) \to S$:
		\begin{align*}
			\Map_{\dSt/S}(\Spec(A),f_* Z) &\simeq \Map_{\dSt/X}(f^* \Spec(A),Z)\\
			&= \Map_{\dSt/X}(f^* \Spec(A),X \times_S Y)\\
			&= \Map_{\dSt/X}(f^* \Spec(A),f^* Y)\\
			&\simeq \Map_{\dSt/S}(f_! f^* \Spec(A),Y)\\
			&= \Map_{\dSt/S}(\Spec(A) \times_S X,Y)\\
			&\simeq \Map_{\dSt/S}(\Spec(A),\underline{\Map}_S(X,Y))
		\end{align*}
		As the functors of point coincide, $f_* Z \simeq \underline{\Map}_S(X,Y)$ as derived stacks over $S$. Therefore, the assumption that $\underline{\Map}_S(X,Y)$ is a relative derived Artin stack locally of finite presentation is exactly the assumption that $f_* Z$ is a relative derived Artin stack locally of finite presentation.
		As such, the conditions of \cref{thm_pushforward_fol} are satisfied, and there is a canonical push-forward functor $f_{*,\mcal Fol} : \mcal Fol(Z/X) \to \mcal Fol(f_* Z/S) = \mcal Fol(\underline{\Map}_S(X,Y)/S)$.
		Let now $\mcal F \in \mcal Fol(Y/S)$ denote the derived foliation that $Y$ comes equipped with, and finally we denote by $p_Y : X \times_S Y \to Y$ the projection on $Y$. Note that it forms a commutative square with $f : X \to S$. We will denote by $\phi : (X \times_S X = Z \to X) \to (Y \to S)$ this commutative square. Since both $Y \to S$ and $Z \to X$ have perfect cotangents as they are locally of finite presentation, this induces a canonical pull-back functor $\phi^* : \mcal Fol(Y/S) \to \mcal Fol(Z/X)$, by a computation similar to \cref{formula_pull-back}. 
		Therefore, we obtain a derived foliation on $\underline{\Map}_S(X,Y)$ by setting:
		\[\mcal F_{\Map} := f_{*,\mcal Fol}(p_Y^* \mcal F) \in \mcal Fol(\underline{\Map}_S(X,Y)/S) \qedhere\]
	\end{proof}
	We will also prove that the induced derived foliation satisfies the property mentioned in the introduction.
	\restate{prop}{main_prop}
	\begin{proof}
		We will first prove the following equivalence:
		\[\bb T_{\mcal F,\Map} \simeq \pi_* ev^* \bb T_{\mcal F}\]
		Note that by definition in \cref{main_theorem}, we have:
		\[\mcal F_{\Map} := f_{*,\mcal Fol}(p_Y^* \mcal F)\]
		Moreover, by \cref{formula_pull-back}, if $\bb L_{\mcal F}$ is the cotangent of $\mcal F$, the cotangent of $p_Y^* \mcal F$ is given by:
		\[\bb L_{p_Y^* \mcal F} \simeq \bb L_{X \times_S Y/X} \oplus_{p_Y^* \bb L_{Y/S}} p_Y^* \bb L_{\mcal F}\]
		Note that $\bb L_{X \times_S Y/X} \simeq p_Y^* \bb L_{Y/S}$, since $\bb L_{X \times_S Y/S} \simeq p_Y^* \bb L_{Y/S} \oplus p_X^* L_{X/S}$, where $p_X : X \times_S Y \to X$ is the projection, and $\bb L_{X \times_S Y/X}$ is defined as the quotient of $\bb L_{X \times_S Y/S}$ by $p_X^* \bb L_{X/S}$. As such, we obtain that:
		\[\bb L_{p_Y^* \mcal F} \simeq p_Y^* \bb L_{\mcal F}\]
		As noted in \cref{formula_pull-back}, in this case we would also obtain that $\bb L_{p_Y^* \mcal F}$ happens to remain perfect since $p_Y^*$ preserves perfects, without the need for further assumptions. More generally, pull-back of non-perfect derived foliations on cartesian square will always satisfy this relation for their cotangent, and thus preserve perfect derived foliations in these cases.
		
		Applying the formula for the cotangent of the derived foliation induced by $f_{*,\mcal Fol}$, which is provided by \cref{pushforward_cotangent_formula}, we obtain that:
		\[\bb L_{\mcal F,\Map} \simeq \pi_*(ev'^* p_Y^* \bb L_{\mcal F}^\vee)^\vee\]
		where $ev' : X \times_S \underline{\Map}_S(X,Y) \to X \times_S Y$ is the co-unit of the adjunction $f^* \dashv f_*$ at $Z = X \times_S Y$. Note that $p_Y \circ ev' \simeq ev$, so that $ev'^* \circ p_Y^* \simeq ev^*$. Therefore, after dualizing, we do obtain the expected formula:
		\[\bb T_{\mcal F,\Map} \simeq \pi_* ev^* \bb T_{\mcal F}\]
		We will now show that this formula induces the one at an $S$-point of $\underline{\Map}_S(X,Y)$. Let $u : S \to \underline{\Map}_S(X,Y)$ be an $S$-linear map, corresponding to an $S$-linear map $g : X \to Y$, as can be seen in the following commutative diagram:
		\[\begin{tikzcd}
			& X & S \\
			Y & {X \times_S \underline{\Map}_S(X,Y)} & {\underline{\Map}_S(X,Y)} \\
			& X & S
			\arrow["f", from=1-2, to=1-3]
			\arrow["g"', from=1-2, to=2-1]
			\arrow["{u'}", from=1-2, to=2-2]
			\arrow["\lrcorner"{anchor=center, pos=0.125}, draw=none, from=1-2, to=2-3]
			\arrow["u", from=1-3, to=2-3]
			\arrow["ev", from=2-2, to=2-1]
			\arrow["\pi", from=2-2, to=2-3]
			\arrow[from=2-2, to=3-2]
			\arrow["\lrcorner"{anchor=center, pos=0.125}, draw=none, from=2-2, to=3-3]
			\arrow[from=2-3, to=3-3]
			\arrow["f"', from=3-2, to=3-3]
		\end{tikzcd}\]
		Note that the square at the bottom right is a pull-back by definition, and the outer square is a pull-back since the compositions $S \to \underline{\Map}_S(X,Y) \to S$ and $X \to X \times_S \underline{\Map}_S(X,Y) \to X$ are the identity maps on $S$ and $X$ respectively, so by pasting the top square is a pull-back too.
		Also, by definition we have that $\bb T_g \mcal F_{\Map} := u^* \bb T_{\mcal F,\Map}$, so by the previous equivalence we have:
		\[\bb T_g \mcal F_{\Map} \simeq u^* \pi_* ev^* \bb T_{\mcal F}\]
		Moreover, since $\pi$ is proper-schematic as it is the pull-back of $f$ which is proper-schematic by assumption, the base-change formula holds, so $u^* \pi_* \simeq f_* u'^*$ for all quasi-coherent sheaves, in particular we obtain:
		\[\bb T_g \mcal F_{\Map} \simeq f_* u'^* ev^* \bb T_{\mcal F}\]
		Notice that, in our notation for this proposition, the functor $f_* : \QCoh(X) \to \QCoh(S)$ is denoted by $\Gamma_S(X,-)$, and that the composition $ev \circ u'$ is equivalent to $g$, so we obtain the desired formula:
		\[\bb T_g \mcal F_{\Map} \simeq \Gamma_S(X,g^* \bb T_{\mcal F})\]
	\end{proof}
	\begin{rmk} \label[rmk]{tangent_pushforward_fol}
		The same reasoning also applies to describe the tangent of the derived foliation induced by $f_{*,\mcal Fol}$ in \cref{thm_pushforward_fol}. In the setting of this theorem, if $\mcal F$ is a derived foliation on $Z$ relative to $S$, and $\mcal F_{\mathbf{Res}}$ is its image by $f_{*,\mcal Fol}$, we have:
		\[\bb T_{\mcal F,\mathbf{Res}} \simeq \pi_* ev^* \bb T_{\mcal F}\]
		where, in this case, we denote by $ev : f^* f_* Z \to Z$ the co-unit of $f^* \dashv f_*$ itself, and $\pi : f^* f_* Z \to Z$ is the projection. Moreover, for $u : S \to f_* Z$ an $S$-point of $f_* Z$, corresponding to an $S$-linear section $g : X \to Z$ of the structural map $Z \to X$, we have:
		\[\bb T_g \mcal F_{\mathbf{Res}} \simeq \Gamma_S(X,g^* \bb T_{\mcal F})\]
		Note that in this case, the map $g : X \to Z$ is defined in the following commutative diagram:
		\[\begin{tikzcd}
			& X & S \\
			Z & {f^* f_* Z} & {f_* Z} \\
			& X & S
			\arrow["f", from=1-2, to=1-3]
			\arrow["g"', from=1-2, to=2-1]
			\arrow["{u'}", from=1-2, to=2-2]
			\arrow["\lrcorner"{anchor=center, pos=0.125}, draw=none, from=1-2, to=2-3]
			\arrow["u", from=1-3, to=2-3]
			\arrow["{Id_S}", shift left=5, equals, from=1-3, to=3-3]
			\arrow[from=2-1, to=3-2]
			\arrow["ev"', from=2-2, to=2-1]
			\arrow["\pi", from=2-2, to=2-3]
			\arrow[from=2-2, to=3-2]
			\arrow["\lrcorner"{anchor=center, pos=0.125}, draw=none, from=2-2, to=3-3]
			\arrow[from=2-3, to=3-3]
			\arrow["f"', from=3-2, to=3-3]
		\end{tikzcd}\]
	\end{rmk}
	\begin{rmk}
		While we have provided a description of the tangent and thus also of the cotangent complex of a derived foliation induced by \cref{main_theorem} or \cref{thm_pushforward_fol}, we do not provide a description of the induced mixed differential on the corresponding graded mixed algebra. It should be noted first that the mixed differential comes from the $\mcal H$-action on the perfect linear stack corresponding to the induced derived foliation, so the mixed differential would be computed via this action. However, it should also be noted that the explicit computation of mixed differentials for derived foliations is very difficult, as discussed in \cite[remark preceding Def.~2.1.1.3]{toen2025derivedfoliations}. This is because the graded mixed algebra of a derived foliation is only quasi-isomorphic to the symmetric algebra of its cotangent. As such, while the mixed differential of the graded mixed algebra itself may have an explicit description, its homotopical transport along the quasi-isomorphism only induces a weaker structure on the symmetric algebra of its cotangent. This structure is not the data of a differential whose square is equal to zero, but the data of a differential equipped with a homotopy from its square to zero, and moreover with infinitely many higher homotopy coherences. The explicit description of this homotopical data is outside the scope of this paper, but would be required to compute the resulting differential on the derived foliations induced by our theorems. Nonetheless, it is often the case that derived foliations first arise from their tangent bundles rather than their cotangent bundles, in which case the differential on the cotangent bundle corresponds to the dual of the Lie bracket on the tangent of the derived foliation, which moreover often is the restriction of the Lie bracket of the ambient space to this tangent bundle. In those cases, this provides an intuitive description of the mixed differential, which we will not make more precise.
	\end{rmk}
	\section{Applications} \label[section]{applications}
	Let us now prove and discuss some applications of this paper.
	\subsection{Foliations for Gromov-Witten Theory}
	Let us start by providing more details to \cref{example_application_1}.
	\begin{reminder}
		As a first application, let us consider the moduli stacks of pre-stable and stable curves of genus $g$ with $n$ marked points. We take the following constructions and properties from \cite[Section 2.2]{schurg2015derived}.
		\begin{itemize}
			\item There is a derived stack $\bb R \mathbf{M}_{g,n}^{pre}$, called the derived moduli stack of pre-stable curves of genus $g$ with $n$ marked points, that comes equipped with a universal family $\bb R \mcal C_{g,n}^{pre}$ defined over $\bb R \mathbf{M}_{g,n}^{pre}$, whose classical truncations coincide with the usual $\mathbf{M}_{g,n}$ and $\mcal C_{g,n}$ respectively.
			\item As such, for $X$ a smooth (possibly derived) projective scheme, we can define the derived stack $\bb R \mathbf{M}_{g,n}^{pre}(X)$ as the following derived mapping stack, canonically defined over $\bb R \mathbf{M}_{g,n}^{pre}$:
			\[\bb R \mathbf{M}_{g,n}^{pre}(X) := \underline{\Map}_{\bb R \mathbf{M}_{g,n}^{pre}}(\bb R \mcal C_{g,n}^{pre},X \times \bb R \mathbf{M}_{g,n}^{pre})\]
			Note that points of this stack morally correspond to families of pre-stable curves of genus $g$ with $n$ marked points in $X$. Its truncation recovers the usual non-derived moduli stack $\mathbf{M}_{g,n}^{pre}(X)$ of pre-stable curves of genus $g$ with $n$ marked points. As such, it is Artin locally of finite presentation, and thus $\bb R \mathbf{M}_{g,n}^{pre}(X)$ is derived Artin locally of finite presentation too. See \cite[Cor.~2.4(ii),Rmk.~2.5]{HR15} and \cite[Prop.~2]{B97} to obtain this last fact on the classical stack.
			\item Moreover, $\bb R \mathbf M_{g,n}^{pre}(X)$ contains an open sub-stack $\bb R \overline{\mathbf M}_{g,n}(X;\beta)$, which constitutes similarly the derived moduli stack of families of stable curves of genus $g$ with $n$ marked points and of class $\beta$ on $X$, and whose classical truncation also coincides with the usual $\overline{\mathbf{M}}_{g,n}(X;\beta)$ non-derived moduli stack of such objects. As it is an open substack of the previous one, it is also Artin locally of finite presentation, and thus $\bb R \overline{\mathbf M}_{g,n}(X;\beta)$ is too.
			\item Finally, note that the map $\bb R \mcal C_{g,n}^{pre} \to \bb R \mathbf M_{g,n}^{pre}$ is proper-schematic (i.e.\ representable by proper derived schemes), flat, and local complete intersection.
		\end{itemize}
	\end{reminder}
	\begin{cor} \label[cor]{foliation_Mgn}
		Let $X$ be a smooth (possibly derived) projective scheme equipped with a derived foliation. Then $\bb R \mathbf M_{g,n}^{pre}(X)$ and $\bb R \overline{\mathbf M}_{g,n}(X;\beta)$ inherit derived foliations.
	\end{cor}
	\begin{proof}
		The derived foliation on $X$ induces one on $X \times \bb R \mathbf{M}_{g,n}^{pre}$ by pull-back along the projection to $X$. The derived foliation on $\bb R \mathbf M_{g,n}^{pre}(X)$ is then obtained as a direct application of \cref{main_theorem} to this derived foliation, as all the conditions of the theorem are satisfied in this case by the above reminder. It also induces one on $\bb R \overline{\mathbf M}_{g,n}(X;\beta)$ by restriction (i.e.\ pull-back) along the inclusion.
	\end{proof}
	\begin{rmk}
		In the study of Gromov-Witten invariants, where the new data of a derived foliation could lead to the definition of new invariants relative to it, the space $\bb R \overline{\mathbf M}_{g,n}(X;\beta)$ is typically considered when $X$ is a smooth projective scheme, which is why this result was stated in this case. While our construction would also apply when $X$ is a derived Deligne-Mumford stack locally of finite presentation, Gromov-Witten invariants on such an $X$ can involve different definitions and constructions, which we will not discuss in this paper.
	\end{rmk}
	\subsection{Foliations on Hilbert Schemes}
	We will now provide more details to \cref{example_application_2}.
	\begin{reminder} \label[reminder]{def_Hilb}
		As a second application, let us now rephrase the construction of $\bb R \mathbf{Hilb}^{lci}(Y)$ from \cite{ciocan2002derived} in terms of mapping stacks. See also the last example of \cite[\S3.2]{toenDAGSurvey} which provides exactly the same construction in a slightly different presentation.
		
		Let $Sch_{prop,flat,lci} : \dAff^\op \to \infty-Grpd$ be the functor defined by:
		\[Sch_{prop,flat,lci}(S) := \{\textup{derived schemes } Z \to S \textup{ that are proper, flat and lci over } S\}\]
		This functor is a derived stack, and moreover it is equipped with a universal object $\mcal S \to Sch_{prop,flat,lci}$ which is a relative derived scheme, proper-schematic, flat, and local complete intersection. Let $Y$ be a derived Deligne-Mumford stack locally of finite presentation, and consider the derived stack:
		\[\underline{\Map}_{Sch_{prop,flat,lci}}(\mcal S,Y \times Sch_{prop,flat,lci})\]
		As mentioned in the last example of \cite[\S3.2]{toenDAGSurvey}, this is a derived Artin stack locally of finite presentation over $Sch_{prop,flat,lci}$.
		Note that if $Y$ is equipped with a derived foliation, following a proof analogous to \cref{foliation_Mgn}, since the conditions of \cref{main_theorem} also apply in this case, this derived stack thus inherits a derived foliation. Moreover, for a derived affine $S$, the $S$-points of this derived stack are equivalent to the data of a derived $S$-scheme $Z$ which is proper, flat, and local complete intersection over $S$, and moreover equipped with a map $Z \to Y$. As such, it contains $\bb R \mathbf{Hilb}^{lci}(Y)$, the derived Hilbert stack, whose $S$-points consist in those $Z \to S$ such that the equipped map $Z \to Y$ is moreover a closed immersion. This substack also inherits a derived foliation by restriction.
	\end{reminder}
	\begin{cor} \label[cor]{foliation_Hilb}
		Let $Y$ be a derived Deligne-Mumford stack locally of finite presentation, equipped with a derived foliation. Then the derived moduli stack of closed derived subschemes of $Y$ that are flat, proper, and local complete intersection over their base, denoted by $\bb R \mathbf{Hilb}^{lci}(Y)$, inherits a derived foliation.
	\end{cor}
	As mentioned in the introduction, we expect integral points of this inherited derived foliation in a sense analogous to \cite[Def.~1.9.2(1)]{EDSbook}, i.e.\ the locus where its cotangent vanishes, to correspond to closed derived subschemes of $Y$ contained in the algebraic leaves of its derived foliation. A rigorous statement and proof of this result is outside the scope of this paper, but will be studied instead in a future work, as mentioned in the next section.
	\separation
	More generally, these two examples of applications can be thought of as illustrations of the transfer of derived foliations to a wide class of derived moduli stacks.
	\section{Ongoing and Future Work} \label[section]{future_work}
	Let us now talk about ongoing and future work following this paper. The first subsection of this section concerns the ongoing work, while the next subsections concern future work.
	
	Note that in the author's thesis \cite[Chapter~7]{Alfieri2026ShiftedSymplecticFoliations}, shifted symplectic structures on derived foliations were defined and studied. More specifically, in the same context as this paper, shifted symplectic forms on a derived foliation were defined, as well as Lagrangian structures for derived foliations in shifted symplectic stacks. Then, a mapping theorem/AKSZ-construction for the former was proven, as well as an intersection theorem for the latter. These results will be presented in a self-contained paper.
	\subsection{Extension beyond the Deligne-Mumford Case}
	In an ongoing work, we expect the results of this paper to extend to the Artin case beyond the Deligne-Mumford case. In fact, in an ongoing work, we expect the right assumption to be that $Z \to X$ is any map of derived pre-stacks which admits a relative cotangent $\bb L_{Z/X} \in \QCoh(Z)$ which is moreover perfect, and we will say that the map $Z \to X$ admits a perfect relative cotangent when it is the case. Note that this is the case for example when $Z \to X$ is a relative derived Artin stack locally of finite presentation, but we expect our results to hold even without the derived Artin assumption. We can also hope to extend the results when the map $f : X \to Y$ along which we push forward is not representable, and in particular it seems that a relative version of the condition of $\mcal O$-compactness of \cite[Def.~2.1]{PTVV} is sufficient. Note similarly that this condition is implied by the current proper-schematic and lci assumption. In particular, in this situation we should be able to remove the flat assumption on $f$ entirely, as in this paper it was only needed to ensure that $f_* Z \to Y$ is again derived Deligne-Mumford, but we will no longer need this condition. As such, we expect to prove the following theorems.
	\begin{expected}
		Let $Z \to X$ be a map of derived pre-stacks which admits a perfect relative cotangent, and $f : X \to Y$ an $\mcal O$-compact map of derived stacks. Then there is a push-forward functor:
		\[f_{*,\mcal Fol} : \mcal Fol(Z/X) \to \mcal Fol(f_*Z / Y)\]
	\end{expected}
	We will now discuss future work.
	\subsection{Exterior Differential Systems}
	In differential geometry, exterior differential systems can be seen as generalizations of foliations. Foliations correspond to partial differential equation systems that only involve $1$-forms called Pfaffian systems, while exterior differential systems can involve $p$-forms for multiple $p \geq 1$, with finitely many generators and relations. In this context, obstructions to the integrability of an exterior differential system are measured by Spencer cohomology, see \cite{EDS}. In a work in preparation, Eyssidieux and Kwon study exterior differential systems in classical algebraic geometry, and plan to express the moduli space of their solutions as a subspace of the Hilbert scheme of closed subschemes.
	
	In a future work, we will study the definition of a derived enhancement of exterior differential systems. Essentially, for a map of derived stacks $Z \to X$, an exterior differential system on $Z$ relative to $X$ is a graded mixed algebra $F \in \epsilon-cdga^{gr}(Z/X)$ which is a compact object in the category of graded algebras, a condition which expresses the fact that $F$ has a finite presentation in some sense. We plan to prove a mapping theorem for such structures by extending the above push-forward of derived foliations to such algebras, and investigate the connection between the dg-Lie algebroid of the derived enhancement of the moduli space of their solutions and Spencer cohomology. In particular, a mapping theorem for exterior differential systems should allow us to induce an exterior differential system on the derived Hilbert scheme constructed above in \cref{def_Hilb}, and we should be able to relate integral points of this system to the moduli space of solutions of the initial system. Moreover, we can study the formal integration of such systems in derived geometry, possibly using the formal integration results from \cite{brantner2025formalintegrationderivedfoliations}.
	\appendix
	\section{Beck-Chevalley Conditions} \label[appendix]{beck-chevalley-section}
	In \cref{equiv-affine-case}, we use the formalism of Beck-Chevalley squares to obtain the equivalence between derived foliations and perfect linear stacks in the affine case. This formalism originates from base-change formulas in algebraic geometry, which appeared first in \cite[(6.1.1)]{EGAIII2}, which are particular instances of Beck-Chevalley conditions that also reappear here in \cref{beck-chevalley-dag}. The term itself was only coined later by Bénabou and Roubaud in 1970, which discussed the property in \cite{BeckChevalleyDef}, referring to Beck and Chevalley who also worked on the notion around the same time without publishing about it. In the modern context of $\infty$-category theory, the notion was defined by Lurie in \cite[Def.~4.7.4.13]{HA}, which he also refers to as adjointable squares. In this appendix, we will recall these definitions and basic properties about them in \cref{beck-chevalley-defs}, and their appearance in our work in derived algebraic geometry as well as sufficient conditions for them to hold in \cref{beck-chevalley-dag}.
	\subsection{Generalities on Beck-Chevalley Conditions} \label[subsection]{beck-chevalley-defs}
	In this subsection, we will focus on recalling the definition of the Beck-Chevalley condition in $\infty$-category theory and results following \cite[Def.~4.7.4.13]{HA}.
	\begin{reminder} \label[reminder]{def_pre_square}
		Consider a diagram of $\infty$-categories of the following form (which is not assumed to be commutative in any way):
		\[\begin{tikzcd}
			{\mcal C'} & {\mcal D'} \\
			{\mcal C} & {\mcal D}
			\arrow["U"', from=1-1, to=2-1]
			\arrow["V", from=1-2, to=2-2]
			\arrow[""{name=0, anchor=center, inner sep=0}, "{R'}"', shift right, from=1-1, to=1-2]
			\arrow[""{name=1, anchor=center, inner sep=0}, "R"', shift right, from=2-1, to=2-2]
			\arrow[""{name=2, anchor=center, inner sep=0}, "{L'}"', shift right, from=1-2, to=1-1]
			\arrow[""{name=3, anchor=center, inner sep=0}, "L"', shift right, from=2-2, to=2-1]
			\arrow["\dashv"{anchor=center, rotate=-90}, draw=none, from=2, to=0]
			\arrow["\dashv"{anchor=center, rotate=-90}, draw=none, from=3, to=1]
		\end{tikzcd}\]
		where $L \dashv R$ and $L' \dashv R'$ are adjunctions. We will moreover denote by $\eta$ (resp.\ $\eta'$) the unit $Id_{\mcal D} \implies L \circ R$ of the adjunction $L \dashv R$ (resp.\ $L' \dashv R'$), and $\epsilon$ (resp.\ $\epsilon'$) its co-unit.
		\begin{itemize}
			\item Natural transformations of the form $V \circ R' \implies R \circ U$ (resp.\ $L \circ V \implies U \circ L'$) will be called \textup{right transformations} (resp.\ \textup{left transformations}).
			\item Given a right transformation $\phi : V \circ R' \implies R \circ U$, we can construct the following left transformation: $\phi^L := L \circ V \implies L \circ V \circ R' \circ L' \implies L \circ R \circ U \circ L' \implies U \circ L'$, where the first map is induced by the unit of $L' \dashv R'$, the second is induced by $\phi$, and the last is induced by the co-unit of $L \dashv R$.
			More precisely, $\phi^L := (\epsilon \star (U \circ L')) \circ (L \star \phi \star L') \circ ((L \circ V) \star \eta')$.
			\textup{We will call $\phi^L$ the left adjoint of $\phi$.}
			\item Conversely, given a left transformation $\psi : L \circ V \implies U \circ L'$, we can construct a right transformation $\psi^R : V \circ R' \implies R \circ U$ using this time the unit of $L \dashv R$, $\psi$, and the co-unit of $L' \dashv R'$.
			More precisely, $\psi^R := ((R \circ U) \star \epsilon') \circ (R \star \psi \star R') \circ (\eta \star (V \circ R'))$.
			\textup{We will call $\psi^R$ the right adjoint of $\psi$.}
			\item The operations $-^L$ and $-^R$ are inverse to each other (by the triangle identities of $\eta,\epsilon$ and $\eta',\epsilon'$) and promote to an equivalence between the spaces of right and left transformations: $-^L : Nat(V \circ R',R \circ U) \simeq Nat(L \circ V,U \circ L') : -^R$.
			\item Starting from the square formed by $U,V$ and $R,R'$ only (resp.\ $L,L'$ only) equipped with a given right transformation (resp.\ left transformation), the data of the square formed by $U,V$ and $L,L'$ (resp.\ $R,R'$) equipped with the left adjoint (resp.\ right adjoint) of the transformation is called the \textup{transposed diagram}.
		\end{itemize}
	\end{reminder}
	Note that neither $-^L$ nor $-^R$ necessarily preserve natural equivalences in general, so the transpose of a transformation which was an equivalence might not still be an equivalence. This motivates the definition of the Beck-Chevalley condition (as in \cite[Rmk.~4.7.4.15]{HA}).
	\begin{reminder} \label[reminder]{def_beck_chevalley}
		In the context of \cref{def_pre_square}, given a left or a right transformation, we say that it satisfies the Beck-Chevalley condition if it is an equivalence and its adjoint also is an equivalence.
		The data of a square equipped with such a transformation will be called a Beck-Chevalley square.
	\end{reminder}
	In this subsection we will use in particular the following properties of Beck-Chevalley squares also present in \cite{HA} (resp.\ Def. 4.7.4.16 and Cor. 4.7.4.18).
	\begin{prop} \label[prop]{prop_beck_chevalley_compose}
		There is a canonical $\infty$-category $\mcal{BC}$ such that:
		\begin{itemize}
			\item Its objects are quadruples $(\mcal C, \mcal D, R, L)$ where $\mcal C$ and $\mcal D$ are $\infty$-categories, and $R : \mcal C \to \mcal D$ and $L : \mcal D \to \mcal C$ form an adjunction $L \dashv R$.
			\item Given two objects $(\mcal C', \mcal D', R', L')$ and $(\mcal C, \mcal D, R, L)$, a 1-morphism from the first to the second is the data of Beck-Chevalley square:
			\[\begin{tikzcd}
				{\mcal C'} & {\mcal D'} \\
				{\mcal C} & {\mcal D}
				\arrow["U"', from=1-1, to=2-1]
				\arrow["V", from=1-2, to=2-2]
				\arrow[""{name=0, anchor=center, inner sep=0}, "{R'}"', shift right, from=1-1, to=1-2]
				\arrow[""{name=1, anchor=center, inner sep=0}, "R"', shift right, from=2-1, to=2-2]
				\arrow[""{name=2, anchor=center, inner sep=0}, "{L'}"', shift right, from=1-2, to=1-1]
				\arrow[""{name=3, anchor=center, inner sep=0}, "L"', shift right, from=2-2, to=2-1]
				\arrow["\dashv"{anchor=center, rotate=-90}, draw=none, from=2, to=0]
				\arrow["\dashv"{anchor=center, rotate=-90}, draw=none, from=3, to=1]
			\end{tikzcd}\]
		\end{itemize}
		In particular, note that the fact it is an $\infty$-category implies that such squares can be composed. More precisely, consider two consecutive Beck-Chevalley squares:
		\[\begin{tikzcd}
			{\mcal C''} & {\mcal D''} \\
			{\mcal C'} & {\mcal D'} \\
			{\mcal C} & {\mcal D}
			\arrow["{U'}"', from=1-1, to=2-1]
			\arrow["U"', from=2-1, to=3-1]
			\arrow["V", from=2-2, to=3-2]
			\arrow["{V'}", from=1-2, to=2-2]
			\arrow[""{name=0, anchor=center, inner sep=0}, "R"', shift right, from=3-1, to=3-2]
			\arrow[""{name=1, anchor=center, inner sep=0}, "L"', shift right, from=3-2, to=3-1]
			\arrow[""{name=2, anchor=center, inner sep=0}, "{R'}"', shift right, from=2-1, to=2-2]
			\arrow[""{name=3, anchor=center, inner sep=0}, "{L'}"', shift right, from=2-2, to=2-1]
			\arrow[""{name=4, anchor=center, inner sep=0}, "{R''}"', shift right, from=1-1, to=1-2]
			\arrow[""{name=5, anchor=center, inner sep=0}, "{L''}"', shift right, from=1-2, to=1-1]
			\arrow["\dashv"{anchor=center, rotate=-90}, draw=none, from=5, to=4]
			\arrow["\dashv"{anchor=center, rotate=-90}, draw=none, from=3, to=2]
			\arrow["\dashv"{anchor=center, rotate=-90}, draw=none, from=1, to=0]
		\end{tikzcd}\]
		where we call $\phi$ the Beck-Chevalley right transformation for the bottom square, and $\psi$ its adjoint (resp.\ $\phi'$ and $\psi'$ for the top square).
		Then the exterior pasting square also is Beck-Chevalley when equipped with the right transformation $\Phi := (\phi \star U') \circ (V \star \phi') : V \circ V' \circ R'' \implies V \circ R' \circ U' \implies R \circ U \circ U'$, and its adjoint is $\Psi := (\psi \star V') \circ (U \star \psi')$.
	\end{prop}
	\begin{proof}
		It suffices to apply the definition 4.7.4.16 from \cite{HA} to the simplicial set $S = \Delta^1$ to produce the required $\infty$-category.
	\end{proof}
	\begin{prop} \label[prop]{prop_beck_chevalley_limit}
		Let $(\mcal C_i, \mcal D_i, L_i, R_i)_{i \in I}$ be a small diagram in the $\infty$-category $\mcal{BC}$ (i.e.\ each $L_i \dashv R_i$ forms an adjunction between $\mcal C_i$ and $\mcal D_i$, and they are related by Beck-Chevalley squares indexed by the edges of $I$).
		Then this diagram admits a limit which is moreover computed component-wise on the objects.
		In particular, for each vertex $i \in I$, there is a canonical Beck-Chevalley square:
		\[\begin{tikzcd}[cramped]
			{\textup{lim}_i \mcal C_i} & {\textup{lim}_i \mcal D_i} \\
			{\mcal C_i} & {\mcal D_i}
			\arrow[""{name=0, anchor=center, inner sep=0}, "{\textup{lim}_i R_i}"', shift right, from=1-1, to=1-2]
			\arrow[""{name=1, anchor=center, inner sep=0}, "{\textup{lim}_i L_i}"', shift right, from=1-2, to=1-1]
			\arrow["{U_i}"', from=1-1, to=2-1]
			\arrow["{V_i}", from=1-2, to=2-2]
			\arrow[""{name=2, anchor=center, inner sep=0}, "{R_i}"', shift right, from=2-1, to=2-2]
			\arrow[""{name=3, anchor=center, inner sep=0}, "{L_i}"', shift right, from=2-2, to=2-1]
			\arrow["\bot"{description}, draw=none, from=1, to=0]
			\arrow["\bot"{description}, draw=none, from=3, to=2]
		\end{tikzcd}\]
		where $U_i$ and $V_i$ are the canonical projection maps from the respective limits of $\mcal C_i$ and $\mcal D_i$.
	\end{prop}
	\begin{proof}
		It suffices to apply the corollary 4.7.4.18 from \cite{HA} to the simplicial set $S = \Delta^1$.
	\end{proof}
	Moreover, we will also need the following lemma which we will prove.
	\begin{lem} \label[lem]{bc_unit}
		In the context of \cref{def_beck_chevalley}, let $\phi$ be a right transformation satisfying the Beck-Chevalley condition, and let $\psi$ be its adjoint.
		Then, the following diagram is commutative:
		\[\begin{tikzcd}[cramped]
			{U \circ L' \circ R'} & {L \circ V \circ R'} & {L \circ R \circ U} \\
			& U
			\arrow["\simeq", from=1-2, to=1-1]
			\arrow["\simeq"', from=1-2, to=1-3]
			\arrow["{U \star \epsilon'}"', from=1-1, to=2-2]
			\arrow["{\epsilon \star U}", from=1-3, to=2-2]
			\arrow["{\psi \star R'}"', from=1-2, to=1-1]
			\arrow["{L \star \phi}", from=1-2, to=1-3]
		\end{tikzcd}\]
		In other words, $U$ sends the co-unit of $L' \dashv R'$ to the co-unit of $L \dashv R$ applied to $U$.
		Similarly, $V$ sends the unit of $L' \dashv R'$ to the unit of $L \dashv R$ applied to $V$.
	\end{lem}
	\begin{proof}
		We will prove the statement for $U$ since the proof for $V$ is completely analogous.
		More precisely, it means that we will prove that $(U \star \epsilon') \circ (\psi \star R') = (\epsilon \star U) \circ (L \star \phi)$ as natural transformations $L \circ V \circ R' \implies U$.
		Since by assumption $\psi$ is the left adjoint to $\phi$, by definition we have $\psi := (\epsilon \star (U \circ L')) \circ (L \star \phi \star L') \circ ((L \circ V) \star \eta')$.
		In particular, using the associativity and distributivity laws of the vertical and horizontal compositions of natural transformations, and the triangle identity $(R' \star \epsilon') \circ (\eta' \star R') = Id_{R'}$, we can compute as required:
		\begin{align*}
			&(U \star \epsilon') \circ (\psi \star R') \\
			=\ &(U \star \epsilon') \circ (((\epsilon \star (U \circ L')) \circ (L \star \phi \star L') \circ ((L \circ V) \star \eta')) \star R') \\
			\simeq\ &(U \star \epsilon') \circ (\epsilon \star (U \circ L' \circ R')) \circ (L \star \phi \star (L' \circ R')) \circ ((L \circ V) \star \eta' \star R') \\
			\simeq\ &(\epsilon \star U) \circ (L \star \phi) \circ ((L \circ V \circ R') \star \epsilon') \circ ((L \circ V) \star \eta' \star R') \\
			\simeq\ &(\epsilon \star U) \circ (L \star \phi) \circ ((L \circ V) \star ((R' \star \epsilon') \circ (\eta' \star R'))) \\
			\simeq\ &(\epsilon \star U) \circ (L \star \phi) \circ ((L \circ V) \star Id_{R'}) \\
			\simeq\ &(\epsilon \star U) \circ (L \star \phi) \circ (Id_{L \circ V \circ R'}) \\
			\simeq\ &(\epsilon \star U) \circ (L \star \phi)
		\end{align*}
		The reader familiar with string diagrams computations is invited to use them to follow along.
	\end{proof}
	\subsection{Beck-Chevalley Squares in Derived Algebraic Geometry} \label[subsection]{beck-chevalley-dag}
	In this subsection, we will look at instances of Beck-Chevalley conditions in the context of derived algebraic geometry. We start from a general adjunction between quasi-coherent algebras and relatively affine derived stacks, and arrive at two conditions, one similar to the base-change formula (in \cref{df_f_right}) and the other similar to relative affineness (in \cref{df_f_left}). We will also find sufficient conditions on maps for these Beck-Chevalley conditions to hold, which will apply in the context of \cref{equiv-affine-case}.
	\begin{reminder} \label[reminder]{rmk_f_beck_chevalley}
		Let $f : V \to U$ be an arbitrary map between derived stacks, and consider the following diagram, as described in \cref{qcoh_dst_adj} and \cref{base_chg_adj}, in which a priori no commutativity holds, which we emphasize by inserting a question mark in its center:
		\[\begin{tikzcd}[cramped]
			{\CAlg(\QCoh(V))^\op} && {\dSt/V} \\
			\\
			{\CAlg(\QCoh(U))^\op} && {\dSt/U}
			\arrow[""{name=0, anchor=center, inner sep=0}, "{\Spec_V}"', shift right, from=1-1, to=1-3]
			\arrow[""{name=1, anchor=center, inner sep=0}, "{\mcal O_{/V}}"', shift right, from=1-3, to=1-1]
			\arrow[""{name=2, anchor=center, inner sep=0}, "{\Spec_U}"', shift right, from=3-1, to=3-3]
			\arrow[""{name=3, anchor=center, inner sep=0}, "{\mcal O_{/U}}"', shift right, from=3-3, to=3-1]
			\arrow[""{name=4, anchor=center, inner sep=0}, "{f_!}"', shift right, from=1-3, to=3-3]
			\arrow[""{name=5, anchor=center, inner sep=0}, "{f^*}"', shift right, from=3-3, to=1-3]
			\arrow[""{name=6, anchor=center, inner sep=0}, "{f_*}"', shift right, from=1-1, to=3-1]
			\arrow[""{name=7, anchor=center, inner sep=0}, "{f^*}"', shift right, from=3-1, to=1-1]
			\arrow["{?}"{description}, draw=none, from=3-1, to=1-3]
			\arrow["\dashv"{anchor=center}, draw=none, from=4, to=5]
			\arrow["\dashv"{anchor=center}, draw=none, from=6, to=7]
			\arrow["\dashv"{anchor=center, rotate=-90}, draw=none, from=1, to=0]
			\arrow["\dashv"{anchor=center, rotate=-90}, draw=none, from=3, to=2]
		\end{tikzcd}\]
		As stated in \cref{base_chg_adj}, in the above situation, the following diagrams are always commutative:
		\[\begin{tikzcd}[cramped]
			{\CAlg(\QCoh(V))^\op} & {\dSt/V} && {\CAlg(\QCoh(V))^\op} & {\dSt/V} \\
			{\CAlg(\QCoh(U))^\op} & {\dSt/U} && {\CAlg(\QCoh(U))^\op} & {\dSt/U}
			\arrow["{\Spec_V}"', from=1-4, to=1-5]
			\arrow["{\Spec_U}"', from=2-4, to=2-5]
			\arrow["{f^*}"', from=2-5, to=1-5]
			\arrow["{f^*}"', shift left, from=2-4, to=1-4]
			\arrow["{f_!}", from=1-2, to=2-2]
			\arrow["{f_*}"', from=1-1, to=2-1]
			\arrow["{\mcal O_{/U}}"', from=2-2, to=2-1]
			\arrow["{\mcal O_{/V}}"', from=1-2, to=1-1]
		\end{tikzcd}\]
	\end{reminder}
	\begin{rmk}
		In particular, we can wonder when the left square is right adjointable (resp.\ the right square is left adjointable), i.e.\ whether or not the two following squares, which by construction come equipped with canonical natural transformations obtained as adjoints of the above commutative squares, also commute:
		\[\begin{tikzcd}[cramped]
			{\CAlg(\QCoh(V))^\op} & {\dSt/V} && {\CAlg(\QCoh(V))^\op} & {\dSt/V} \\
			{\CAlg(\QCoh(U))^\op} & {\dSt/U} && {\CAlg(\QCoh(U))^\op} & {\dSt/U}
			\arrow["{f^*}"', from=2-5, to=1-5]
			\arrow["{f^*}"', shift left, from=2-4, to=1-4]
			\arrow["{f_!}", from=1-2, to=2-2]
			\arrow["{f_*}"', from=1-1, to=2-1]
			\arrow["{\Spec_V}"', from=1-1, to=1-2]
			\arrow["{\Spec_U}"', from=2-1, to=2-2]
			\arrow["{\mcal O_{/V}}"', from=1-5, to=1-4]
			\arrow["{\mcal O_{/U}}"', from=2-5, to=2-4]
			\arrow["{?}"{description}, draw=none, from=2-1, to=1-2]
			\arrow["{?}"{description}, draw=none, from=2-4, to=1-5]
		\end{tikzcd}\]
	\end{rmk}
	We will therefore consider two definitions.
	\begin{df} \label[df]{df_f_left}
		We say that a map of derived stacks $f : V \to U$ \textup{preserves relative spectra by push-forward} if the canonical natural transform $f_! \circ \Spec_V \to \Spec_U \circ f_*$, induced by \cref{rmk_f_beck_chevalley}, is an equivalence (i.e.\ if the left square of this remark is Beck-Chevalley).
		More explicitly, for any algebra $A \in \CAlg(\QCoh(V))$, there is always a canonical commutative diagram:
		\[\begin{tikzcd}
			{\Spec_V(A)} & {\Spec_U(f_* A)} \\
			V & U
			\arrow["{\phi_A}", from=1-1, to=1-2]
			\arrow[from=1-1, to=2-1]
			\arrow[from=1-2, to=2-2]
			\arrow["f"', from=2-1, to=2-2]
		\end{tikzcd}\]
		and we require that the top map $\phi_A$ is an equivalence for all $A$. In particular, when $A = \mcal O_V$, we obtain that $V \simeq \Spec_U(f_* \mcal O_V)$ as derived stacks over $U$, i.e. $V \to U$ is relatively affine.
	\end{df}
	\begin{df} \label[df]{df_f_right}
		We say that a map of derived stacks $f : V \to U$ \textup{preserves relative structural sheaves by pull-back} if the canonical natural transform $f^* \circ \mcal O_{/U} \to \mcal O_{/V} \circ f^*$, induced by \cref{rmk_f_beck_chevalley}, is an equivalence (i.e.\ if the right square of this remark is Beck-Chevalley).
		More explicitly, for any $q : T \to U$, there is always a canonical pull-back diagram:
		\[\begin{tikzcd}
			{f^*T} & T \\
			V & U
			\arrow["{f'}", from=1-1, to=1-2]
			\arrow["p"', from=1-1, to=2-1]
			\arrow["\lrcorner"{anchor=center, pos=0.125}, draw=none, from=1-1, to=2-2]
			\arrow["q", from=1-2, to=2-2]
			\arrow["f"', from=2-1, to=2-2]
		\end{tikzcd}\]
		and a canonical map $f^* q_* \mcal O_T \to p_* f'^* \mcal O_T$ of algebras in $\QCoh(V)$, which we require to be an equivalence for all $T$. This is related to a base-change formula.
	\end{df}
	We will now provide conditions on a map $f$ that are sufficient to satisfy these properties. While the first condition can hold under reasonable conditions, the second one will require a refinement, as it is too strict when formulated as above. Let us thus start with the first one.
	\begin{lem} \label[lem]{affine_pushforward_spec}
		If $f : V \to U$ is a quotient of a map between derived stacks relatively affine over $\Spec(k)$\footnote{As noted in \cref{equiv_eps_gr_BH}, we mean here that $U$ and $V$ are affine stacks in the sense of \cite[\S2.1]{AffineStacks}, so we mean that $U = \Spec_k(B)$ and $V = \Spec_k(C)$ for $B$ and $C$ two possibly non-connective $k$-cdga's. Equivalently, we mean that the terminal maps $U,V \to \Spec(k)$ are relatively affine in the usual sense.} by an affine group scheme, then it preserves relative spectra by push-forward (see \cref{df_f_left}).
	\end{lem}
	\begin{proof}
		Let us first only assume that $f : V \to U$ is a map between derived stacks relatively affine over $\Spec(k)$. We will later generalize it to the case of a quotient by an affine group scheme.
		We will now assume that $U = \Spec_k(B)$ and $V = \Spec_k(C)$ for some possibly non-connective $k$-cdga's $B$ and $C$, that $f = \Spec_k(\phi)$ for some $\phi : B \to C$ a map of $k$-cdgas. In this context, let $A$ be a $C$-algebra.
		We have to show that the canonical map $f_! \Spec_C(A) \to \Spec_B(f_* A)$ is an equivalence of derived stacks over $\Spec(B)$, i.e.\ that the top map of the following canonical commutative square is an equivalence:
		\[\begin{tikzcd}
			{\Spec_C(A)} & {\Spec_B(f_* A)} \\
			{\Spec(C)} & {\Spec(B)}
			\arrow["{\simeq ?}", from=1-1, to=1-2]
			\arrow[from=1-1, to=2-1]
			\arrow[from=1-2, to=2-2]
			\arrow[from=2-1, to=2-2]
		\end{tikzcd}\]
		In this case, we can give an explicit description of these stacks and this map.
		Let $\Spec(O)$ be an affine, by definition we have:
		\begin{itemize}
			\item $\Map_{\dSt}(\Spec(O),\Spec(B)) = \Map_{cdga}(B,O)$
			\item $\Map_{\dSt}(\Spec(O),\Spec(C)) = \Map_{cdga}(C,O)$
			\item Given a structural map $\Spec(O) \to \Spec(C) \in \dSt/\Spec(C)$, we have:\\
			$\Map_{\dSt/\Spec(C)}(\Spec(O),\Spec_C(A)) = \Map_{C-Alg}(A,O)$
			\item Given a structural map $\Spec(O) \to \Spec(B) \in \dSt/\Spec(B)$, we have:\\
			$\Map_{\dSt/\Spec(B)}(\Spec(O),\Spec_B(f_* A)) = \Map_{B-Alg}(f_* A,O)$
		\end{itemize}
		In the last two items, some data is in fact redundant. Let us look at $\Spec_C(A)$ for example, in which case the required data is two dashed arrows forming a commutative triangle (the plain arrow is fixed, it is the structural map $C \to A$):
		\[\begin{tikzcd}
			O & A \\
			& C
			\arrow[dashed, from=1-2, to=1-1]
			\arrow[dashed, from=2-2, to=1-1]
			\arrow[from=2-2, to=1-2]
		\end{tikzcd}\]
		This data is equivalent to the data of a single map $A \to O$, since by commutativity the induced map $C \to O$ has to be the composition $C \to A \to O$.
		Moreover, in the case of $\Spec_B(f_* A)$, note that $f_* A$ is equivalent to $A$ as a cdga, as in this context the notation $f_* A$ only serves to indicate that we equip it with the fixed structural map $B \to A$ given by the composition $B \to C \to A$.
		Therefore, in this context where everything is affine, $f_! \Spec_C(A)$ and $\Spec_B(f_* A)$ are equivalent to $\Spec(A)$ as derived stacks, and under this equivalence the canonical map $f_! \Spec_C(A) \to \Spec_B(f_* A)$ is the identity, i.e.\ the following square of derived stacks commutes:
		\[\begin{tikzcd}
			{\Spec(A)} & {\Spec(A)} \\
			{f_!\Spec_C(A)} & {\Spec_B(f_* A)}
			\arrow[from=2-1, to=2-2]
			\arrow[Rightarrow, no head, from=1-1, to=1-2]
			\arrow["\simeq", from=2-1, to=1-1]
			\arrow["\simeq"', from=2-2, to=1-2]
		\end{tikzcd}\]
		This shows that the bottom map is also an equivalence as required.
		
		Let us assume now that $f : V \to U$ is not a map of derived stacks relatively affine over $\Spec(k)$, but a quotient of such a map $g : V' \to U'$ by an affine group scheme $G$ (i.e.\ $U = [U'/G]$, $V = [V'/G]$, $f = [g/G]$, and $U',V'$ and $G$ are derived stacks relatively affine over $\Spec(k)$).
		We are thus assuming that we have the following colimit diagram, where each of the two columns is a simplicial object:
		\[\begin{tikzcd}
			{[V'/G] = V} & {U = [U'/G]} \\
			{V'} & {U'} \\
			{V' \times G} & {U' \times G} \\
			{V' \times G^2} & {U' \times G^2} \\
			\cdots & \cdots
			\arrow["f", from=1-1, to=1-2]
			\arrow[from=2-1, to=1-1]
			\arrow["{g_0 = g}"', from=2-1, to=2-2]
			\arrow[from=2-2, to=1-2]
			\arrow[from=3-1, to=2-1]
			\arrow["{g_1}"', from=3-1, to=3-2]
			\arrow[from=3-2, to=2-2]
			\arrow[from=4-1, to=3-1]
			\arrow["{g_2}"', from=4-1, to=4-2]
			\arrow[from=4-2, to=3-2]
			\arrow[from=5-1, to=4-1]
			\arrow[from=5-2, to=4-2]
		\end{tikzcd}\]
		Recall from \cref{prop_beck_chevalley_compose} that the $\infty$-category $\mcal{BC}$ of adjunctions as objects and Beck-Chevalley squares as morphisms admits limits that are computed pointwise, so that its arrow $\infty$-category also has similar limits.
		For $n \in \bb N$, we will denote by $V_n$ (resp.\ $U_n$) the following adjunctions, which are seen as objects of $\mcal{BC}$:
		\[\begin{tikzcd}
			& {V_n :} \\
			{\CAlg(\QCoh(V' \times G^n))^\op} && {\dSt/V' \times G^n} \\
			\\
			& {U_n :} \\
			{\CAlg(\QCoh(U' \times G^n))^\op} && {\dSt/U' \times G^n}
			\arrow["{\Spec_{V' \times G^n}}"', shift right, from=2-1, to=2-3]
			\arrow["{\mcal O_{/V' \times G^n}}"', shift right, from=2-3, to=2-1]
			\arrow["{\Spec_{U' \times G^n}}"', shift right, from=5-1, to=5-3]
			\arrow["{\mcal O_{/U' \times G^n}}"', shift right, from=5-3, to=5-1]
		\end{tikzcd}\]
		Now at each level $n \in \bb N$, since $G$ and the map $g$ are relatively affine, $g_n$ is relatively affine too, and the following square is Beck-Chevalley by the preceding case of this lemma:
		\[\begin{tikzcd}
			{\CAlg(\QCoh(V' \times G^n))^\op} & {\dSt/V' \times G^n} \\
			{\CAlg(\QCoh(U' \times G^n))^\op} & {\dSt/U' \times G^n}
			\arrow["{\Spec_{V' \times G^n}}"', shift right, from=1-1, to=1-2]
			\arrow["{(g_n)_*}"', from=1-1, to=2-1]
			\arrow["{\mcal O_{/V' \times G^n}}"', shift right, from=1-2, to=1-1]
			\arrow["{(g_n)_!}", from=1-2, to=2-2]
			\arrow["{\Spec_{U' \times G^n}}"', shift right, from=2-1, to=2-2]
			\arrow["{\mcal O_{/U' \times G^n}}"', shift right, from=2-2, to=2-1]
		\end{tikzcd}\]
		We will call $C_n$ this Beck-Chevalley square, thus $C_n : V_n \to U_n$ is an arrow in the $\infty$-category of $\mcal{BC}$, or equivalently it is an object in the arrow $\infty$-category of $\mcal{BC}$.
		Moreover, we will now describe transition maps from $C_n$ to $C_{n+1}$.
		Note that the following square is Beck-Chevalley too by the same argument since both $V' \times G^{n+1}$ and $V' \times G^n$ are relatively affine over $\Spec(k)$ (the same is true in the analogous case for $U'$):
		\[\begin{tikzcd}
			{\CAlg(\QCoh(V' \times G^{n+1}))^\op} & {\dSt/V' \times G^{n+1}} \\
			{\CAlg(\QCoh(V' \times G^n))^\op} & {\dSt/V' \times G^n}
			\arrow["{\Spec_{V' \times G^{n+1}}}"', shift right, from=1-1, to=1-2]
			\arrow[from=1-1, to=2-1]
			\arrow["{\mcal O_{/V' \times G^{n+1}}}"', shift right, from=1-2, to=1-1]
			\arrow[from=1-2, to=2-2]
			\arrow["{\Spec_{V' \times G^n}}"', shift right, from=2-1, to=2-2]
			\arrow["{\mcal O_{/V' \times G^n}}"', shift right, from=2-2, to=2-1]
		\end{tikzcd}\]
		This provides arrows $V_n \to V_{n+1}$ and $U_n \to U_{n+1}$, and moreover the following diagram commutes:
		\[\begin{tikzcd}
			{V_n} & {U_n} \\
			{V_{n+1}} & {U_{n+1}}
			\arrow["{C_n}", from=1-1, to=1-2]
			\arrow[from=1-1, to=2-1]
			\arrow[from=1-2, to=2-2]
			\arrow["{C_{n+1}}"', from=2-1, to=2-2]
		\end{tikzcd}\]
		Therefore, we have the following co-simplicial object in the $\infty$-category of Beck-Chevalley squares (the arrow $\infty$-category of $\mcal{BC}$):
		\[\begin{tikzcd}[cramped]
			{C_0} \\
			{C_1} \\
			{C_2} \\
			\vdots
			\arrow[from=1-1, to=2-1]
			\arrow[from=2-1, to=3-1]
			\arrow[from=3-1, to=4-1]
		\end{tikzcd}\]
		By \cref{prop_beck_chevalley_limit}, it admits a limit $C$, computed pointwise, which is also a Beck-Chevalley square. Moreover, by definition and using the pointwise limit, we can immediately see that $C$ is the following square:
		\[\begin{tikzcd}
			{\CAlg(\QCoh(V)^\op} & {\dSt/V} \\
			{\CAlg(\QCoh(U)^\op} & {\dSt/U}
			\arrow["{\Spec_{V}}"', shift right, from=1-1, to=1-2]
			\arrow["{f_*}"', from=1-1, to=2-1]
			\arrow["{\mcal O_{/V}}"', shift right, from=1-2, to=1-1]
			\arrow["{f_!}", from=1-2, to=2-2]
			\arrow["{\Spec_{U}}"', shift right, from=2-1, to=2-2]
			\arrow["{\mcal O_{/U}}"', shift right, from=2-2, to=2-1]
		\end{tikzcd}\]
		The fact that this square is Beck-Chevalley concludes this lemma.
	\end{proof}
	As explained earlier, we will now turn to the second condition, but we will have to refine it.
	\begin{rmk}
		Let $f : V \to U$ be a map of derived stacks, and suppose we want to prove that it preserves relative structural sheaves by pull-back (see \cref{df_f_right}).
		Let $q : T \to U$ be any map of derived stacks, we can form the following canonical pull-back square:
		\[\begin{tikzcd}
			{f^* T} & T \\
			V & U
			\arrow["{f'}", from=1-1, to=1-2]
			\arrow["{q'}"', from=1-1, to=2-1]
			\arrow["\lrcorner"{anchor=center, pos=0.125}, draw=none, from=1-1, to=2-2]
			\arrow["q", from=1-2, to=2-2]
			\arrow["f"', from=2-1, to=2-2]
		\end{tikzcd}\]
		We have to show that the following canonical map of commutative algebras in $\QCoh(V)$ is an equivalence:
		\[f^* p'_* \mcal O_T \to p_* f'^* \mcal O_T\]
		This statement is a base-change formula, and we are requiring it to be true along any arbitrary map $q : T \to U$, which need not be true under reasonable assumptions in this level of generality.
		However, as is well-known in algebraic geometry, we will show that it is true when the statement is restricted to maps $q : T \to U$ that are relatively affine, and $f : V \to U$ is finite.
	\end{rmk}
	\begin{df}
		Given a derived stack $S$, the essential image of $\Spec_S$ will be denoted by:
		\[\textup{RelAff}/S\]
	\end{df}
	\begin{rmk}
		As in \cref{rmk_f_beck_chevalley}, given $f : U \to V$ any map of derived stacks, we can construct the following square:
		\[\begin{tikzcd}
			{\CAlg(\QCoh(U))^\op} & {\dSt/U} \\
			{\CAlg(\QCoh(V))^\op} & {\dSt/V}
			\arrow[""{name=0, anchor=center, inner sep=0}, "{\Spec_U}"', shift right, from=1-1, to=1-2]
			\arrow["{f^*}", from=1-1, to=2-1]
			\arrow[""{name=1, anchor=center, inner sep=0}, "{\mcal O_{/U}}"', shift right, from=1-2, to=1-1]
			\arrow["{f^*}"', from=1-2, to=2-2]
			\arrow[""{name=2, anchor=center, inner sep=0}, "{\Spec_V}"', shift right, from=2-1, to=2-2]
			\arrow[""{name=3, anchor=center, inner sep=0}, "{\mcal O_{/V}}"', shift right, from=2-2, to=2-1]
			\arrow["\dashv"{anchor=center, rotate=-90}, draw=none, from=1, to=0]
			\arrow["\dashv"{anchor=center, rotate=-90}, draw=none, from=3, to=2]
		\end{tikzcd}\]
		Moreover, the following diagram commutes:
		\[\begin{tikzcd}
			{\CAlg(\QCoh(U))^\op} & {\dSt/U} \\
			{\CAlg(\QCoh(V))^\op} & {\dSt/V}
			\arrow["{\Spec_U}"', from=1-1, to=1-2]
			\arrow["{f^*}", from=1-1, to=2-1]
			\arrow["{f^*}"', from=1-2, to=2-2]
			\arrow["{\Spec_V}"', from=2-1, to=2-2]
		\end{tikzcd}\]
		As a consequence, given $q : T \to U$, if $T$ lies in $\textup{RelAff}/U$, i.e.\ if there exists $A \in \CAlg(\QCoh(U))$ such that $T \simeq \Spec_U(A)$ as derived stacks over $U$, then $f^* T \simeq \Spec_V(f^* A)$ as derived stacks over $V$, i.e.\ $f^* T$ lies in $\textup{RelAff}/V$.
		In conclusion, the previous square restricts to the following one:
		\[\begin{tikzcd}
			{\CAlg(\QCoh(U))^\op} & {\textup{RelAff}/U} \\
			{\CAlg(\QCoh(V))^\op} & {\textup{RelAff}/V}
			\arrow[""{name=0, anchor=center, inner sep=0}, "{\Spec_U}"', shift right, from=1-1, to=1-2]
			\arrow["{f^*}", from=1-1, to=2-1]
			\arrow[""{name=1, anchor=center, inner sep=0}, "{\mcal O_{/U}}"', shift right, from=1-2, to=1-1]
			\arrow["{f^*}"', from=1-2, to=2-2]
			\arrow[""{name=2, anchor=center, inner sep=0}, "{\Spec_V}"', shift right, from=2-1, to=2-2]
			\arrow[""{name=3, anchor=center, inner sep=0}, "{\mcal O_{/V}}"', shift right, from=2-2, to=2-1]
			\arrow["\dashv"{anchor=center, rotate=-90}, draw=none, from=1, to=0]
			\arrow["\dashv"{anchor=center, rotate=-90}, draw=none, from=3, to=2]
		\end{tikzcd}\]
	\end{rmk}
	\begin{df} \label[df]{df_f_left_relaff}
		We will say that a map of derived stacks $f : V \to U$ preserves the relative structural sheaves of relative affines by push-forward if the above square is Beck-Chevalley, i.e.\ if when restricted to $\textup{RelAff}/U$, the canonical natural transform:
		\[f^* \circ \mcal O_U \to \mcal O_V \circ f^*\]
		is an equivalence.
	\end{df}
	We will now give a reasonable condition on $f$ for this weaker condition to hold.
	\begin{lem} \label[lem]{f_finite}
		If $f : V \to U$ is a finite map of derived stacks, it preserves the relative structural sheaves of relative affines by push-forward.
	\end{lem}
	\begin{proof}
		Let $q : T \to U$ be a relatively affine map of derived stacks, we can form the following canonical pull-back square:
		\[\begin{tikzcd}
			{f^* T} & T \\
			V & U
			\arrow["{f'}", from=1-1, to=1-2]
			\arrow["{q'}"', from=1-1, to=2-1]
			\arrow["\lrcorner"{anchor=center, pos=0.125}, draw=none, from=1-1, to=2-2]
			\arrow["q", from=1-2, to=2-2]
			\arrow["f"', from=2-1, to=2-2]
		\end{tikzcd}\]
		We have to show that the following canonical map of commutative algebras in $\QCoh(V)$ is an equivalence:
		\[f^* p'_* \mcal O_T \to p_* f'^* \mcal O_T\]
		However, since $q$ is relatively affine, it is quasi-compact and quasi-separated, and by the assumption that $f : V \to U$ is finite, we can apply \cite[Lem.~A.1.3]{CategoricalProperness}, which yields the result.
	\end{proof}
	\begin{rmk}
		While the first condition provided a Beck-Chevalley square with the $\infty$-categories $\dSt/V$ and $\dSt/U$ appearing on the right-hand side (see \cref{df_f_right}), the previous condition provides a Beck-Chevalley square where the right-hand side is restricted to $\textup{RelAff}/V$ and $\textup{RelAff}/U$ (see \cref{df_f_left_relaff}). As such, these squares cannot be composed as-is. However, we can see that the first condition restricts to relative affines. Indeed, let $f : V \to U$ be a map of derived stacks which preserves relative spectra by push-forward. This means that the following square is Beck-Chevalley:
		\[\begin{tikzcd}
			{\CAlg(\QCoh(V))^\op} & {\dSt/V} \\
			{\CAlg(\QCoh(U))^\op} & {\dSt/U}
			\arrow[""{name=0, anchor=center, inner sep=0}, "{\Spec_V}"', shift right, from=1-1, to=1-2]
			\arrow["{f_*}"', from=1-1, to=2-1]
			\arrow[""{name=1, anchor=center, inner sep=0}, "{\mcal O_{/V}}"', shift right, from=1-2, to=1-1]
			\arrow["{f_!}", from=1-2, to=2-2]
			\arrow[""{name=2, anchor=center, inner sep=0}, "{\Spec_U}"', shift right, from=2-1, to=2-2]
			\arrow[""{name=3, anchor=center, inner sep=0}, "{\mcal O_{/U}}"', shift right, from=2-2, to=2-1]
			\arrow["\dashv"{anchor=center, rotate=-90}, draw=none, from=1, to=0]
			\arrow["\dashv"{anchor=center, rotate=-90}, draw=none, from=3, to=2]
		\end{tikzcd}\]
		In particular, for any $A \in \CAlg(\QCoh(V))$, we have that $f_! \Spec_V(A) \simeq \Spec_U(f_* A)$ as derived stacks over $U$.
		As a consequence, it restricts to the following Beck-Chevalley square:
		\[\begin{tikzcd}
			{\CAlg(\QCoh(V))^\op} & {\textup{RelAff}/V} \\
			{\CAlg(\QCoh(U))^\op} & {\textup{RelAff}/U}
			\arrow[""{name=0, anchor=center, inner sep=0}, "{\Spec_V}"', shift right, from=1-1, to=1-2]
			\arrow["{f_*}"', from=1-1, to=2-1]
			\arrow[""{name=1, anchor=center, inner sep=0}, "{\mcal O_{/V}}"', shift right, from=1-2, to=1-1]
			\arrow["{f_!}", from=1-2, to=2-2]
			\arrow[""{name=2, anchor=center, inner sep=0}, "{\Spec_U}"', shift right, from=2-1, to=2-2]
			\arrow[""{name=3, anchor=center, inner sep=0}, "{\mcal O_{/U}}"', shift right, from=2-2, to=2-1]
			\arrow["\dashv"{anchor=center, rotate=-90}, draw=none, from=1, to=0]
			\arrow["\dashv"{anchor=center, rotate=-90}, draw=none, from=3, to=2]
		\end{tikzcd}\]
	\end{rmk}
	In conclusion of this appendix, we arrived at the following property:
	\begin{conclusion} \label[conclusion]{conclusion-bc}
		Let $f : V \to U$ be a map of derived stacks.
		\begin{enumerate}
			\item By \cref{affine_pushforward_spec}, if $f$ is a quotient of a map between derived stacks relatively affine over $\Spec(k)$ by an affine group scheme, then the following square is Beck-Chevalley:
			\[\begin{tikzcd}
				{\CAlg(\QCoh(V))^\op} & {\textup{RelAff}/V} \\
				{\CAlg(\QCoh(U))^\op} & {\textup{RelAff}/U}
				\arrow[""{name=0, anchor=center, inner sep=0}, "{\Spec_V}"', shift right, from=1-1, to=1-2]
				\arrow["{f_*}"', from=1-1, to=2-1]
				\arrow[""{name=1, anchor=center, inner sep=0}, "{\mcal O_{/V}}"', shift right, from=1-2, to=1-1]
				\arrow["{f_!}", from=1-2, to=2-2]
				\arrow[""{name=2, anchor=center, inner sep=0}, "{\Spec_U}"', shift right, from=2-1, to=2-2]
				\arrow[""{name=3, anchor=center, inner sep=0}, "{\mcal O_{/U}}"', shift right, from=2-2, to=2-1]
				\arrow["\dashv"{anchor=center, rotate=-90}, draw=none, from=1, to=0]
				\arrow["\dashv"{anchor=center, rotate=-90}, draw=none, from=3, to=2]
			\end{tikzcd}\]
			\item By \cref{f_finite}, if $f$ is finite, then the following square is Beck-Chevalley:
			\[\begin{tikzcd}
				{\CAlg(\QCoh(U))^\op} & {\textup{RelAff}/U} \\
				{\CAlg(\QCoh(V))^\op} & {\textup{RelAff}/V}
				\arrow[""{name=0, anchor=center, inner sep=0}, "{\Spec_U}"', shift right, from=1-1, to=1-2]
				\arrow["{f^*}", from=1-1, to=2-1]
				\arrow[""{name=1, anchor=center, inner sep=0}, "{\mcal O_{/U}}"', shift right, from=1-2, to=1-1]
				\arrow["{f^*}"', from=1-2, to=2-2]
				\arrow[""{name=2, anchor=center, inner sep=0}, "{\Spec_V}"', shift right, from=2-1, to=2-2]
				\arrow[""{name=3, anchor=center, inner sep=0}, "{\mcal O_{/V}}"', shift right, from=2-2, to=2-1]
				\arrow["\dashv"{anchor=center, rotate=-90}, draw=none, from=1, to=0]
				\arrow["\dashv"{anchor=center, rotate=-90}, draw=none, from=3, to=2]
			\end{tikzcd}\]
		\end{enumerate}
	\end{conclusion}
	\printbib
\end{document}